\documentclass[11pt,reqno]{amsart}
\usepackage{amssymb}
\usepackage{verbatim}
\usepackage{amsfonts}
\usepackage{mathrsfs}
\usepackage{amsmath}
\usepackage{graphicx}
\usepackage{hyperref}
\usepackage{float}
\usepackage{epstopdf}
\usepackage{color}
\usepackage{bm}
\usepackage{comment}
\usepackage{soul}
\usepackage{cite}
\usepackage{subfigure}

\allowdisplaybreaks
\newtheorem{theorem}{Theorem}[section]
\newtheorem{proposition}[theorem]{Proposition}
\newtheorem{lemma}[theorem]{Lemma}
\newtheorem{corollary}[theorem]{Corollary}
\newtheorem{remark}[theorem]{Remark}
\newtheorem{remarks}[theorem]{Remarks}

\newtheorem{definition}[theorem]{Definition}

\newtheorem{question}[theorem]{Question}

\def\mcP{\mathcal{P}}

\numberwithin{equation}{section}

\begin{document}
	\title[$L^q$-spectra of box-like graph-directed self-affine measures]{$L^q$-spectra of box-like graph-directed self-affine measures: closed forms, with rotation}

	\author{Hua Qiu}
	\address{Department of Mathematics, Nanjing University, Nanjing, 210093, P. R. China.}
	\thanks{The research of Qiu was supported by the National Natural Science Foundation of China, grant 12071213, and the Natural Science Foundation of Jiangsu Province in China, grant BK20211142.}
	\email{huaqiu@nju.edu.cn}
	
	\author{Qi Wang}
	\address{Department of Mathematics, Nanjing University, Nanjing, 210093, P. R. China.}
	\thanks{}
	\email{1378893849@qq.com}

	\subjclass[2010]{Primary 28A80}
	
	\date{}
	
	\keywords{$L^q$-spectra, box dimension, self-affine sets, graph-directed construction, closed form expression. }
	
	\maketitle

	\begin{abstract}
We consider  $L^q$-spectra of planar graph-directed self-affine measures  generated by diagonal or anti-diagonal matrices. Assuming the directed graph is strongly connected and the system satisfies the rectangular open set condition, we obtain a general closed form expression for the $L^q$-spectra. Consequently, we obtain a closed form expression for box dimensions of associated planar graph-directed box-like self-affine sets. We also provide a precise answer to a question of Fraser in 2016 concerning the $L^q$-spectra of planar self-affine measures generated by diagonal matrices.

	\end{abstract}

 \tableofcontents
	
	\section{Introduction}\label{sec1}
    Let $\{T_1, \cdots, T_N\}$ be a finite collection of affine contracting non-singular  matrices, and let $\Psi=\{\psi_i(\cdot)=T_i(\cdot)+t_i\}_{i=1}^N$ be an \textit{iterated function system (IFS)} with $t_i\in \mathbb{R}^n$ for all $i\in \{1,\cdots, N\}$. It is well known that there exists a unique non-empty compact set $X$ such that 
	$$
	X=\bigcup_{i=1}^N \psi_i(X).
	$$
	We call $\Psi$ a \textit{self-affine IFS} and $X$ a \textit{self-affine set}.
	In the special case when $T_i$'s are all similarities, call $\Psi$  a \textit{self-similar IFS} and $X$  a \textit{self-similar set}.  
 
 For a positive probability vector $\mathcal{P}=(p_i)_{i=1}^N$, there exists a probability measure $\mu$ satisfying 
 $$
 \mu= \sum_{i=1}^N p_i \cdot \mu\circ \psi_i^{-1}.
 $$
 Call $\mu$ a \textit{self-affine measure} (resp. \textit{self-similar measure}) when $\Psi$ is a self-affine IFS (resp. self-similar IFS).
	

	The dimension theory of self-affine sets or measures is one of central problems in fractal geometry. Historically, there are two basic strands to determine the Hausdorff  and box dimensions of self-affine sets, one of which is to study  generic self-affine sets basing on the \textit{singular value functions}, and to make \textit{almost sure} statements 
 $$
	\dim_H X=\dim_B X=d(T_1,\cdots, T_N) \text{ for Lebesgue-almost sure } \textbf{t}=(t_1, \cdots, t_N)\in \mathbb{R}^{nN},
	$$
 pioneered by Falconer \cite{F88}. The critical number $d(T_1, \cdots , T_N)$, called \textit{affinity dimension}, is determined in terms of \textit{singular values} of $\{T_1, \cdots, T_N\}$. The original consideration of Falconer requires that all norms of $T_i$'s are less than $1/3$, which was later improved by Solomyak \cite{S98} to $1/2$ and the constant $1/2$ is proved to be sharp in \cite{E92,S98}.
 Along this direction, the study is thriving, see \cite{K04,BHR19,FK18,MS19,M16,BM18,FS14} and the references therein.
	
The other strands of study is to focus on special classes of self-affine sets, and to determine \textit{sure} statements for the dimensions of attractors, which was pioneered by McMullen \cite{M84} and Bedford \cite{B84}, considering planar box-like self-affine sets with homogeneous grid structure.
Their approaches were further developed by  Lalley and Gatzouras \cite{LG92} and Bara\'{n}ski \cite{B07} to box-like sets with certain geometric arrangement or  general grid structure. See \cite{KP96,DS17,K23,FJ21}  for extensions to high dimensions.


 The planar box-like self-affine sets without grid structure were firstly considered by Feng and Wang \cite{FW05}, and later extended by Fraser \cite{F12,F16} allowing the IFS's have non-trivial rotations and reflections (later called \textit{ self-affine carpets}), i.e. linear parts of maps were allowed to be  \textit{diagonal} or \textit{anti-diagonal}. All these works \cite{FW05,F12,F16} on self-affine carpets focus on computing the $L^q$-spectra of their associated self-affine measures. See also \cite{FFL21} for an extension to non-conformal measures.

 In this paper, we continue to study the \textit{$L^q$-spectra} of self-affine measures.  Let $\nu$ be a compactly supported Borel probability measure on $\mathbb{R}^n$ with $n\geq 1$. For $\delta>0$, let $\mathcal{M}_\delta$ be the collection of closed cubes in the $\delta$-mesh of $\mathbb{R}^n$. For $q\geq 0$, write 
	$$ 
	\mathcal{D}_\delta^q(\nu)=\sum_{Q\in \mathcal{M}_\delta}\nu(Q)^q.
	$$
	\begin{definition}
		For $q\geq 0$, the \emph{upper and lower $L^q$-spectra} of $\nu$ are defined to be 
		$$\overline{\tau}_\nu(q)=\limsup_{\delta\to 0+}\frac{\log \mathcal{D}_\delta^q(\nu)}{-\log \delta}
		$$
		and 
		$$\underline{\tau}_\nu(q)=\liminf_{\delta\to 0+}\frac{\log \mathcal{D}_\delta^q(\nu)}{-\log \delta},
		$$
		respectively. If these two values coincide, we define the \emph{$L^q$-spectra} of $\nu$ to be their common value, and denote it as $\tau_\nu(q)$.
	\end{definition}
	
	It is known that as functions of $q$, both $\underline{\tau}_\nu(q),\overline{\tau}_\nu(q)$ are decreasing, and equal to zero at $q=1$. Also,  they are convex, continuous on $(0,\infty)$, and Lipschitz on $[\lambda,\infty)$ for any $\lambda>0$. Note that when $q=0$, the upper and lower $L^q$-spectra are equal to the upper and lower box dimensions of $supp\nu$, respectively. Another important property of $L^q$-spectra is that if it is differential at $q=1$, then the measure $\nu$ is exactly dimensional, and the Hausdorff dimension of $\nu$ equals to $-\tau'(1)$.  The concept of $L^q$-spectra is an important fundamental  ingredient in the study of fractal geometry, particularly in multifractal analysis. See \cite{LN98,FLN00,F05,LW05,F03,K95,O98,F99} and references therein for more details.
 \vspace{0.2cm}

 For a self-similar measure $\mu$ with probability vector $\mathcal{P}=(p_i)_{i=1}^N$, for $q\in \mathbb{R}$, the $L^q$-spectrum of $\mu$ is given by a \textit{closed form}  expression that 
 \begin{equation}\label{t14}
 \sum_{i=1}^N p_i^q r_i^{\tau_{\mu}(q)}=1,
 \end{equation}
 where $r_i$ is the contraction ratio of $\psi_i$. See Cawley and  Mauldin \cite{CM92} and Olsen \cite{O95}. 
 
 For self-affine measures, Feng and Wang \cite[Theorem 2]{FW05} obtained the analogous closed form expression for diagonal self-affine carpets in terms of the $L^q$-spectra of the projections of measures onto $y$-axis providing that the contraction ratios on $x$-axis are less than on $y$-axis for all elements in the IFS's. 
 In a different way, Fraser \cite{F16} introduced the concept of \textit{modified singular value functions} (modified from Falconer's original definition \cite{F88}), and used which to compute the  closed form expression for $L^q$-spectra of self-affine measures on  self-affine carpets without limitation of relative sizes of contraction ratios on $x$-axis or $y$-axis, but still requiring that all $T_i$'s are diagonal. 
 
 In their setting,  all $T_i$'s are of the form $T_i=diag\{\pm a_i, \pm b_i\}$ with $0<a_i,b_i<1$. Let $\gamma_A(q),\gamma_B(q)$ be the unique solutions of 
 $$
	\sum_{i=1}^N p_i^q a_i^{\tau_{\mu^x}(q)} b_i^{\gamma_A(q)-\tau_{\mu^x}(q)}=1
	$$
	and 
	$$
	\sum_{i=1}^N p_i^q a_i^{\gamma_B(q)-\tau_{\mu^y}(q)} b_i^{\tau_{\mu^y}(q)}=1,
	$$
 where $\mu^x $(resp. $\mu^y$) is the projection of $\mu$ onto $x$-axis (resp. $y$-axis). The result in \cite{F16} states  that $\tau_{\mu}(q)= \max \{\gamma_A(q),\gamma_B(q)\}$ if $\max \{\gamma_A(q),\gamma_B(q)\}\leq \tau_{\mu^x}(q)+\tau_{\mu^y}(q)$, and 
 $$
 \tau_{\mu}(q)\leq \min \{\gamma_A(q),\gamma_B(q)\}
 $$ if $\min \{\gamma_A(q),\gamma_B(q)\}\geq \tau_{\mu^x}(q)+\tau_{\mu^y}(q)$ and equality occurs if 
\begin{subequations}\label{e90}
\begin{equation}\label{e90a}
\textit{ either }\sum_{i=1}^N p_i^q a_i^{\tau_{\mu^x}(q)} b_i^{\gamma_A(q)-\tau_A(q)} \log a_i/b_i\geq 0,
\end{equation}
\begin{equation}\label{e90b}
\textit{ or }\sum_{i=1}^N p_i^q b_i^{\tau_{\mu^y}(q)} a_i^{\gamma_B(q)-\tau_{\mu^y}(q)} \log a_i/b_i\leq 0.
\end{equation}
\end{subequations}
  Naturally it remains  a question  \cite[Question 2.14]{F16} raised by Fraser that whether the additional condition \eqref{e90} can be removed, i.e.
 \begin{question}[{\cite[Question 2.14]{F16}}]\label{que}
     When $\min \{ \gamma_A(q),\gamma_B(q)\} \geq \tau_{\mu^x}(q)+\tau_{\mu^y}(q)$, is $$
     \tau_{\mu}(q)=\min \{ \gamma_A(q),\gamma_B(q)\}
     $$
     still true if \eqref{e90} dose not hold?
 \end{question}
 
 This question was answered by Fraser, Lee, Morris and Yu \cite{FLMY21} in the negative by a special family of counterexamples. In particular, they consider a family of diagonal systems consisting of two maps equipped with a Bernoulli-$(1/2,1/2)$ measure. For this family, it may really happen that 
 \begin{equation}\label{e91}
 \tau_{\mu}(q)<\min \{\gamma_A(q),\gamma_B(q)\}
 \end{equation}
for all $q>1$, and the exact  expression of $\tau_\mu(q)$ was obtained recently by Kolossv\'{a}ry \cite[Proposition 4.4]{K23} in the setting that   grid structure of carpets (could be in high dimensional) are required. 

Nevertheless, it remains unclear that: \vspace{0.2cm}

\textit{$\bullet$ What is the general exact expression of $\tau_\mu(q)$ when $\min \{ \gamma_A(q),\gamma_B(q)\} \geq \tau_{\mu^x}(q)+\tau_{\mu^y}(q)$? }

\textit{$\bullet$ What is the general comparison between the values of $\tau_\mu(q)$ and $\min\{\gamma_A(q),\gamma_B(q)\}$?}\vspace{0.2cm}
 
 All the above considerations require that maps in IFS's are diagonal. \vspace{0.2cm}
 
 \textit{$\bullet$ What would it be when allowing maps in IFS's to be anti-diagonal?}\vspace{0.2cm}

 Along this direction, Morris \cite[Proposition 5]{M19} derived a closed form expression for box dimensions (taking $q=0$ in $\tau_\mu(q)$)
for  self-affine carpets,  requiring that \textit{at least one of  $T_i$'s in IFS's is  anti-diagonal}.
\vspace{0.2cm}

 \begin{figure}[htbp]
		\centering
		\subfigure{		\includegraphics[width=0.35\textwidth]{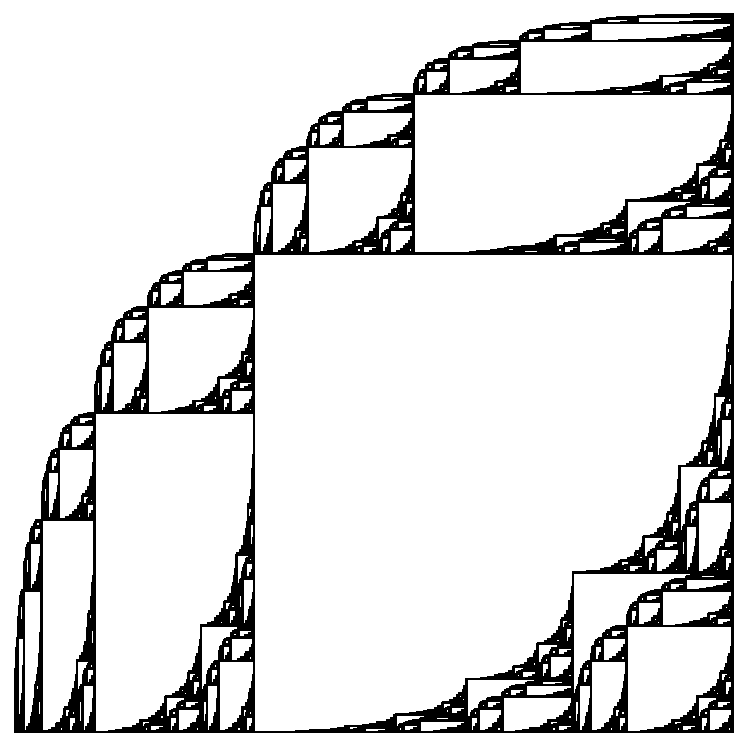}}\hspace{0.6cm}
        \subfigure{\includegraphics[width=0.35\textwidth]{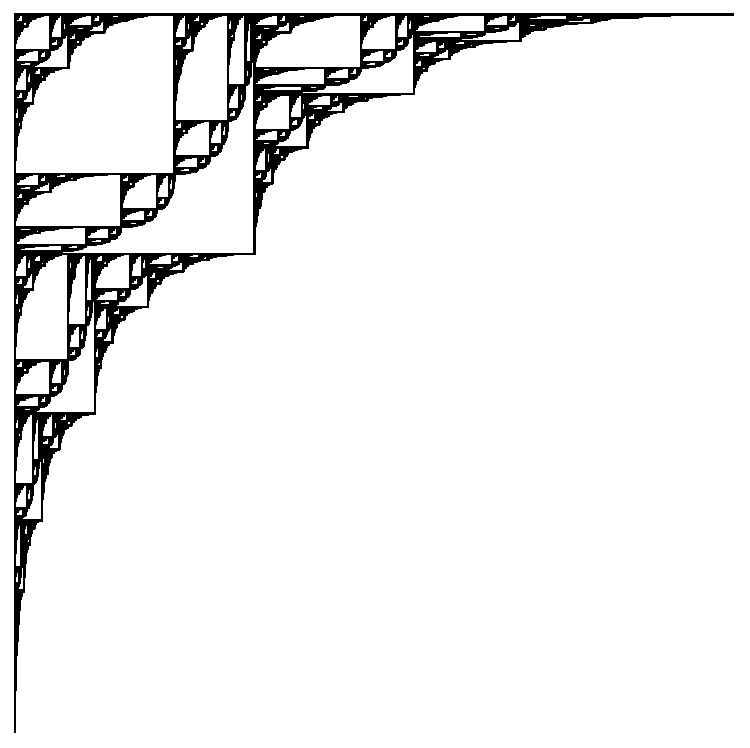}}
		\caption{An example of graph-directed self-affine carpet families.}
		\label{f4}
	\end{figure}
 
Our main aim in this paper is to answer the above questions. We will extend the consideration from the IFS setting to the more general graph-directed IFS (GIFS) setting, allowing contracting maps to be either diagonal or anti-diagonal, i.e. each associated matrix $T_i$ is of the form$$
	T_i= \left(\begin{array}{cc}
		\pm a_i & 0 \\
		0   & \pm b_i \\	\end{array}\right) 
	\text{ or } 
	\left(\begin{array}{cc}
		0 & \pm a_i \\
		\pm b_i & 0 \\
	\end{array}\right).
	$$ See Figure \ref{f4} for an example of  associated graph-directed self-affine carpet families.  We will obtain a general exact closed form expression for $L^q$-spectra of graph-directed self-affine measures, for general $q\geq 0$. Specifically,  returning to the diagonal IFS setting concerned by Question \ref{que}, our result will state that the strict inequality \eqref{e91} generally holds when \eqref{e90} does not hold. Indeed, we will prove 
when $\min \{ \gamma_A(q),\gamma_B(q)\} \geq \tau_{\mu^x}(q)+\tau_{\mu^y}(q)$, \begin{equation}\label{e99}
				\begin{aligned}
					\tau_\mu(q)&=\min \{x+y: \sum_{i=1}^N p_i^q a_i^x b_i^y=1, \tau_A(q)\leq x \leq \gamma_B(q)-\tau_B(q) 
					\}\\
					&= \left \{
		\begin{aligned}
			&\gamma_A(q) &\quad  &\text{ if \eqref{e90a} holds}, \\
			&\gamma_B(q) &\quad&\text{ if \eqref{e90b} holds},
			\\
			&<\min\{\gamma_A(q),\gamma_B(q)\} & \quad & \textit{otherwise. }
		\end{aligned} \right .
				\end{aligned}
			\end{equation}
			Not only that, we will illustrate that the above expression can alternatively be directly derived from Feng and Wang's  original result \cite[Theorem 1]{FW05} by using a careful Lagrange multipliers method.
   Another improvement of \eqref{e99} is that it specifies the necessary and sufficient condition that  $\tau_\mu(q)$  equals to $\gamma_A(q)$ (resp.  $\gamma_B(q)$), compared with that in \cite{F16}.

 When allowing some maps to be anti-diagonal, our result is also a non-trivial extension of that of Morris's   \cite{M19} for box dimension (the $q=0$ case) to all $q\geq 0$ and to the GIFS setting. In his  IFS setting, the graph-directed self-affine measure family degenerates to a single measure $\mu$. The requirement that at least one of $T_i$'s in the IFS is anti-diagonal ensures that  $\tau_{\mu^x}(0)=\tau_{\mu^y}(0)$  (taking $q=0$) since $\{\mu^x,\mu^y\}$ becomes a strongly connected graph-directed self-similar measure family. 
 However, in the GIFS setting, the graph-directed self-affine measure family $\{\mu_v\}_v$ will generate a collection of projection measures $\{\mu_v^x, \mu_v^y\}_v$, which will be proved to be a disjoint union of one or two  strongly connected self-similar measure families, and consequently it may happen that $\tau_{\mu_v^x}(0)\neq \tau_{\mu_v^y}(0)$. This will cause the main difficulty in GIFS setting. Another main difficulty is to properly divide the consideration into distinct cases for distinct $q\geq 0$. 

 The motivation that we extend the consideration to the GIFS setting is the potential application that we can use which to consider box-like self-affine IFS's of finite overlapping types, analogous to that of Ngai and Wang \cite{NW01} and the extension \cite{NL07} for self-similar IFS's. We illustrate this in a recent paper \cite{QWW23} concerning the $L^q$-spectra for lower triangular planar non-conformal  measures. We mention that there are also some  previous works in box-like self-affine GIFS setting. In \cite{KP962}, Kenyon and Peres extended the results of Bedford \cite{B84} and McMullen \cite{M84} computing the Hausdorff and box dimensions of graph-directed self-affine carpets with homogeneous grid structure. In \cite{NW09}, Ni and Wen considered the $L^q$-spectra for graph-directed self-affine measures on Feng and Wang's sets \cite{FW05}, in the setting that contraction ratios on $x$-axis are always less than that on $y$-axis for all maps in the GIFS's, but additionally requiring contraction ratios on $x$-axis to be arithmetic. 






	\vspace{0.2cm}
\noindent{\textit{Basic setting and notations.}}
\vspace{0.2cm}

	The concept of \textit{graph-directed iterated function system (GIFS)} was firstly introduced by Mauldin and Williams \cite{MW88}. 
 
 Let $(V,E)$ be a \textit{finite directed graph} with $V$ being the \textit{vertex set} and $E$ being the \textit{directed edge set}, allowing loops and multiple edges. For $e\in E$, denote by $i(e)$ the \textit{initial vertex}, $t(e)$ the \textit{terminal vertex} of $e$, and sometimes write this as $i(e) \stackrel{e}{\rightarrow} t(e)$. We always assume that for each $v\in V$, there exists at least one edge $e\in E$ satisfying $i(e)=v$.
 
	Denote the collection of all \textit{finite admissible words} by
	$$
	E^*=\{w=w_1 \cdots w_k:t(w_{i-1})=i(w_{i}) , \forall 1<i\leq k, k \in \mathbb{N}\}.
	$$	
	For $w=w_1\cdots w_k\in E^*$, denote $|w|=k$ the \textit{length} of $w$, $i(w)=i(w_1)$, $t(w)=t(w_k)$ the \textit{initial} and \textit{terminal} vertices of $w$, and also write $i(w) \stackrel{w}{\rightarrow} t(w)$.  For $w,w'\in E^*$ with $t(w)=i(w')$, write $ww'$  the \textit{concatenation} of $w$ and $w'$,
	and call $w$ a \textit{prefix} of $ww'$.
	Denote $E^k$ the collection of all admissible words of length $k\geq 1$.
	
	For $v,v'\in V$, we say that there exists a \textit{directed path} from $v$ to $v'$ if there exists $w\in E^*$ satisfying $v \stackrel{w}{\rightarrow}v'$ (write simply $v \to v'$ when we do not emphasize $w$).  Write $v\nrightarrow v'$ if there is not a directed path form $v$ to $v'$. We say $(V,E)$ is \textit{strongly connected} (or \textit{irreducible}), if $v\to v'$ for all pairs $v,v'\in V$. 

	For each $e\in E$, we assume that there exists a contraction $\psi_e$ in the form of $\psi_e(\cdot)=T_e (\cdot) +t_e$, where $T_e$ is an $n\times n$ affine contracting matrix and $t_e \in \mathbb{R}^n$. We write $\Psi =\{\psi_e\}_{e\in E}$ the collection of  all contractions $\psi_e$'s. Call the triple $(V,E,\Psi)$ a \textit{self-affine GIFS}. It is well known that there exists a unique family of compact sets $\{X_v\}_{v\in V}$ satisfying 
	$$
	X_v=\bigcup_{e\in E:i(e)=v} \psi_e(X_{t(e)}),\quad \text{for all } v\in V.
	$$
	We call $\{X_v\}_{v\in V}$ a \textit{graph-directed self-affine set family} associated with $(V,E,\Psi)$. Note that if $V$ is a singleton, $(V,E,\Psi)$ degenerates to a self-affine IFS and $\{X_v\}_{v\in V}$ degenerates to a self-affine set $X$.
 
	Let $\mcP=(p_e)_{e\in E}$ be a positive vector satisfying \begin{equation}\label{psum}
		\sum_{e\in E:i(e)=v}p_e=1, \quad \text{for all } v\in V.
	\end{equation}
	It is known that there exists a unique finite family of probability measures $\{\mu_v\}_{v\in V} :=\{\mu_{\mcP,v}\}_{v\in V}$ supported on $\{X_v\}_{v\in V}$ such that 
	$$
	\mu_v =\sum_{e\in E:i(e)=v} p_e \cdot  \mu_{t(e)}\circ\psi_e^{-1}, \quad \text{for all } v\in V.
	$$
	We call $\{\mu_v\}_{v\in V}$ a \textit{graph-directed self-affine measure family} associated with $\mcP$. 
	\vspace{0.2cm}
	
\begin{figure}[htbp]
		\centering
		\subfigure{		\includegraphics[width=0.35\textwidth]{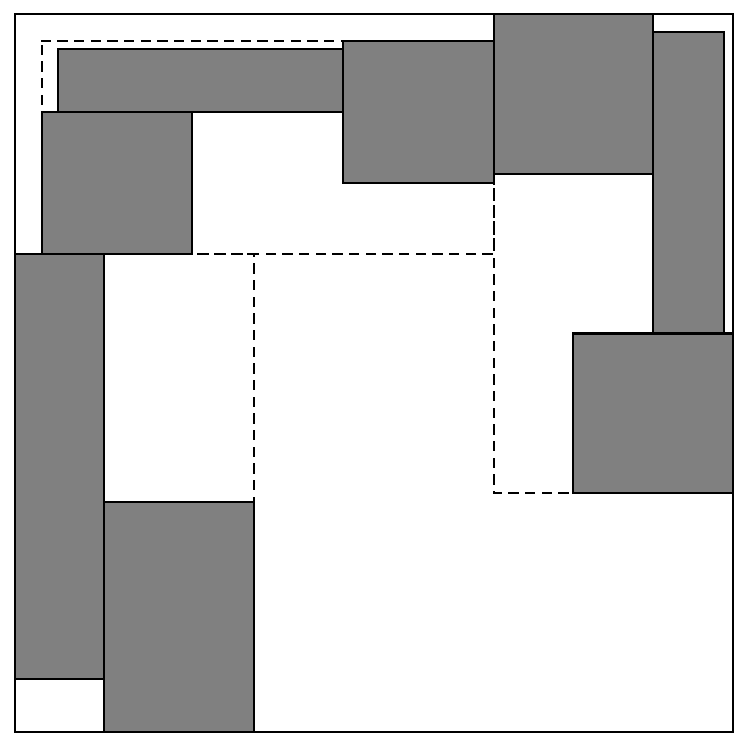}}\hspace{0.6cm}
        \subfigure{\includegraphics[width=0.35\textwidth]{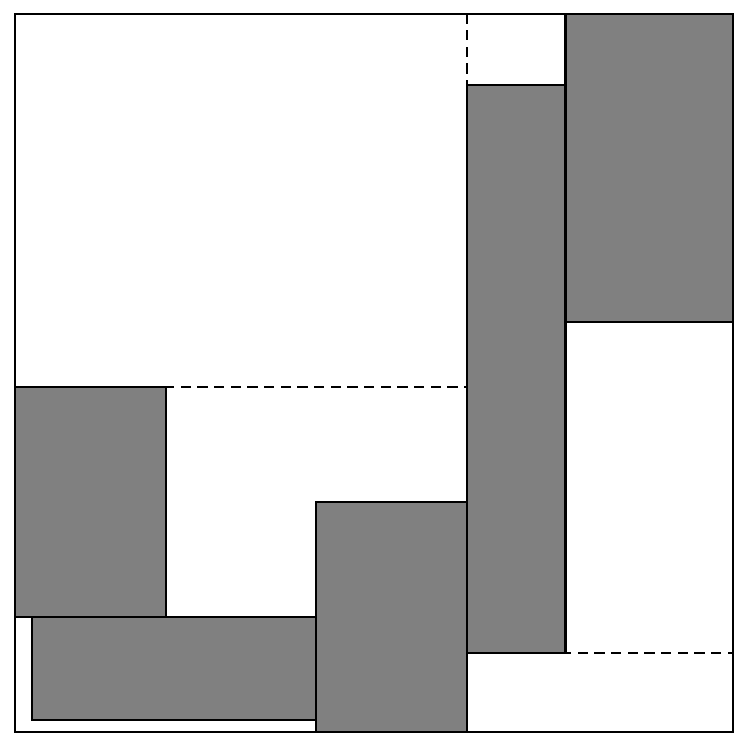}}
		\caption{A planar box-like self-affine GIFS with $\#V=2,\#E=5$. Images of $[0,1]^2$ under the first and second level iterations of maps in the GIFS.}
		\label{f2}
	\end{figure}

\begin{figure}[htbp]
		\centering
		\subfigure{		\includegraphics[width=0.35\textwidth]{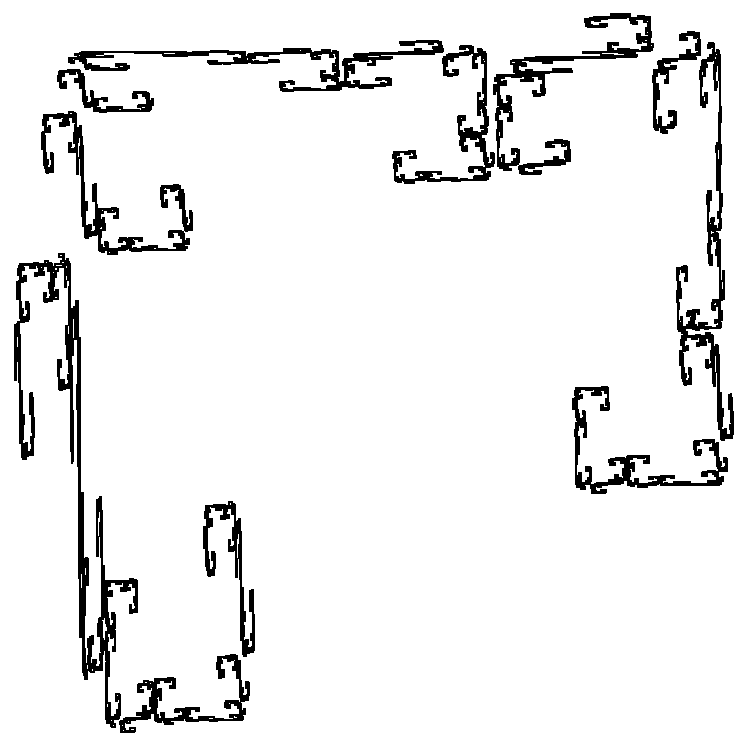}}\hspace{0.6cm}
        \subfigure{\includegraphics[width=0.35\textwidth]{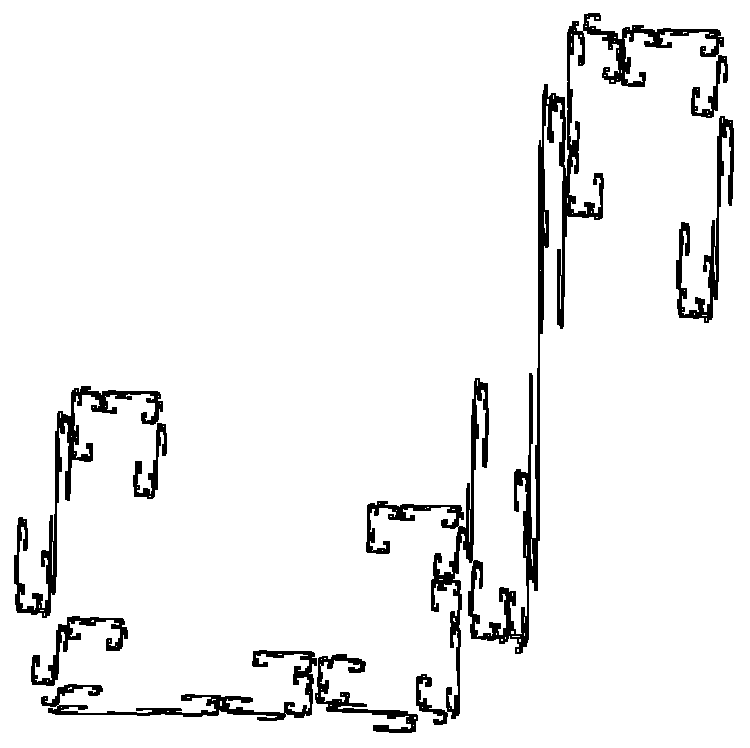}}
		\caption{The graph-directed self-affine carpet family generated by the GIFS in Figure \ref{f2}.}
		\label{f3}
	\end{figure}

	Throughout the paper, we assume that for each $e\in E$, $T_e$ is a  $2\times 2$ diagonal or anti-diagonal matrix of the form 
	\begin{equation}\label{t15}
	T_e= \left(\begin{array}{cc}
		\pm a_e & 0 \\
		0   & \pm b_e \\	\end{array}\right) 
	\text{ or } 
	\left(\begin{array}{cc}
		0 & \pm a_e \\
		\pm b_e & 0 \\
	\end{array}\right)
	\end{equation}
	where $0<a_e,b_e <1$. We call $(V,E,\Psi)$ a \textit{planar box-like self-affine GIFS},   $\{X_v\}_{v\in V}$ a \textit{graph-directed self-affine carpet family} associated with $(V,E,\Psi)$  and $\{\mu_v\}_{v\in V}$ a \textit{graph-directed box-like self-affine measure family}. See Figures \ref{f2}-\ref{f3} for  an example.

 In this paper, we care about the  \textit{ $L^q$-spectra of strongly connected planar graph-directed box-like self-affine measures} $\{\mu_v\}_{v\in V}$.
 For calculating the $L^q$-spectra, we need the following separating condition for the planar box-like self-affine GIFS's, which was firstly proposed by Feng and Wang \cite{FW05} and 
 plays crucial roles in subsequent works \cite{F12,F16,FLMY21,FFL21}.
	
	\begin{definition}[Rectangular open set condition] 
		We say a planar box-like self-affine GIFS $(V,E,\Psi)$ satisfies the \emph{rectangular open set condition} (ROSC) if for all $v$  in $V$, 
		$$
		\bigcup_{e\in E:i(e)=v}\psi_e((0,1)^2)\subseteq(0,1)^2
		$$
		and the union is disjoint.
	\end{definition}
	
	Our results rely on the $L^q$-spectra of the projections of measures onto $x$-axis or $y$-axis.
	Let $\pi_x$, $\pi_y:\mathbb{R}^2\to \mathbb{R}$ be defined by $\pi_x(\xi_1,\xi_2)=\xi_1$ and $\pi_y(\xi_1,\xi_2)=\xi_2$ for $(\xi_1,\xi_2) \in \mathbb{R}^2$, respectively. 
 For each $v$ in $V$, define 
	$$\mu_v^x=\mu_v\circ \pi_x^{-1}, \quad \mu_v^y=\mu_v\circ \pi_y^{-1}
	$$
	the \textit{projection measures} of $\mu_v$ onto $x$-axis and $y$-axis.
	Note that $\{\mu_v^x,\mu_v^y\}_{v\in V}$ is a family of graph-directed self-similar measures. 
	
	\begin{proposition}[{\cite[Theorem 2.1]{F16}}] \label{prop1}
		For a  strongly connected graph-directed self-similar measure family $\{\nu_v\}_{v\in V}$, for all $q\geq 0$ and $v,v'\in V$, we have 
		$$
		\underline{\tau}_{\nu_v}(q)=\overline{\tau}_{\nu_v}(q) =\underline{\tau}_{\nu_{v'}}(q)=\overline{\tau}_{\nu_{v'}}(q),
		$$
		i.e. the $L^q$-spectra exist and are the same for all $\nu_v$'s.
	\end{proposition}
	
	When  $(V,E)$ is strongly connected and all $T_e$'s are diagonal, both $\{\mu_v^x\}_{v\in V}$ and $\{\mu_v^y\}_{v\in V}$ are two strongly connected graph-directed self-similar measure families; when $V$ is a singleton and $T_e$ is anti-diagonal for some $e\in E$, $\{\mu_v^x, \mu_v^y\}_{v\in V}$ is a strongly connected graph-directed self-similar measure family; for general strongly connected case, $\{\mu_v^x, \mu_v^y\}_{v\in V}$ can be divided into one or two families of strongly connected graph-directed self-similar measures (see Proposition \ref{thmain4}).
	By Proposition \ref{prop1}, we always have the $L^q$-spectra exist for all $\mu_v^x$'s, $\mu_v^y$'s.
	
	\vspace{0.2cm}
	
	Throughout the paper, we will write $a\lesssim b$ for two variables (functions) if there is a constant $C>0$ such that $a\leq C\cdot b$, and write $a\asymp b$ if both $a\lesssim b$ and $b\lesssim a$ hold. We write $a\lesssim_\theta b$ to mean that the constant depends on some parameter $\theta$. Similarly, write $a \asymp_\theta b$ if both $a\lesssim_\theta b$ and $b\lesssim_\theta a$ hold.
	For two vectors $u=(u_i)$ and $v=(v_i)$, we write $u\geq v$ if all $u_i\geq v_i$. Also for two matrices $A=(a_{ij})$ and $B=(b_{ij})$, we write $A\geq B$ if all $a_{ij}\geq b_{ij}$. For a $N\times N$ matrix $A$, for indices $\{i_1,\cdots,i_k\}, \{j_1,\cdots, j_l\}$ with $k,l\leq N$, we write 
	$A[\{i_1,\cdots,i_k\},\{j_1,\cdots, j_l\}]$ for a submatrix of $A$ (lying in
	rows $\{i_1,\cdots,i_k\}$ and columns $\{j_1,\cdots, j_l\}$).
	We always denote $\Vert \cdot \Vert$  the \textit{$1$-norm} of a matrix, i.e. $\Vert A \Vert =\sum_{i,j} |a_{ij}|$.

	\section{Results}\label{sec7}
	In this section, we list the results in the paper but postpone their proofs to later sections. Our main aim is to obtain the closed form expression for the 
	 $L^q$-spectra of planar graph-directed box-like self-affine measures. 
  
	Throughout the following, we always let $(V,E,\Psi)$ be a strongly connected planar box-like self-affine  GIFS, $\mathcal{P}$ be a positive vector satisfying \eqref{psum}, and $\{\mu_v\}_{v\in V}$ be a graph-directed box-like self-affine measure family associated with $(V,E,\Psi)$ and $\mathcal{P}$. Note that for each $e\in E$, there exists a contraction $\psi_e$ in the form of $\psi_e(\cdot)=T_e(\cdot)+t_e$ for some $2\times 2$  diagonal or anti-diagonal contracting matrix $T_e$ and $t_e\in \mathbb{R}^2$. For $v\in V$,
	we use  $\mu_v^x$ (resp. $\mu_v^y$) to denote the projection of $\mu_v$ onto $x$-axis (resp. $y$-axis). 
 
 We will separate our consideration basing on two basic settings:\vspace{0.2cm}
 
 \textit{first, assume all $T_e$'s are diagonal}; 
 
 \textit{then, extend the consideration to general case, i.e. allowing some $T_e$'s to be anti-diagonal}. 
	\vspace{0.2cm}
 
 Before proceeding, we will prove in general that the $L^q$-spectra of measures $\tau_{\mu_v^x}(q)$, $\tau_{\mu_v^y}(q)$ and $\tau_{\mu_v}(q)$ exist for all $q\geq 0$ and $v\in V$, which will play  fundamental roles in our later consideration. We will achieve this through the following Proposition \ref{thmain4} and Theorem \ref{thmain1} dealing with projection measures and original measures separately. 
 
	\begin{proposition}
		\label{thmain4}
		Let $(V,E,\Psi)$ be a strongly connected planar box-like self-affine  GIFS, $\mathcal P$ be an associated positive vector, and $\{\mu_v\}_{v\in V}$ be an associated graph-directed box-like self-affine measure family. Then  
		\begin{equation}\label{s66}
			\tau_{\mu_v^x}(q) \text{ and } \tau_{\mu_v^y}(q) \text{ exist } \quad \text{ for all } q\geq 0, v\in V.
		\end{equation}
		Moreover, $\{\mu_v^x,\mu_v^y\}_{v\in V}$ can be divided into two disjoint families $A$ and $B$ so that $$\#A=\#B=\#V,$$  for all $q\geq 0$ there exist $\tau_A(q),\tau_B(q)$ satisfying 
		\begin{equation}\label{s67}
			\begin{aligned}
				\tau_{\nu}(q)=\tau_A(q)\quad \text{ for all } \nu \in A,\\
				\tau_{\nu}(q)=\tau_B(q)\quad \text{ for all } \nu \in B,
			\end{aligned}
		\end{equation}
		and for all $v\in V$, 
		\begin{equation}\label{s65}
			\text{ either }\quad  \mu_v^x \in A, \mu_v^y\in B \quad \text{ or } \quad  \mu_v^x \in B, \mu_v^y \in A.
		\end{equation}
		In particular,
		when all $T_e$'s are diagonal, we could take
		$$
		A=\{ \mu_v^x\}_{v\in V} \text{ and } B=\{ \mu_v^y\}_{v\in V}.
		$$
	\end{proposition}
	Due to Proposition \ref{thmain4}, for  $q\geq 0$, $v\in V$, $e\in E$, throughout the paper,  we will  write $$
	\begin{aligned}
		\tau_{x,v}(q)&:=\tau_{\mu_v^x}(q)\quad  & & \text{ and }\quad & & \tau_{y,v}(q):=\tau_{\mu_v^y}(q),
		\\
		\tau_{x,e}(q)&:=\tau_{x,t(e)}(q) & & \text{ and } & &\tau_{y,e}(q):=\tau_{y,t(e)}(q).
	\end{aligned}
	$$ for short. Also, write $$t(q):=\tau_A(q)+\tau_B(q)$$ for later use. Clearly, for all $v\in V$, $\tau_{x,v}(q)+\tau_{y,v}(q)=t(q)$.\vspace{0.2cm}
	
	Next, with Proposition \ref{thmain4} in hand, inspired by Fraser's works \cite{F16,F12} dealing with the self-affine carpets, i.e. the case that $(V,E)$ degenerates to a singleton, we will introduce (in Section \ref{sec2}) a \textit{pressure function} in graph-directed setting, $$P:\mathbb{R}\times[0,\infty) \to \mathbb{R},$$ basing on certain modified singular value function matrices. For each $q\geq 0$, as a function of $s$, $P(s,q)$ will be strictly decreasing and continuous, tending to $0$ as $s \to +\infty$ and  to $+\infty$ as $s\to -\infty$. Using this we will define a function $\gamma:[0,\infty)\to \mathbb{R}$ by setting $$P(\gamma(q),q)=1.$$
	
	\begin{theorem}\label{thmain1}
		Let $(V,E,\Psi)$ be a strongly connected planar box-like self-affine  GIFS satisfying ROSC, 
$\mathcal P$ be an associated positive vector,  and $\{\mu_v\}_{v\in V}$ be an associated graph-directed box-like self-affine measure family.  Then for all $q\geq 0$ and $v\in V$, $\tau_{\mu_v}(q)$ exist and equal to $\gamma(q)$. 
	\end{theorem}

 Proposition \ref{thmain4} will be proved in Section \ref{sec2}. The details for pressure function $P$ will be presented also in Section \ref{sec2}. The proof of Theorem \ref{thmain1} will be postponed to Section \ref{sec3}.\vspace{0.2cm}
	
	From definition, for $q\geq 0$, $\gamma(q)$ does not seem to be able to be explicitly computed through a finite  amount of steps. So our next aim attributes to find a closed form expression for  $\gamma(q)$, we will separate our consideration into two parts.\vspace{0.2cm}

 \noindent{\textit{Non-rotational setting.}}\vspace{0.2cm}
 
 First we will assume that all $T_e$'s are  diagonal matrices, i.e. each $T_e$ is of the form
	$$
	T_e= \left(\begin{array}{cc}
		\pm a_e & 0 \\
		0   & \pm b_e \\	\end{array}\right) .
	$$
In this setting, by Proposition \ref{thmain4}, $A=\{\mu_v^x\}_{x\in V}$, $B=\{\mu_v^y\}_{v\in V},$ and for all $v\in V$, $\tau_{x,v}(q)=\tau_A(q)$, $\tau_{y,v}(q)=\tau_B(q)$ and $\tau_{x,v}(q)+\tau_{y,v}(q)=t(q).$ 
	For $q\geq 0$, $x,y\in \mathbb{R}$, we will introduce (in Section \ref{sec4}) a $\#E \times \#E$ function matrix $F^{(q)}_{x,y}$ with entries defined by
	\begin{equation}\label{s14}
		F^{(q)}_{x,y}(e,e') =\left \{ \begin{aligned}
			&p_{e'}^q a_{e'}^x b_{e'}^y &\quad& \text{ if } t(e)=i(e'),\\
			&0 &\quad& \text{ otherwise}.
		\end{aligned} \right.
	\end{equation}
	Then define two functions $\gamma_A,\gamma_B:[0,\infty]\to \mathbb{R}$ such that for $q\geq 0$, $\gamma_A(q)$ and $\gamma_B(q)$ are the unique solutions of  
	\begin{equation}\label{t7}
		\rho(F^{(q)}_{\tau_A(q),\gamma_A(q)-\tau_A(q)})=1
	\end{equation}
	and 
	\begin{equation}\label{t8}
		\rho(F^{(q)}_{\gamma_B(q)-\tau_B(q),\tau_B(q)})=1,
	\end{equation}
	respectively (to be  well-defined in Section \ref{sec4}), where $\rho(\cdot)$ is the \textit{spectral radius} of a matrix.
	For fixed $q\geq 0$, for $x\in \mathbb{R}$, we will prove that there exists a unique $y(x)\in \mathbb{R}$ satisfying $\rho(F^{(q)}_{x,y(x)})=1$ and introduce a positive unit row vector $f^{(q)}(x)=(f^{(q)}_e(x))_{e\in E}$ (in Section \ref{sec4}).
	
	\begin{theorem}\label{thmain2}
		Let $(V,E,\Psi)$, $\mcP$, $\{\mu_v\}_{v \in V}$ and $\gamma$ be same as in Theorem \ref{thmain1}. Assume that all $T_e$'s are diagonal. Then for $q\geq 0 $,
	\begin{subequations}	\begin{equation}\label{e1}
		\text{ either }	\max \{\gamma_A(q), \gamma_B(q)\} \leq t(q)
		\end{equation}
\begin{equation}\label{e2}
	\text{ \hspace{0.2cm} or \hspace{0.4cm}}			\min \{\gamma_A(q), \gamma_B(q)\} \geq t(q).
		\end{equation}
  \end{subequations}
		\begin{itemize}
			\item[(a).] If \eqref{e1} holds,
			\begin{equation}\label{s20}
				\begin{aligned} 
					\gamma(q)&=\max \{\gamma_A(q), \gamma_B(q)\}\\
					&= \max \{x+y: \rho(F^{(q)}_{x,y})=1, \gamma_B(q)-\tau_B(q) \leq x \leq \tau_A(q)
					\}\\
					&=\max \{x+y: \rho(F^{(q)}_{x,y})=1, \gamma_A(q)-\tau_A(q) \leq y \leq \tau_B(q)
					\}.
				\end{aligned}
			\end{equation}
			\item[(b).] If \eqref{e2} holds, \begin{equation}\label{s21}
				\begin{aligned}
					\gamma(q)&=\min \{x+y: \rho(F^{(q)}_{x,y})=1, \tau_A(q)\leq x \leq \gamma_B(q)-\tau_B(q) 
					\}\\
					&=\min \{x+y: \rho(F^{(q)}_{x,y})=1, \tau_B(q)\leq y \leq \gamma_A(q)-\tau_A(q) 
					\}.\\
				\end{aligned}
			\end{equation} 
			Moreover,
			\begin{itemize}
				\item[(b1).] if $\sum_{e\in E}f^{(q)}_e(\tau_A(q)) \log (a_e/b_e) \geq 0$, $$\gamma(q)=\gamma_A(q),$$ 
				\item[(b2).] if $\sum_{e\in E}f^{(q)}_e(\gamma_B(q)-\tau_B(q)) \log (a_e/b_e) \leq  0$, $$\gamma(q)=\gamma_B(q),$$ 
				\item[(b3).] otherwise, there exist $x\in [\tau_A(q),\gamma_B(q)-\tau_B(q)]$ and  $y\in [\tau_B(q), \gamma_A(q)-\tau_A(q)]$ with $\rho(F^{(q)}_{x,y})=1$ satisfying  $\sum_{e\in E} f^{(q)}_e(x) \log (a_e/b_e) =0$ and 
				$$
				\gamma(q)=x+y.
				$$
			\end{itemize}
		\end{itemize}
		
	\end{theorem}
	
	Recall that when $q=0$, the  $L^q$-spectrum of a measure $\nu$ is equal to the box dimension of $supp \nu$. 
	\begin{corollary}\label{cor2}
		
		Let $(V,E,\Psi)$ be same as in Theorem \ref{thmain2}. Let $\{X_v\}_{v\in V}$ be the unique graph-directed  self-affine carpet family associated with $(V,E,\Psi)$. Then for all $v\in V$, 
		$$
		\dim_B X_v =\max \{ \gamma_A(0),\gamma_B(0)\}.
		$$
	\end{corollary}
	
	When $V$ is a singleton, $(V,E,\Psi)$ degenerates to a  box-like self-affine  IFS, $\mathcal P$ degenerates to a positive probability vector, and $\{\mu_v\}_{v\in V}$ degenerates to a single measure $\mu$. The directed edge set $E$ can be written as $\{1,\cdots, N\}$. At this time, all rows of the matrix $F^{(q)}_{x,y}$ are same. So by Perron-Frobenius Theorem, $\rho(F^{(q)}_{x,y})= \sum_{i=1}^N p_i^q a_i ^x b_i^y$. Thus $\gamma_A,\gamma_B$ can be reduced to the unique solutions satisfying
	$$
	\sum_{i=1}^N p_i^q a_i^{\tau_A(q)} b_i^{\gamma_A(q)-\tau_A(q)}=1
	$$
	and 
	$$
	\sum_{i=1}^N p_i^q a_i^{\gamma_B(q)-\tau_B(q)} b_i^{\tau_B(q)}=1,
	$$
	respectively. Also, for  $q\geq 0,x\in \mathbb{R}$, $f^{(q)}(x)$ will be reduced to $f^{(q)}(x)=(p_i^q a_i^x b_i^y)_{i=1}^N$ where $y$ is the unique solution satisfying $\sum_{i=1}^{N} p_i^q a_i^x b_i^y=1$ (see details in Section \ref{sec4}). 
	Therefore, the conditions (b1), (b2) of Theorem \ref{thmain2} will become 
	\begin{subequations}
\begin{equation}\label{c1}
		\sum_{i=1}^N p_i^q a_i^{\tau_A(q)} b_i^{\gamma_A(q)-\tau_A(q)} \log a_i/b_i\geq 0,
	\end{equation} 
 and
\begin{equation}\label{c2}
		\sum_{i=1}^N p_i^q b_i^{\tau_B(q)} a_i^{\gamma_B(q)-\tau_B(q)} \log a_i/b_i\leq 0.
	\end{equation}
 \end{subequations}
	Then we will obtain the following corollary,  a precise answer to \cite[Question 2.14]{F16}.
	\begin{corollary}\label{cor1}
		Let $\{\psi_i(\cdot) =T_i (\cdot)+t_i\}_{i=1}^N$ be a box-like self-affine IFS satisfying ROSC. Assume that all $T_i$'s are diagonal. Let $\mathcal{P}=(p_i)_{i=1}^N$ be a positive probability vector. Let $\mu$ be the self-affine measure associated with $\mathcal{P}$.  For $q\geq 0$, \begin{subequations}
  \begin{equation}\label{e3}
	\text{ either }		\max \{\gamma_A(q), \gamma_B(q)\} \leq t(q)
		\end{equation}
\begin{equation}\label{e4}
		\text{ \hspace{0.2cm} or \hspace{0.4cm}}	\min \{\gamma_A(q), \gamma_B(q)\} \geq t(q).
		\end{equation}
  \end{subequations}
		\begin{itemize}
			\item[(a).] If \eqref{e3} holds,
			\begin{equation}
				\nonumber 
				\begin{aligned} 
					\tau_\mu(q)&=\max \{\gamma_A(q), \gamma_B(q)\}\\
					&= \max \{x+y: \sum_{i=1}^N p_i^q a_i^x b_i^y=1, \gamma_B(q)-\tau_B(q) \leq x \leq \tau_A(q)
					\}\\
					&= \max \{x+y: \sum_{i=1}^N p_i^q a_i^x b_i^y=1, \gamma_A(q)-\tau_A(q) \leq y \leq \tau_B(q)
					\}.
				\end{aligned}
			\end{equation}
			\item[(b).] If \eqref{e4} holds, \begin{equation}\nonumber 
				\begin{aligned}
					\tau_\mu(q)&=\min \{x+y: \sum_{i=1}^N p_i^q a_i^x b_i^y=1, \tau_A(q)\leq x \leq \gamma_B(q)-\tau_B(q) 
					\}\\
					&=\min \{x+y: \sum_{i=1}^N p_i^q a_i^x b_i^y=1, \tau_B(q)\leq y \leq \gamma_A(q)-\tau_A(q) 
					\}.
				\end{aligned}
			\end{equation}
			Moreover,
			\begin{itemize}
				\item[(b1).] if \eqref{c1} holds, $\tau_\mu(q)=\gamma_A(q),$ 
				\item[(b2).] if \eqref{c2} holds, $\tau_\mu(q)=\gamma_B(q),$ 
				\item[(b3).] otherwise, there exist $x\in [\tau_A(q),\gamma_B(q)-\tau_B(q)]$ and $y \in [\tau_B(q),\gamma_A(q)-\tau_A(q)]$ satisfying $\sum_{i=1}^N p_i^q a_i^x b_i^y=1$ such that  $\sum_{i=1}^N p_i^q a_i^x b_i^y\log (a_i/b_i) =0$ and
				$$
				\tau_\mu(q)=x+y < \min \{\gamma_A(q),\gamma_B(q)\}.
				$$
			\end{itemize}
		\end{itemize}
	\end{corollary}
	\begin{remarks}\label{r1}
		\emph{(a).} There is a slight difference between the statements in Theorem \ref{thmain2} and Corollary \ref{cor1}. For the case (b3) in Corollary \ref{cor1}, we can further know that $$\tau_\mu(q)<\min\{\gamma_A(q),\gamma_B(q)\}.$$

		\emph{(b).} Fraser \cite[Theorem 2.10]{F16} proved  that  if \eqref{e4} holds, 
        $
        \tau_\mu(q)\leq \min \{ \gamma_A(q),\gamma_B(q)\}
        $
        with equality if either \eqref{c1} or \eqref{c2} holds. Indeed, he proved that $\tau_{\mu}(q)=\gamma_A(q)$ if \eqref{c1} holds, and $\tau_{\mu}(q)=\gamma_B(q)$ if \eqref{c2} holds, but did not specify the sufficient condition for $\tau_\mu(q)=\gamma_A(q)$ (resp. $\gamma_B(q)$).
		
		\emph{(c).}
		Corollary \ref{cor1}  can also  be directly derived by using \cite[Theorem 1]{FW05} without using Theorem \ref{thmain2}, see Subsection \ref{subsec32}.
	\end{remarks}

 The non-rotational setting will be considered in Section \ref{sec4}, where Theorem \ref{thmain2} and Corollaries \ref{cor2}, \ref{cor1} will be proved. Particularly, we will provide an alternative proof of Corollary \ref{cor1} in Subsection \ref{subsec32}.\vspace{0.2cm}


 \noindent{\textit{General setting.}}\vspace{0.2cm}
	
	Next we will turn to the general setting by allowing some  $T_e$ to be anti-diagonal, i.e. some of $T_e$'s may be of the form
 $$
	T_e= \left(\begin{array}{cc}
		0 & \pm a_e \\
		\pm b_e   & 0 \\	\end{array}\right) .
	$$ Note that when no such $T_e$ exists, this reduces to the non-rotational setting.
 
 Before proceeding, we point out that,  considering a box-like self-affine  IFS $\{\psi_i\}_{i=1}^N$ (i.e. the case that $(V,E)$ degenerates to a singleton), requiring that their exists at least one anti-diagonal $T_e$, Morris \cite{M19} has derived  a closed form expression for  the box dimension of the associated self-affine  carpet. In his setting, $\{\mu_v\}_{v\in V}$ degenerates to a single measure $\mu$, $\tau_A(0)=\tau_B(0)$ (taking $q=0$) since $\{\mu^x,\mu^y\}$ becomes a strongly connected graph-directed self-similar measure family. To deal with the general $q\geq 0$ and general GIFS setting, inspired by his work, we will replace the $\#E\times \#E$ matrix $F_{x,y}^{(q)}$ considered in non-rotational setting by an $\#E\times \#E$  block matrix $\mathcal{G}_{x,y}^{(q)}$ with entries being $2\times 2$ matrices according to the rotational or anti-rotational choice of each  $T_e$. However, the main difficulties emerge from two aspects:  firstly, it is non-trivial to adapt and extend some ideas of the proof for the non-rotational setting, in particular, to properly divide the consideration into distinct cases for distinct  $q\geq 0$; secondly, due to Proposition \ref{thmain4}, it may happen $\tau_A(q)\neq \tau_B(q)$ and so the $(2\#E)\times (2\#E)$ matrix $\mathcal{G}_{x,y}^{(q)}$ is not always irreducible. 
 
    \vspace{0.2cm}
	
	For $q\geq 0,x,y\in \mathbb{R}$, $e,e'\in E$, define a $2\times 2$ matrix with indices $\{e(1),e(2)\}\times \{e'(1),e'(2)\}$,
	\begin{equation}\nonumber
		G^{(q)}_{x,y,e,e'} = \left \{
		\begin{aligned}
			&\left(\begin{array}{cc}
				p_{e'}^q a_{e'}^{x+\tau_{x,e'}(q)} b_{e'}^{y-\tau_{x,e'}(q)} & 0 \\
				0   & p_{e'}^q b_{e'}^{x+\tau_{y,e'}(q)} a_{e'}^{y-\tau_{y,e'}(q)} \\	\end{array}\right) &\quad  &\text{ if }  t(e)=i(e') \text{ and }T_{e'} \text{ is diagonal, } \\
			&\left(\begin{array}{cc}
				0 & p_{e'}^q a_{e'}^{x+\tau_{y,e'}(q)} b_{e'}^{y-\tau_{y,e'}(q)} \\ p_{e'}^q
				b_{e'}^{x+\tau_{x,e'}(q)} a_{e'}^{y-\tau_{x,e'}(q)}   &0 \\	\end{array}\right)  &\quad&\text{ if } t(e)=i(e') \text{ and } T_{e'} \text{ is anti-diagonal,}
			\\
			&\left(\begin{array}{cc}
				0 & 0 \\
				0   &0 \\	\end{array}\right) & \quad & \textit{otherwise. }
		\end{aligned} \right .
	\end{equation}
	Define a $\# E \times \# E$ block function matrix $\mathcal{G}_{x,y}^{(q)}$ with entries being $2\times 2$ matrices,
	\begin{equation}\nonumber
		\mathcal{G}_{x,y}^{(q)}[\{e(1),e(2)\},\{e'(1),e'(2)\}] =
		G_{x,y,e,e'}^{(q)},
	\end{equation}
	and regard $\mathcal{G}_{x,y}^{(q)}$ as a $(2\#E)\times (2\#E)$ matrix. Define a function $\hat{\gamma}(q):[0,+\infty) \to \mathbb{R}$ so that for each $q$, $\hat{\gamma}(q)$  is the unique solution of 
	\begin{equation}\label{s28}
		\rho(\mathcal{G}^{(q)}_{0,\hat{\gamma}(q)})=\rho(\mathcal{G}^{(q)}_{\hat{\gamma}(q)-t(q),t(q)})=1
	\end{equation}
	(to be well-defined in Section \ref{subsec33}). The following theorem is our main result.
	\begin{theorem}\label{thmain3}
		Let $(V,E,\Psi),\mathcal{P}, \{\mu_v\}_{v\in V}$ and $\gamma$ be same as in Theorem \ref{thmain1}.  For $q\geq 0$, we have 
		\begin{itemize}
			\item[(a).] if $\hat{\gamma}(q)\leq t(q)$,
			\begin{equation}\label{s50}
				\begin{aligned} 
					\gamma(q)&=\hat{\gamma}(q)\\
					&= \max \{x+y: \rho(\mathcal{G}^{(q)}_{x,y})=1, \hat{\gamma}(q)-t(q) \leq  x \leq 0
					\}\\
					&=\max \{x+y: \rho(\mathcal{G}^{(q)}_{x,y})=1, \hat{\gamma}(q) \leq  y \leq t(q)
					\},
				\end{aligned}
			\end{equation}
			\item[(b).] if $\hat{\gamma}(q)>t(q)$,
			\begin{equation}\label{s58}
				\begin{aligned}
					\gamma(q)&=\min \{x+y: \rho(\mathcal{G}^{(q)}_{x,y})=1, 0 \leq x \leq \hat{\gamma}(q)-t(q)
					\}\\
					&=\min \{x+y: \rho(\mathcal{G}^{(q)}_{x,y})=1, t(q)\leq y \leq \hat{\gamma}(q)\}.\\
				\end{aligned}
			\end{equation} 
		\end{itemize}
	\end{theorem}
	\begin{corollary}\label{cor3}
		Let $(V,E,\Psi)$ be same as in Theorem \ref{thmain3}. Let $\{X_v\}_{v\in V}$ be the unique graph-directed  self-affine carpet family associated with $(V,E,\Psi)$. Then for all $v\in V$, 
		$$
		\dim_B X_v  =\hat{\gamma}(0).
		$$
	\end{corollary}
 
 Return to the case that $V$ is a singleton.
Let $E=\{1,\cdots ,N\}$ and $\mu$ be the associated self-affine measure.  Without loss of generality, by rearranging the order of $\{\psi_i\}_{i=1}^N$, we can assume that there is a $k\in \{1,\cdots ,N+1\}$ so that 
	$$
	T_i=\left\{ \begin{aligned} &\left(\begin{array}{cc}
			\pm a_i & 0 \\
			0   & \pm b_i \\	\end{array}\right) \quad& & \text{ if } 1\leq i<  k,
		\\
		&\left(\begin{array}{cc}
			0 & \pm a_i \\
			\pm b_i & 0 \\	\end{array}\right) \quad& & \text{ if } k\leq i\leq N.
	\end{aligned}\right .
	$$
	Note that when $k=N+1$, all $T_i$'s are diagonal.
 
	For $q\geq 0, x,y\in \mathbb{R}$, define a $2 \times 2 $ function matrix 
	\begin{equation}\label{s61}
		H^{(q)}_{x,y}=\left( \begin{array}{cc}
			\sum\limits_{i=1}^{k-1} p_i^q a_i^{x+\tau_{\mu^x}(q)} b_i^{y-\tau_{\mu^x}(q)} & \sum\limits_{i=k}^N p_i^q a_i^{x+\tau_{\mu^y}(q)} b_i^{y-\tau_{\mu^y}(q)}\\
			\sum\limits_{i=k}^N p_i^q b_i^{x+\tau_{\mu^x}(q)} a_i^{y-\tau_{\mu^x}(q)} & \sum\limits_{i=1}^{k-1} p_i^q b_i^{x+\tau_{\mu^y}(q)} a_i^{y-\tau_{\mu^y}(q)}\\	\end{array}\right).
	\end{equation}
	\begin{corollary}\label{cor4}
		Let $\{\psi_i(\cdot) =T_i (\cdot)+t_i\}_{i=1}^N$ be a box-like self-affine IFS  satisfying ROSC.
		Let $\mathcal{P}=(p_i)_{i=1}^N$ be a positive probability vector. Let $\mu$ be the self-affine measure associated with $\mathcal{P}$.  For $q\geq 0$, $\hat{\gamma}(q)$ satisfies
		\begin{equation}\label{s59}
			\rho(H^{(q)}_{0,\hat{\gamma}(q)})=\rho(H^{(q)}_{\hat{\gamma}(q)-t(q),t(q)})=1.
		\end{equation}
  In addition,
		\begin{itemize}
			\item[(a).] if $\hat{\gamma}(q)\leq t(q)$,
			\begin{equation}\nonumber
				\begin{aligned} 
					\tau_{\mu}(q)&=\hat{\gamma}(q)\\
					&= \max \{x+y: \rho(H^{(q)}_{x,y})=1, \hat{\gamma}(q)-t(q)\leq  x \leq 0
					\}\\
					&=\max \{x+y: \rho(H^{(q)}_{x,y})=1, \hat{\gamma}(q)  \leq  y \leq t(q)
					\},
				\end{aligned}
			\end{equation}
			\item[(b).] if $\hat{\gamma}(q)>t(q)$,
			\begin{equation}\nonumber
				\begin{aligned}
					\tau_{\mu}(q)&=\min \{x+y: \rho(H^{(q)}_{x,y})=1, 0\leq x \leq \hat{\gamma}(q)-t(q) 
					\}\\
					&=\min \{x+y: \rho(H^{(q)}_{x,y})=1, t(q)\leq y \leq \hat{\gamma}(q)\}.
				\end{aligned}
			\end{equation} 
		\end{itemize}
	\end{corollary}	
	\begin{remark}
		Morris' closed form expression \cite[Proposition 5]{M19} for box dimensions of self-affine carpets can be seen by combining Corollaries \ref{cor3} and \ref{cor4} together and taking  $q=0$.
	\end{remark}
	
	We will consider the general setting in Section \ref{subsec33}, where we will prove Theorem \ref{thmain3} and Corollaries \ref{cor3}, \ref{cor4}.

	\section{Pressure functions}\label{sec2}
	Let $(V,E,\Psi)$, $\mcP=(p_e)_{e\in E}$  and $\{\mu_v\}_{v \in V}$ be same as before.
Firstly, we  prove the existence of $L^q$-spectra of measures in $\{\mu_v^x,\mu_v^y\}_{v\in V}.$
	
	\begin{proof}[Proof of Proposition \ref{thmain4}]
		Regard $\bar{V}=\{\mu_v^x,\mu_v^y\}_{v\in V}$ as a vertex set. Write $v_x=\mu_v^x$, $v_y=\mu_v^y$ for short, i.e. $\bar{V}=\{v_x,v_y\}_{v\in V}$. For each $e\in E$ with $v=i(e),v'=t(e)$,
		if $T_e$ is diagonal, we associate an edge $e_x$ so that   $v_x\stackrel{e_x}{\longrightarrow}  v'_x$ (resp.  $e_y$ so that   $v_y \stackrel{e_y}{\longrightarrow} v'_y$); if $T_e$ is anti-diagonal, we associate an edge $e_x$ so that   $v_x\stackrel{e_x}{\longrightarrow}  v'_y$ (resp.  $e_y$ so that   $v_y \stackrel{e_y}{\longrightarrow} v'_x$).
		Let $\psi_{e_x}=\pi_x(\psi_e)$, $\psi_{e_y}=\pi_y(\psi_e)$ and $p_{e_x}=p_{e_y}=p_e$. Denote $\bar{E}=\{e_x,e_y\}_{e\in E}$, $\bar{\Psi}=\{\psi_{e_x},\psi_{e_y}\}_{e\in E}$ and $\bar{\mathcal{P}}=(p_{e_x},p_{e_y})_{e\in E}.$
		Then  $(\bar{V}, \bar{E},\bar{\Psi})$ becomes a self-similar GIFS, and $\{\mu_v^x,\mu_v^y\}_{v\in V}$ is the unique (but not necessarily strongly connected) graph-directed self-similar measure family associated with $(\bar{V},\bar{E},\bar{\Psi})$ and $\bar{\mathcal{P}}$.

		Consider the  \textit{adjacency matrix} 
		$\mathcal{A}$ associated with $(\bar{V},\bar{E})$, i.e.
		$$
		\mathcal{A}(\bar{v},\bar{v}')=\left\{
		\begin{aligned}
			&1 \quad & & \text{if there exists } \bar{e}\in \bar{E} \text{ such that } \bar{v}\stackrel{\bar{e}}{\to} \bar{v}',
			\\
			&0 & & \text{ otherwise}.
		\end{aligned} \right .
		$$
		
		\vspace{0.2cm}
		
		When $\mathcal{A}$ is \textit{irreducible} (i.e. for $\bar{v},\bar{v}'\in \bar{V}, \mathcal{A}^k(\bar{v},\bar{v}') >0$ for some $k\in \mathbb{N}$), $(\bar{V},\bar{E},\bar{\Psi})$ is a strongly connected self-similar GIFS. By Proposition \ref{prop1}, we know that all $\tau_{\mu_v^x}(q),\tau_{\mu_v^y}(q)$ exist and equal to a common value.
		\vspace{0.2cm}
		
		When $\mathcal{A}$ is not irreducible, pick a pair $\bar{v}',\bar{v}''\in \bar{V}$ such that 
		\begin{equation}\label{s63}
			\bar{v}'\nrightarrow\bar{v}''.
		\end{equation}
		Define 
		$$
		\begin{aligned}
			\bar{V}'&=\{\bar{v}\in \bar{V}: \bar{v}' \to \bar{v}\},
			\quad 
			\bar{V}''&=\{\bar{v}\in \bar{V}: \bar{v}\to \bar{v}''\}.
		\end{aligned}
		$$
		Clearly $\bar{V}'\cap \bar{V}''=\emptyset.$ For each $v\in V$,
		noticing that there exist $k_1,k_2 \in \mathbb{N}$ such that 
		$$
		\mathcal{A}^{k_1}[\{\bar{v}'\},\{v_x,v_y\}] \text{ and } \mathcal{A}^{k_2}[\{v_x,v_y\},\{\bar{v}''\}] \text{ are non-zero matrix},
		$$
		using \eqref{s63}, we have 
		$$
		0=\mathcal{A}^{k_1+k_2}(\bar{v}',\bar{v}'') \geq \mathcal{A}^{k_1}[\{\bar{v}'\},\{v_x,v_y\}] \cdot \mathcal{A}^{k_2}[\{v_x,v_y\},\{\bar{v}''\}].
		$$
		Thus
		\begin{equation}\label{s70}
			\text{ either }\quad v_x \in \bar{V}', v_y\in \bar{V}'' \quad \text{ or }\quad v_x \in \bar{V}'', v_y\in \bar{V}',
		\end{equation}
		which implies that 
		\begin{equation}\label{s69}
			\bar{V}'\cup \bar{V}''=\bar{V},\quad \text{ and } \quad \#\bar{V}'= \#\bar{V}''=\# \bar{V}/2.
		\end{equation}

		For the above $v$, we continue to  consider two cases. 
		\vspace{0.2cm}
		
		\textit{Case 1:  $v_x \in \bar{V}', v_y\in \bar{V}''$}.
		\vspace{0.2cm}
		
		We have 
		$$
		v_x \nrightarrow v_y.
		$$
		Define 
		$$
		\begin{aligned}
			\bar{V}_v'&=\{\bar{v}\in \bar{V}: v_x\to \bar{v}\},\quad 
			\bar{V}_v''&=\{\bar{v}\in \bar{V}: \bar{v}\to v_y\}.
		\end{aligned}
		$$
		We can see that $\bar{V}_v''\cap \bar{V}'= \emptyset$, which by \eqref{s69} immediately implies that $ \bar{V}_v'=\bar{V}'$ and  $\bar{V}_v''=\bar{V}''.$
		Indeed, suppose that $\bar{V}_v''\cap \bar{V}'\neq \emptyset$, then  we have 
		$$
		\bar{v}' \to \tilde{v} \to v_y \to 
		\bar{v}'' \text{ for some } \tilde{v}\in \bar{V}_v''\cap \bar{V}',
		$$
		which contradicts to \eqref{s63}. Also, we have 
		$\bar{v}'\in \bar{V}_v'=\bar{V}' $,  since  if $\bar{v}'\in \bar{V}_v''$, then
		$
		\bar{v}'\to v_y \in \bar{V}',
		$
		which contradicts to $v_y \in \bar{V}''$.
		\vspace{0.2cm}
		
		\textit{Case 2: $v_x\in \bar{V}'',v_y\in \bar{V}'$}.
		\vspace{0.2cm}
		
		Define 
		$$
		\begin{aligned}
			\bar{V}_v'&=\{\bar{v}\in \bar{V}: v_y\to \bar{v}\},\quad 
			\bar{V}_v''&=\{\bar{v}\in \bar{V}: \bar{v} \to v_x\}.
		\end{aligned}
		$$
		By a similar argument as above,
		we also have $\bar{v}'\in \bar{V}_v'=\bar{V}'$ and $\bar{V}_v''=\bar{V}''$.
		
		Thus in both cases, $\mathcal{A}[\bar{V}',\bar{V}']$ is irreducible and $\mathcal{A}[\bar{V}',\bar{V}''] $ is a zero matrix.
		
		Let $\kappa:\bar{V} \to \bar{V}$ be a one-to-one map  defined as  $\kappa(v_x)=v_y$, $\kappa(v_y)=v_x$, for any $v\in V$. By \eqref{s70} and the definition of $(\bar{V},\bar{E})$, we know that $\kappa(\bar{V}')=\bar{V}''$
		and $\mathcal{A}(\tilde{v},\tilde{v}')=\mathcal{A}({\kappa(\tilde{v}), \kappa(\tilde{v}')})$. Thus $\mathcal{A}[\bar{V}'',\bar{V}'']$ is also irreducible and  $\mathcal{A}[\bar{V}'',\bar{V}']$ is a zero matrix.

		Let $A= \bar{V}'$ and $B= \bar{V}''$. Then $A,B$ satisfy \eqref{s65} and  are two strongly connected self-similar measure families. By Proposition \ref{prop1}, \eqref{s66} and \eqref{s67} holds.
		
		Finally, if all $T_e$'s are diagonal. It is easy to see that   $\{\mu_v^x\}_{v\in V}$ and $\{\mu_v^y\}_{v\in v}$ are two strongly connected graph-directed self-similar measure families. So we may choose $A=\{\mu_v^x\}_{v\in V}$ and $B=\{\mu_v^y\}_{v\in V}$.
	\end{proof}
	
		For $w=w_1\cdots w_k\in E^*$,  denote $p_w$ the product $p_{w_1} \cdots p_{w_k}$, $T_w$ the product $T_{w_1} \cdots T_{w_k}$, and  $\psi_w$ the composition $\psi_{w_1} \circ \cdots \circ \psi_{w_k}$.
	Let $$c_w=\big|\pi_x\left(\psi_w([0,1]^2)\right)\big| \quad \text{ and } \quad d_w=\big|\pi_y\left(\psi_w([0,1]^2)\right)\big|
 $$ 
 denote the width and height of the rectangle $\psi_w([0,1]^2)$. Define
	\begin{equation}
		\label{s10}
		\pi_w= \left \{ \begin{aligned}
			&\pi_x \quad & &\text{ if } c_w\geq d_w \text{ and } T_w \text{ is diagonal},\\
			&\pi_y & &\text{ if } c_w <d_w \text{ and } T_w \text{ is diagonal},\\
			&\pi_y & &\text{ if } c_w\geq d_w \text{ and } T_w \text{ is anti-diagonal},\\
			&\pi_x & &\text{ if } c_w<d_w \text{ and } T_w \text{ is anti-diagonal}.
		\end{aligned}\right.
	\end{equation}
        For $q\geq 0,$ define $$\tau_w(q)=\tau_{\pi_w(\mu_{t(w)})}(q).$$
	In other words, $\tau_w(q)$ is  the $L^q$-spectrum of the projection of $\mu_{t(w)}\circ \psi_{w}^{-1}$ onto the longest side of the rectangle $\psi_w([0,1]^2)$, and it always equals to either $\tau_A(q)$ or $\tau_B(q)$ by Proposition \ref{thmain4}.
	
	For $i=1,2$, denote $\alpha_i(T)$ the $i$-th \textit{singular value} of a $2\times 2$ non-singular matrix $T$, i.e. the positive square root of the $i$-th (in decreasing order) eigenvalue of $T^*T$, where $T^*$ is the transpose of $T$. For $w \in E^*$, we write $\alpha_i(w)$ instead of $\alpha_i(T_w)$ for short. Now we are able to give the definition of \textit{modified}  singular value function matrices (differs from the original definition \cite{F88}), inspired by Fraser \cite{F12, F16} dealing with the IFS setting.
	
	\begin{definition}[modified singular value function matrices]
		For $s\in \mathbb{R}$ and $q\geq 0$, define a function  $\varphi^{s,q}:E^*\to (0,\infty)$ by
		$$
		\varphi^{s,q}(w)=p_w^q \alpha_1(w)^{\tau_w(q)} \alpha_2(w)^{s-\tau_w(q)}.
		$$
		For $k\in \mathbb{N}$, define a $\#V \times \# V$ matrix $A_k^{s,q}$ with entries 
		$$
		A_k^{s,q}(v,v')=\sum_{w\in E^k: v \stackrel{w}{\rightarrow} v'} \varphi^{s,q}(w),
		$$ where the empty sum is taken to be 0. Denote $A_0^{s,q}=id$, i.e. the identity matrix, for convention.
		Call $\{A_k^{s,q}\}$ a sequence of  \emph{modified singular value function matrices}.
	\end{definition}
 
	\begin{lemma}\label{lemma7}
	    For $w=w_1 \cdots w_k,w'=w'_1 \cdots w'_l\in E^*$ with $t(w)=i(w'),$
     \begin{equation}\nonumber
        \begin{aligned}
				&\tau_{x,w_k}(q)=\tau_{x,w'_l}(q), \quad \tau_{y,w_k}(q)=\tau_{y,w'_l}(q) \quad & &\text{ if  } T_{w'} \text{ is diagonal,}\\
				&\tau_{x,w_k}(q)=\tau_{y,w'_l}(q), \quad \tau_{y,w_k}(q)=\tau_{x,w'_l}(q) \quad & & \text{ if } T_{w'} \text{ is anti-diagonal.}
			\end{aligned}
        \end{equation}
	\end{lemma}
 \begin{proof}
    It suffices to  prove that for $e,e'\in E$ with $t(e)=i(e')$,
        \begin{equation}\label{s81}
        \begin{aligned}
				&\tau_{x,e}(q)=\tau_{x,e'}(q), \quad \tau_{y,e}(q)=\tau_{y,e'}(q) \quad & &\text{ if  } T_{e'} \text{ is diagonal,}\\
				&\tau_{x,e}(q)=\tau_{y,e'}(q), \quad \tau_{y,e}(q)=\tau_{x,e'}(q) \quad & & \text{ if } T_{e'} \text{ is anti-diagonal.}
			\end{aligned}
        \end{equation}
        It follows from the proof of Proposition \ref{thmain4}, we know that if $T_{e'}$ is diagonal, either $$\mu_{t(e)}^x ,\mu_{t(e')}^x\in A,\quad  \mu_{t(e)}^y ,\mu_{t(e')}^y\in B,$$
		or 
		$$\mu_{t(e)}^x ,\mu_{t(e')}^x\in B,\quad  \mu_{t(e)}^y ,\mu_{t(e')}^y\in A,$$
		unless $\tau_A(q)=\tau_B(q)$. 
		By Proposition \ref{thmain4}, \eqref{s81} holds in this case.
		Also, if $T_{e'}$ is anti-diagonal, using a same argument, we still have \eqref{s81}.
 \end{proof}
	
	\begin{lemma} \label{lem23}
		Let $s\in \mathbb{R},q\geq 0$.
		\begin{itemize}
			
			\item[(a).] For $w, w'\in E^*$ with $t(w)=i(w')$, we have 
			\begin{itemize}
				\item[(a1).] if $s<t(q)$, $\varphi^{s,q}(ww')\leq \varphi^{s,q}(w)\varphi^{s,q}(w')$,
				\item[(a2).]   if $s=t(q)$, $\varphi^{s,q}(ww')= \varphi^{s,q}(w)\varphi^{s,q}(w')$,
				\item[(a3).] if $s>t(q)$, $\varphi^{s,q}(ww')\geq \varphi^{s,q}(w)\varphi^{s,q}(w')$.
			\end{itemize}
			\item[(b).] For $k,l \in \mathbb{N}$, we have 
			\begin{itemize}
				\item[(b1).]  if $s\leq t(q)$,  $\Vert A_{k+l}^{s,q} \Vert \leq \Vert A_k^{s,q} \Vert \cdot \Vert A_l^{s,q} \Vert$,
				\item[(b2).] if $s>  t(q)$,
				 $\Vert A_{k+1}^{s,q}\Vert\gtrsim_{s,q}\Vert A_k^{s,q} \Vert$ and there exists  $J\in \mathbb{N}$ independent of $k,l$ such that $\max_{0\leq j\leq J}\Vert A_{k+j+l}^{s,q} \Vert \gtrsim_{s,q} \Vert A_k^{s,q} \Vert \cdot \Vert A_l^{s,q} \Vert$.
			\end{itemize}
			
		\end{itemize}
	\end{lemma}
	\begin{proof}
		(a1). Write $w=w_1 \cdots w_k$ and $w'=w'_1 \cdots w_l'$. By Proposition \ref{thmain4}, we always have 
$$
\tau_{x,e}(q)+\tau_{y,e}(q)=t(q), \quad \text{ for all } e\in E.
$$
        We only prove the case that $c_w\geq d_w$ while the  case $c_w<d_w$ is similar. We consider two possible subcases.\vspace{0.2cm}

        \textit{Case 1: $c_{ww'}=c_w c_{w'}, d_{ww'}=d_{w}d_{w'}$.}
        \vspace{0.2cm}
        
        Note that in this case $T_w$ is diagonal, so that  $\tau_w(q)=\tau_{x,w_k}(q)$. If $c_{w'}\geq d_{w'}$ and $c_{ww'}\geq d_{ww'}$, using Lemma \ref{lemma7}, we always have 
        $$
        \tau_{w'}(q)=\tau_{ww'}(q)=\tau_{w}(q),
        $$
        by checking $T_{w'}$ is diagonal or not separately.
        Thus
        \begin{equation}\nonumber
        \frac{\varphi^{s,q}(ww')}{\varphi^{s,q}(w)\varphi^{s,q}(w')}=\frac{(c_w c_{w'})^{\tau_w(q)} (d_w d_{w'})^{s-\tau_w(q)}}{c_w^{\tau_w(q)} d_w^{s-\tau_w(q)} c_{w'}^{\tau_w(q)} d_{w'}^{s-\tau_w(q)}} =1.
        \end{equation}
 If $c_{w'}<d_{w'}$ and $c_{ww'}\geq d_{ww'}$, considering similarly as above, we always have 
        $$
        \tau_{w'}(q)+\tau_{w}(q)=t(q) \quad \text{ and }\quad \tau_{ww'}(q)=\tau_{w}(q).
        $$
        Thus 
        $$
        \frac{\varphi^{s,q}(ww')}{\varphi^{s,q}(w)\varphi^{s,q}(w')}=\frac{(c_w c_{w'})^{\tau_w(q)} (d_w d_{w'})^{s-\tau_w(q)}}{c_w^{\tau_w(q)} d_w^{s-\tau_w(q)} d_{w'}^{t(q)-\tau_w(q)} c_{w'}^{s-t(q)+\tau_w(q)}}=\left(\frac{d_{w'}}{c_{w'}}\right)^{s-t(q)}\leq 1.
        $$
Otherwise,  $c_{w'}<d_{w'}$ and $c_{ww'}<d_{ww'}$, and similarly,
        $$
        \tau_{w'}(q)+\tau_{w}(q)=t(q) \quad \text{ and }\quad \tau_{ww'}(q)=\tau_{w'}(q)
        $$
        always holds. Thus 
        $$
        \frac{\varphi^{s,q}(ww')}{\varphi^{s,q}(w)\varphi^{s,q}(w')}=\frac{(d_w d_{w'})^{t(q)-\tau_w(q)} (c_w c_{w'})^{s-t(q)+\tau_w(q)}}{c_w^{\tau_w(q)} d_w^{s-\tau_w(q)} d_{w'}^{t(q)-\tau_w(q)} c_{w'}^{s-t(q)+\tau_w(q)}}=\left(\frac{c_{w}}{d_{w}}\right)^{s-t(q)}\leq 1.
        $$

In summary, (a1) holds in this case.\vspace{0.2cm}

        \textit{Case 2: $c_{ww'}=c_w d_{w'}, d_{ww'}=d_w c_{w'}.$}
        \vspace{0.2cm}
        
        In this case, $T_w$ is anti-diagonal, so $\tau_{w}(q)=\tau_{y,w_k}(q).$ If $d_{w'}\geq c_{w'}$ and $c_{ww'}\geq d_{ww'}$, we always have 
        $$
        \tau_{w'}(q)=\tau_{ww'}(q)=\tau_{w}(q),
        $$
        so
        \begin{equation}\nonumber
        \frac{\varphi^{s,q}(ww')}{\varphi^{s,q}(w)\varphi^{s,q}(w')}=\frac{(c_w d_{w'})^{\tau_w(q)} (d_w c_{w'})^{s-\tau_w(q)}}{c_w^{\tau_w(q)} d_w^{s-\tau_w(q)} d_{w'}^{\tau_w(q)} c_{w'}^{s-\tau_w(q)}} =1.
        \end{equation}
          If $d_{w'}<c_{w'}$ and $c_{ww'}\geq d_{ww'}$, we always have 
        $$
        \tau_{w'}(q)+\tau_{w}(q)=t(q) \quad \text{ and } \quad \tau_{ww'}(q)=\tau_{w}(q),
        $$
        so
        \begin{equation}\nonumber
        \frac{\varphi^{s,q}(ww')}{\varphi^{s,q}(w)\varphi^{s,q}(w')}=\frac{(c_w d_{w'})^{\tau_w(q)} (d_w c_{w'})^{s-\tau_w(q)}}{c_w^{\tau_w(q)} d_w^{s-\tau_w(q)} c_{w'}^{t(q)-\tau_w(q)} d_{w'}^{s-t(q)+\tau_w(q)}} =\left(\frac{c_{w'}}{d_{w'}}\right)^{s-t(q)}\leq 1.
        \end{equation}
Otherwise, 
        $d_{w'}<c_{w'}$ and $c_{ww'}<d_{ww'}$, and we always have
        $$
        \tau_{w'}(q)+\tau_{w}(q)=t(q) \quad \text{ and } \quad \tau_{ww'}(q)=\tau_{w'}(q),
        $$
        so
        \begin{equation}\nonumber
        \frac{\varphi^{s,q}(ww')}{\varphi^{s,q}(w)\varphi^{s,q}(w')}=\frac{(d_w c_{w'})^{t(q)-\tau_w(q)} (c_w d_{w'})^{s-t(q)+\tau_w(q)}}{c_w^{\tau_w(q)} d_w^{s-\tau_w(q)} c_{w'}^{t(q)-\tau_w(q)} d_{w'}^{s-t(q)+\tau_w(q)}} =\left(\frac{c_{w}}{d_{w}}\right)^{s-t(q)}\leq 1.
        \end{equation}
        
        So in summary, (a1) also holds in this case.

  \vspace{0.2cm}

    The proofs of (a2) and (a3)  follow by using a similar argument as above.
  \vspace{0.2cm}

		(b1). If $s\leq t(q)$, then
		\begin{align}
			\Vert A_{k+l}^{s,q} \Vert &= \sum_{v,v'\in V} A_{k+l}^{s,q}(v,v')
			=\sum_{v,v'\in V} \sum_{w\in E^{k+l}: v\stackrel{w}{\rightarrow}v'} \varphi^{s,q}(w)\nonumber\\
			&\leq \sum_{v,v'\in V} \sum_{v''\in V} \sum_{w\in E^k:v\stackrel{w}{\rightarrow} v''} \sum_{w'\in E^l:v''\stackrel{w'}{\rightarrow} v'} \varphi^{s,q}(w) \varphi^{s,q}(w')
			\nonumber \quad \text{ by (a1) and (a2)}\\
			&=\sum_{v,v'\in V} \sum_{v''\in V} A_k^{s,q}(v,v'') A_l^{s,q}(v'',v') =\sum_{v,v'\in V}( A_k^{s,q}\cdot A_l^{s,q} )(v,v')
			\nonumber\\
			&= \Vert A_k^{s,q} \cdot A_l^{s,q} \Vert \leq \Vert A_k^{s,q} \Vert \Vert A_l^{s,q} \Vert.\nonumber
		\end{align} 
		
		(b2). If $s> t(q)$, for $k\geq 0$, we have $A_{k+1}^{s,q}\geq A_k^{s,q} A^{s,q}_1$ by a same argument as above and using $(a3)$. Denote $C=\min_{v\in V} \sum_{v'\in V} A_1^{s,q}(v,v')>0$. It is direct to see that  $\Vert A_k^{s,q} A_1^{s,q} \Vert \geq C \Vert A_k^{s,q} \Vert $ which gives the first part of $(b2)$.
		
		Supposed the second part of $(b2)$ is not true, i.e. for any $\epsilon >0$ and  $J\in \mathbb{N}$, there exist $k,l\in \mathbb{N}$ such that $ \Vert A_{k+j+l}^{s,q} \Vert \leq \epsilon \Vert A_k^{s,q} \Vert \Vert A_l^{s,q} \Vert$ for all $0\leq j\leq J$. Let $u=(1,\cdots, 1)$ be a row vector in $\mathbb{R}^{\#V}$. Noticing that  $A_{k+j+l}^{s,q} \geq A_k^{s,q} \cdot A_j^{s,q} \cdot A_l^{s,q}$, we have
		$$
		u A_k^{s,q} A_j^{s,q} A_l^{s,q} u^* \leq \Vert A_{k+j+l}^{s,q} \Vert  \leq \epsilon \Vert A_k^{s,q} \Vert \Vert A_l^{s,q} \Vert.
		$$
		Define two non-negative unit row vectors ${u'}=\frac{u\cdot A_k^{s,q}}{\Vert A_k^{s,q}\Vert }$ and ${u''}=\frac{ (A_l^{s,q} u^*)^*}{\Vert A_l^{s,q} \Vert}$, then 
		$$
		{u'} A_j^{s,q} {u''}^* \leq \epsilon, \quad  \text{ for all } 0\leq j\leq J .
		$$
		Take two sequence $\{\epsilon_n\}\to 0$ and $\{J_n\} \to \infty$. Then for each $n$, there exist two non-negative unit row vector $u_n',u_n''$ such that 
		$$
		u_n' A_j^{s,q} {u_n''}^* \leq \epsilon_n,\quad  \text{ for all } 0\leq j\leq J_n.
		$$
		Let $(u',u'')$ be a limit point of $\{(u_n',u_n'')\}_{n\in \mathbb{N}}$, then we have 
		$$
		u' A_j^{s,q} u''^*=0, \quad  \text{ for all } j\geq 0.
		$$
		Noticing that $u',u''$ are two non-negative  unit row vectors and all $A_j^{s,q} $ are  non-negative matrices, we have that there exist $v,v'\in V$ such that $A_j^{s,q}(v,v')=0$ for all $j\geq 0$. This contradicts the strongly connectivity of the directed graph $(V,E)$.
	\end{proof}
	
	\begin{remark}
		The property  (b2) was summarised by Feng \cite{F09} which was  fundamental for random matrix product and multifractal analysis, see e.g. \cite{KR14,FK11}.
	\end{remark}
	
	\begin{lemma}[pressure function]\label{lemd}
		Define a function $P:\mathbb{R}\times [0,\infty) \to [0,\infty)$ by 
		\begin{equation}\label{defpressure}
			P(s,q)=\lim_{k\to \infty} \Vert A_k^{s,q} \Vert^{1/k}.
		\end{equation}
		The definition is well-defined. Call $P$ a \emph{pressure function} associated with $(V,E,\Psi)$  and $\mathcal{P}$.
	\end{lemma}
	
	\begin{proof}
		It  suffices to prove the existence of limit in \eqref{defpressure}.
		
		For $q\geq 0$, $s\leq t(q)$, the limit exists by Lemma \ref{lem23}-(b1) and the standard property of submultiplicative sequences.
		
		For  $s> t(q)$, by Lemma \ref{lem23}-(b2), there exist $C>0$ (depending on $s,q$) and $J\in \mathbb{N}$ such that 
		\begin{equation}\nonumber
			\begin{aligned}
				\Vert A_{k+1}^{s,q} \Vert  &\geq C \Vert A_k^{s,q} \Vert, &\quad & \text{ for all } k\geq 0,\nonumber \\
				\max_{0\leq j\leq J} \Vert A_{k+j+l}^{s,q} \Vert &\geq C \Vert A_k^{s,q} \Vert \Vert A_l^{s,q} \Vert, &\quad & \text{ for all } k,l\geq 0.\nonumber
			\end{aligned}
		\end{equation}
		Therefore, there exists $C'>0$ such that $\Vert A_{k}^{s,q} \Vert \Vert A_l^{s,q} \Vert \leq C^{-1}\sum_{0\leq j\leq J}\Vert A_{k+j+l}^{s,q} \Vert \leq C' \Vert A_{k+J+l}^{s,q}\Vert $.
		This implies that the sequence $\{\frac{1}{C'} \Vert A_{k-J}^{s,q} \Vert\}_{k\geq J}$  is supermultiplicative, so the limit \eqref{defpressure} exists.
	\end{proof}
	Write 
	\begin{equation}\label{s16}
		\begin{aligned}
			\alpha_*&=\min \{\alpha_2(e):e\in E\} =\min \{a_e,b_e:e\in E\}, \\ 
			\alpha^*&=\max \{\alpha_1(e):e\in E\}=\max \{a_e,b_e:e\in E\},\\
			p_*&=\min \{p_e:e\in E\},  \\ 
			p^*&=\max \{p_e:e\in E\}. 
		\end{aligned}
	\end{equation}
	Recall that the $L^q$-spectra of a measure is Lipschitz continuous on $[\lambda,\infty)$ for all $\lambda>0$. Let $L_\lambda$ be the larger of the two Lipschitz constants corresponding to $\tau_A$ and $\tau_B$ on $[\lambda,\infty)$.
	
	\begin{lemma}\label{propertyP}
		For $t,r\in \mathbb{R}$ and $\lambda >0$, define 
		$$
		U(t,r,\lambda)=\min \{\alpha_*^tp_*^r, \alpha_*^tp^{*r}, \alpha^{*t}p_*^r, \alpha^{*t}p^{*r}\}\left(\frac{\alpha^*}{\alpha_*}\right)^{\min \{-L_\lambda r,0\}}
		$$
		and 
		$$
		V(t,r,\lambda)=\max \{\alpha_*^tp_*^r, \alpha_*^tp^{*r}, \alpha^{*t}p_*^r, \alpha^{*t}p^{*r}\}\left(\frac{\alpha^*}{\alpha_*}\right)^{\max \{-L_\lambda r,0\}}.
		$$
		Then for all $s,t\in \mathbb{R}$, $\lambda>0$, $q\geq \lambda$ and $r\geq \lambda-q$,
		$$
		U(t,r,\lambda) P(s,q) \leq P(s+t,q+r) \leq V(t,r,\lambda) P(s,q).
		$$
		and for $s,t\in \mathbb{R}$,
		$$
		\min\{\alpha_*^t,\alpha^{*t}\} P(s,0)\leq P(s+t,0) \leq \max\{ \alpha_*^t,\alpha^{*t} \} P(s,0).
		$$
		Also, for all $s\in \mathbb{R}$ and $q\geq 0$,
		$$
		P(s,q)\leq p^{*q}P(s,0).
		$$
		Consequently, we have
		
		(a). $P$ is continuous on $\mathbb{R}\times (0,\infty)$ and on $\mathbb{R}\times \{0\}$,
		
		(b). $P$ is strictly decreasing in $s\in \mathbb{R}$,
		
		(c).  for each $q\geq 0$, there exists a unique $s\in \mathbb{R}$ such that $P(s,q)=1$.
	\end{lemma}
	\begin{proof}
		This is essentially the same as \cite[Lemma 2.3]{F16}.
	\end{proof}
	\begin{remark}\label{re1}
		For $q\geq 0$, we refer $\gamma(q)$ to be the unique $s\in \mathbb{R}$ satisfying $P(s,q)=1$. In Theorem \ref{thmain1}, we will prove that $\tau_{\mu_v}(q)=\gamma(q)$ for all $v\in V$. 
	\end{remark}

	\section{Closed forms in non-rotational setting}\label{sec4}
	
	In this section, we mainly prove Theorem \ref{thmain2} and Corollaries \ref{cor2}, \ref{cor1}. We postpone the proof of Theorem \ref{thmain1} to  Section \ref{sec3}, and assume that it is true in advance. We always let  $(V,E,\Psi)$ be a strongly connected planar box-like self-affine  GIFS with all $T_e$'s diagonal. Let $\mathcal{P}$, $\{\mu_v\}_{v\in V}$ be the associated positive vector and measures as before. Throughout this section,  we always fix a $q\geq 0$, so when define new variables, we may omit $q$.
	
	\subsection{Notations and lemmas}

	For $w=w_1 \cdots w_k\in E^*$, write $a_w =a_{w_1} \cdots a_{w_k}$ and $b_w =b_{w_1} \cdots b_{w_k}$. Recalling the definition of $c_w$ and $d_w$ in Section \ref{sec2}, we have $a_w =c_w$, $b_w=d_w$ are the width and height of the rectangle $\psi_w([0,1]^2)$.
	Therefore, by Proposition \ref{thmain4}, we have
	\begin{equation}\label{e5}
		\tau_w(q)= \left \{ \begin{aligned}
			&\tau_A(q) \quad &\text{ if } a_w\geq b_w,\\
			&\tau_B(q) &\text{ if } a_w <b_w.
		\end{aligned}\right.
	\end{equation}
	As announced in \eqref{s14}, we introduce a $\#E \times \#E$ function matrix $F^{(q)}_{x,y}$ with entries defined by
	\begin{equation}\nonumber
		F^{(q)}_{x,y}(e,e') =\left \{ \begin{aligned}
			&p_{e'}^q a_{e'}^x b_{e'}^y &\quad& \text{ if } t(e)=i(e'),\\
			&0 &\quad& \text{ otherwise},
		\end{aligned} \right.
	\end{equation}
	and write $\rho(F^{(q)}_{x,y})$  the spectra radius of $F^{(q)}_{x,y}$.
	
	\begin{lemma}\label{lemma41} The function $\rho(F^{(q)}_{x,y})$  is continuous in $x,y\in \mathbb{R}$. For fixed $y\in \mathbb{R}$, $\rho(F^{(q)}_{x,y})$ is strictly decreasing in $x\in \mathbb{R}$, and there exists a unique $x\in \mathbb{R}$ such that $\rho(F^{(q)}_{x,y})=1$. This is also true for $\rho(F^{(q)}_{x,y})$ as a function of $y\in \mathbb{R}$ for fixed $x \in \mathbb{R}$.
	\end{lemma}
	\begin{proof}
		For any $\epsilon,\eta \in \mathbb{R}$, let 
		$$
		U(\epsilon,\eta) =\min \{\alpha_*^{\epsilon+\eta}, \alpha_*^\epsilon \alpha^{*\eta}, \alpha^{*\epsilon} \alpha_*^\eta, \alpha^{*(\epsilon+\eta)}\}
		$$
		and 
		$$
		V(\epsilon,\eta) =\max \{\alpha_*^{\epsilon+\eta}, \alpha_*^\epsilon \alpha^{*\eta}, \alpha^{*\epsilon} \alpha_*^\eta, \alpha^{*(\epsilon+\eta)}\},
		$$
		where $\alpha_*, \alpha^*$ defined by \eqref{s16}.
		Then
		$$
		U(\epsilon,\eta)	p_{e}^q a_{e}^x b_{e}^y \leq p_{e}^q a_{e}^{x+\epsilon} b_{e}^{y+\eta} \leq V(\epsilon,\eta)  p_{e}^q a_{e}^x b_{e}^y, \quad  \text{ for all } {e} \in E,
		$$
		which yields 
		$$
		U^k(\epsilon,\eta) 	\big\Vert \left(F^{(q)}_{x,y}\right)^k \Vert  \leq \big\Vert \left(F^{(q)}_{x+\epsilon,y+\eta}\right)^k \big\Vert  \leq V^k(\epsilon,\eta)  \big\Vert  \left(F^{(q)}_{x,y}\right)^k \big\Vert,\quad  \text{ for all } k\in \mathbb{N}.
		$$
		By Gelfand formula, we have 
		\begin{equation}\label{s15}
			U(\epsilon,\eta) \rho(F^{(q)}_{x,y})\leq \rho(F^{(q)}_{x+\epsilon,y+\eta}) \leq V(\epsilon,\eta)  \rho(F^{(q)}_{x,y}),
		\end{equation}
		which gives the first part of this lemma. For fixed $y\in \mathbb{R}$, letting $\eta=0$, still using \eqref{s15}, $\rho(F^{(q)}_{x,y})$ tends to $0$ as $x\to +\infty$, and tends to $+\infty$ as $x \to -\infty$. Therefore, there exists a unique $x\in \mathbb{R}$ such that $\rho(F^{(q)}_{x,y})=1$.
		
	\end{proof}

	By Lemma \ref{lemma41}, we know that  $\gamma_A(q),\gamma_B(q)$ defined by \eqref{t7} \eqref{t8} i.e. 
	\begin{equation}\nonumber
		\rho(F^{(q)}_{\tau_A(q),\gamma_A(q)-\tau_A(q)})=1, \quad
		\rho(F^{(q)}_{\gamma_B(q)-\tau_B(q),\tau_B(q)})=1,
	\end{equation}
	are well-defined. 
	\begin{lemma}\label{lemma42}
		Either 
		\eqref{e1} or \eqref{e2} holds.
	\end{lemma}
	\begin{proof}
		Note that if $\gamma_A(q) \leq t(q)$, then by Lemma \ref{lemma41}, $$
        \rho(F^{(q)}_{\tau_A(q),\tau_B(q)}) =\rho (F^{(q)}_{\tau_A(q),t(q)-\tau_A(q)}) \leq 1,$$
        and so $\gamma_B(q)-\tau_B(q)\leq \tau_A(q)$. So in this case $\max \{\gamma_A(q),\gamma_B(q)\}	\leq t(q)$. By a similar argument, the case
		$\gamma_A(q)\geq t(q)$ could also imply $\gamma_B(q)\geq t(q)$.
	\end{proof}

	By Lemma \ref{lemma41}, we may define a function $y(x):\mathbb{R}\to \mathbb{R}$ by $$\rho(F^{(q)}_{x,y(x)})=1.$$  
	\begin{lemma}\label{l1}
		\begin{equation}
			\begin{aligned}\label{e6}
				\\\big\{(x,y):\rho(F^{(q)}_{x,y})&=1, \min \{ \tau_A(q),\gamma_B(q)-\tau_B(q)\}\leq x\leq  \max \{\tau_A(q),\gamma_B(q)-\tau_B(q)\big\}\\
				=\big\{(x,y):\rho(F^{(q)}_{x,y})&=1,\min \{\tau_B(q),\gamma_A(q)-\tau_A(q)\}\leq y\leq \max \{ \tau_B(q),\gamma_A(q)-\tau_A(q)\big\}.
			\end{aligned}
		\end{equation}
	\end{lemma}
	\begin{proof}
		By Lemma \ref{lemma42}, it suffices to prove \eqref{e6} when either \eqref{e1} or \eqref{e2} holds. We only prove the case that \eqref{e1} holds, and the other case can be achieved by a same argument. 
		If \eqref{e1} holds, $\gamma_B(q)-\tau_B(q)\leq \tau_A(q)$.
		For $x \in [\gamma_B(q)-\tau_B(q), \tau_A(q)]$, it follows that $\rho(F^{(q)}_{\tau_A(q), y(x)}) \leq  1$ by Lemma \ref{lemma41}, and  still by Lemma \ref{lemma41}, $y(x) \geq   \gamma_A(q)-\tau_A(q)$.
		Similarly, we also have $y(x) \leq \tau_B(q)$. Thus $y(x)\in [\gamma_A(q)-\tau_A(q),\tau_B(q)]$. Conversely, for $y\in [\gamma_A(q)-\tau_A(q),\tau_B(q)]$, a same arugment as above yields there exists a unique $x\in [\gamma_B(q)-\tau_B(q),\tau_A(q)]$ with $\rho(F^{(q)}_{x,y})=1$. This gives \eqref{e6}.
	\end{proof}
	
	Noticing that $F^{(q)}_{x,y}$ is irreducible  by the strongly connectivity of $(V,E)$. By Perron-Frobenius Theorem, there exists a unique positive unit column vector  $u(x)=(u_e(x))_{e\in E}$ satisfying 
	$$
	F^{(q)}_{x,y(x)} u(x)=u(x),
	$$
	say  the \textit{right Perron vector} of $F^{(q)}_{x,y(x)}$.
	Similarly, there exists a unique positive left (row) eigenvector $v(x)=(v_e(x))_{e\in E}$ of $F^{(q)}_{x,y(x)}$  satisfying  
	\begin{equation}\label{s19}
		\sum_{e\in E} v_e(x) u_e(x)=1,
	\end{equation}
	say the \textit{left Perron vector} of $F^{(q)}_{x,y}$.
	Define a $\# E\times \#E $ matrix $\tilde F_x$ with entries 
	$$
	\tilde{F}_x(e,e')=F^{(q)}_{x,y(x)}(e,e') \frac{u_{e'}(x)}{u_e(x)}.
	$$ 
	Obviously, $\tilde{F}_x$ is a stochastic matrix, i.e all row sums equal to $1$. Let 
	\begin{equation}\label{s71}
	f^{(q)}(x):=(f^{(q)}_e(x))_{e\in E}=(v_e(x) u_e(x))_{e\in E}
\end{equation}
	be a positive probability row vector by \eqref{s19}. Then 
	\begin{equation}\label{s22}
		f^{(q)}(x) \tilde{F}_x =f^{(q)}(x).
	\end{equation}
	\begin{lemma}\label{lemma43}
		The function $y(x)$ is continuous, decreasing, and  $f^{(q)}(x)$ is continuous in $x\in \mathbb{R}$.
	\end{lemma}
	\begin{proof}
		Assuming that $y(x)$ is not continuous at $x\in \mathbb{R}$, we may find a sequence $\{x_n\} \to x$ such that $y(x_n) \to y' \neq y(x)$. By Lemma \ref{lemma41}, $1 =\rho(F^{(q)}_{x_n,y(x_n)}) \to \rho(F^{(q)}_{x,y'})\neq 1$, a contradiction. Also, by Lemma \ref{lemma41}, $y(x)$ is decreasing.
		
		In order to prove the continuity of $f^{(q)}(x)$, it  suffices to prove $u(x),v(x)$ are continuous. Assume that $u(x)$ is not continuous at $x$. Pick a sequence $\{x_n\}\to x$ satisfying $u(x_n) \to u'\neq u(x)$. Then $u'$ is a non-negative unit vector. By the proof of Lemma \ref{lemma41}, we have 
		$$
		u(x_n)-u(x)=F^{(q)}_{x_n,y(x_n)} u({x_n}) - F^{(q)}_{x,y(x)}u(x) \geq  F^{(q)}_{x,y(x)} \big(U(x_n-x,y(x_n)-y(x)) u(x_n)-u(x)\big)
		$$
		which implies that 
		\begin{equation}\label{s17}
			u'-u(x) \geq F^{(q)}_{x,y(x)}(u'-u(x))
		\end{equation}
		by letting $n\to \infty$.
		Similarly, we have 
		\begin{equation}\label{s18}
			u'-u(x) \leq F^{(q)}_{x,y(x)} (u'-u(x)).
		\end{equation}
		Combining \eqref{s17} and \eqref{s18}, $u'-u(x)$ is a  right eigenvector. However, $u'\neq u(x)$, a contradiction arised by the uniqueness of the right Perron vector of $F^{(q)}_{x,y(x)}$. Similarly, $v(x)$ is continuous.
	\end{proof}
	
	\vspace{0.2cm}

	Before proceeding, we recall some knowledge about \textit{Markov Chains}. Let $X=\{X_i\}_{i\geq 1}$ be a \textit{Markov chain} on a \textit{finite state space} $S$. Suppose $P$ is a \textit{transition probability matrix associated} with $X$, i.e.
	$$
	P(s,s')=\mathbb{P}(X_{i+1}=s'|X_{i}=s), \quad \text{for all } s,s'\in S, i\geq 1.
	$$
	Note that for any $i\geq 1$, $P^i$ gives the $i$-step transition probabilities of the chain $X$, i.e. $P^i(s,s')=\mathbb{P}(X_{i+1}=s'| X_1=s)$, where $P^i$ denotes the $i$-th power of $P$. Also, the matrix $P$ naturally induces an edge set $\Gamma:=\{(s,s'): P(s,s')>0\}$ so that $(S,\Gamma)$ becomes a finite directed graph.

	The marginal distribution of $X_1$ is called the \textit{initial distribution} of $X$. An initial distribution $\lambda=(\lambda_s)_{s\in S}$, together with the transition matrix $P$,  determines a joint distribution $\mathbb P_\lambda$ of the process $X$ by $$\mathbb{P}_\lambda(X_1=s_1,X_2=s_2,\cdots ,X_i=s_i)=\lambda_{s_1}P(s_1,s_2)\cdots P(s_{i-1},s_i).$$
	
	Say the Markov chain $X$ is \textit{irreducible} if  for any $s,s'\in S$, there exists $i\geq 1$ such that $P^i(s,s')>0$. Clearly, $X$ is irreducible if and only if $(S,\Gamma)$ is strongly connected.
	
	Let $\pi=(\pi_s)_{s\in S}$   be the unique \textit{invariant distribution} associated with $X$, i.e. 
	$$
	\pi P=\pi.
	$$
	For any function $h:S \to \mathbb{R}$, write $\pi(h):=\sum_{s\in S}\pi_s h(s)$.
	It is known that an irreducible finite Markov chain with an invariant distribution $\pi$ satisfies the following \textit{central limit theorem}. 
	
	\begin{proposition}\label{prop2}
		For any irreducible  Markov chain $X=\{X_i\}_{i \geq 1}$ with finite state space $S$,  invariant distribution $\pi$, for  any $h:S \to \mathbb{R}$ satisfying $\pi(h^2)<+\infty$, the central limit theorem holds, i.e. there exists $\sigma_h^2<+\infty$, such that  for any initial distribution $\lambda$,
		\begin{equation}\nonumber
			\sqrt{k} \left (\sum_{i=1}^k \frac{h(X_i)}{k} -\pi(h) \right ) \stackrel{\mathbb{P}_\lambda}{\longrightarrow} N(0,\sigma_h^2), \quad  \text{ as } k \to +\infty.
		\end{equation}
	\end{proposition}
	
	See  \cite[Corollary 4.2(ii)]{Cog} or \cite[Theorem 17.0.1]{MT93} for proofs of the above proposition for general Markov processes which are uniform ergodic. Note that an irreducible and aperiodic finite Markov chain is always uniform ergodic \cite[Theorem 16.0.2]{MT93}. In particular, for irreducible finite Markov chains, the proposition still holds without any assumption of aperiodicity, see  \cite[Theorem 23, Proposition 30]{RR04} and the remarks thereafter. Please refer to \cite{MT93} for any unexplained terminologies and details.\vspace{0.2cm}

	In the following, for $x\in \mathbb{R}$, we regard  $E$  as a finite states space,  $X=\{X_i\}_{i\geq 1}$ as  a Markov chain associated with the transition probability matrix $\tilde{F}_x$. Then,  $X$ is irreducible since $(V,E)$ is strongly connected. Also, by \eqref{s22}, $f^{(q)}(x)$ is an invariant distribution of $X$. 
	
	For $w\in E^*$, $e\in E$, denote $\#(w,e) :=\# \{i:w_i=e, 1\leq i\leq |w|\}$ the number of times $e$ appears in $w$. The following lemma is an immediate consequence of Proposition \ref{prop2}.	
	
	\begin{lemma} \label{leclt}
		For fixed $x\in \mathbb{R}$,
		for any initial distribution $\lambda=(\lambda_e)_{e\in E}$, for $\epsilon>0$, there exists $0<C<1$ independent of $\epsilon$ such that $\sum_{w=w_1\cdots w_k \in \mathcal{B}_{x,k}(\epsilon) } \lambda_{w_1} \tilde{F}_x(w_1,w_2) \cdots \tilde{F}_x(w_{k-1},w_k) <C$ for all large enough $k$, where 
		$$	 		\mathcal{B}_{x,k}(\epsilon):=\left \{  w\in E^k: \sum_{e\in E}\left| \frac{\#(w,e)}{k}-f^{(q)}_e(x)\right|>\epsilon \right \}.
		$$		
	\end{lemma}
	\begin{proof}
		For $e\in E$, let $h_e$ be the characteristic function of $\{e\}$ in $E$. Then for $k\geq 1$, $w=w_1 \cdots w_k\in E^k$, $\sum_{i=1}^k h_e(w_i)=\#(w,e)$, and $\pi(h_e)=f_e^{(q)}(x)$. For $\epsilon>0$,  denote 
		$$
		\mathcal{B}_{x,k,e}(\epsilon):=\left \{  w\in E^k: \left| \frac{\#(w,e)}{k}-f^{(q)}_e(x)\right|>\frac{\epsilon}{\#E}  \right \}.
		$$
        For some large $C'>0$, by using Proposition \ref{prop2}, we see that  for large enough $k$, $$\mathcal{B}_{x,k,e}(\epsilon)\subseteq \left \{  w\in E^k: \left| \frac{\#(w,e)}{k}-f^{(q)}_e(x)\right|>\frac{C'}{\sqrt{k}\#E}  \right \},$$ and 
		$$
		\sum_{w=w_1 \cdots w_k\in {\mathcal{B}}_{x,k,e}(\epsilon)} \mathbb{P}_\lambda(X_1=w_1,\cdots, X_k=w_k)=\sum_{w\in {\mathcal{B}}_{x,k,e}(\epsilon)}\lambda_{w_1} \tilde{F}_x(w_1,w_2) \cdots \tilde{F}_x(w_{k-1},w_k)<C'',
		$$
  for some $C''< \frac{1}{\#E}.$
		Since $\mathcal{B}_{x,k}(\epsilon)\subseteq \bigcup_{e\in E}\mathcal{B}_{x,k,e}(\epsilon)$, the lemma follows by taking  $C=(\#E)C''<1$.
	\end{proof}

 \subsection{Proofs of Theorem \ref{thmain2} and Corollaries \ref{cor2}, \ref{cor1}}
	Now we return to the proofs of the main results in this section.
	\begin{proof}[Proof of Theorem \ref{thmain2}] 
		Due to Lemma \ref{lemma42}, it suffices to prove (a) and  (b). We divide the proof into two parts.
		\vspace{0.2cm}
		
		\textit{Part I: When \eqref{e1} holds.}  \vspace{0.2cm}
		
		Recalling the definition of matrix $A_k^{s,q}$ in Section \ref{sec2}, and using \eqref{e5},  
		we have 
		\begin{equation}
			\begin{aligned}\nonumber
				\Vert A_k^{\gamma(q),q} \Vert 
				&=\sum_{w\in E^k} p_w^q \alpha_1(w)^{\tau_w(q)} \alpha_2(w)^{\gamma(q)-\tau_w(q)}\\ &=\sum_{w\in E^k:a_w\geq b_w} p_w^q a_w^{\tau_A(q)} b_w^{\gamma(q)-\tau_A(q)} + \sum_{w\in E^k:a_w <b_w} p_w^q b_w^{\tau_B(q)} a_w^{\gamma(q)-\tau_B(q)}  \\
				&\leq \sum_{w\in E^k} p_w^q a_w^{\tau_A(q)} b_w^{\gamma(q)-\tau_A(q)} + \sum_{w\in E^k} p_w^q b_w^{\tau_B(q)} a_w^{\gamma(q)-\tau_B(q)} \\
				&\asymp \Vert \left(F^{(q)}_{\tau_A(q), \gamma(q)-\tau_A(q)} \right )^k\Vert +\Vert \left ( F^{(q)}_{\gamma(q)-\tau_B(q),\tau_B(q)} \right )^k  \Vert . 
			\end{aligned}
		\end{equation}
		By the definition of $\gamma(q)$, we have 
		$$
		\begin{aligned}
			1 &= \lim_{k\to \infty} \Vert A_k^{\gamma(q),q} \Vert ^{1/k} \\
			&\leq \lim_{k\to \infty} \max \left \{ \Vert \left (F^{(q)}_{\tau_A(q), \gamma(q)-\tau_A(q)} \right )^k\Vert^{1/k} ,\Vert \left ( F^{(q)}_{\gamma(q)-\tau_B(q),\tau_B(q)} \right )^k \Vert^{1/k}\right \}\\
			&= \max \left  \{\rho\left (F^{(q)}_{\tau_A(q),\gamma(q)-\tau_A(q)}\right ), \rho\left (F^{(q)}_{\gamma(q)-\tau_B(q),\tau_B(q)}\right ) \right \}
		\end{aligned}
		$$
		which gives that $\gamma(q) \leq \max \{\gamma_A(q),\gamma_B(q)\}$ by Lemma \ref{lemma41}. 
		
		For $x \in [\gamma_B(q)-\tau_B(q), \tau_A(q)]$, by the proof of Lemma \ref{l1}, we have $y(x)\in [ \gamma_A(q)-\tau_A(q),\tau_B(q) ]$. Therefore
		\begin{equation}\nonumber
			\begin{aligned}
				\Vert \left (F^{(q)}_{x,y(x)}\right )^k \Vert &\asymp \sum_{w\in E^k} p_w^q a_w^x b_w^{y(x)}\\
				&= \sum_{w\in E^k:a_w\geq b_w} p_w^q \alpha_1(w)^x \alpha_2(w)^{y(x)} +\sum_{w\in E^k:a_w <b_w} p_w^q \alpha_2(w)^x \alpha_1(w)^{y(x)}
				\\
				&\leq \sum_{w\in E^k:a_w\geq b_w} p_w^q \alpha_1(w)^{\tau_A(q)} \alpha_2(w)^{x+y(x)-\tau_A(q)} +\sum_{w\in E^k:a_w <b_w} p_w^q \alpha_1(w)^{\tau_B(q)} \alpha_2(w)^{x+y(x)-\tau_B(q)}
				\\
				&=\Vert A_k^{x+y(x),q} \Vert,
			\end{aligned}
		\end{equation}
		which yields $P(x+y(x),q)\geq 1$. It follows that $\gamma(q)\geq x+y(x)$ by Lemma \ref{propertyP}, which gives that \begin{equation}\label{e10}
			\gamma(q)\geq \max \{x+y(x): \gamma_B(q)-\tau_B(q)\leq x \leq \tau_A(q)\}.
		\end{equation}

		Combining the above observations with Lemma \ref{l1}, we obtain  \eqref{s20}.
		
		\vspace{0.2cm}

		\textit{Part II: When \eqref{e2} holds.}  \vspace{0.2cm}
		
		We firstly prove that $\gamma(q)\leq \min \{x+y(x): \tau_A(q)\leq x\leq \gamma_B(q)-\tau_B(q)\}.$
		For $x\in [\tau_A(q),\gamma_B(q)-\tau_B(q)]$, we have $y(x)\in [\tau_B(q), \gamma_A(q)-\tau_A(q)]$. Note that 
		\begin{equation}\nonumber
			\begin{aligned}
				\Vert A_k^{\gamma(q),q} \Vert &= \sum_{w\in E^k:a_w \geq b_w} p_w^q a_w^{\tau_A(q)} b_w^{\gamma(q)-\tau_A(q)} +\sum_{w\in E^k:a_w<b_w} p_w^q b_w^{\tau_B(q)} a_w^{\gamma(q)-\tau_B(q)}\\
				&\leq \sum_{w\in E^k}p_w^q a_w^{x} b_w^{\gamma(q)-x} +\sum_{w\in E^k} p_w^q a_w^{\gamma(q)-y(x)} b_w^{y(x)}\\
				&\asymp \max \left \{\Vert \left (F^{(q)}_{x,\gamma(q)-x}\right )^k \Vert ,\Vert \left (F^{(q)}_{\gamma(q)-y(x),y(x)}\right )^k \Vert \right \}.
			\end{aligned}
		\end{equation}
		It is easy to see that $\max \left \{\rho(F^{(q)}_{x,\gamma(q)-x} ), \rho (F^{(q)}_{\gamma(q)-y(x),y(x)} )\right \} \geq 1$, which gives that $\gamma(q)\leq x+y(x)$. Thus \begin{equation}\label{eq39}\gamma(q)\leq \min \{x+y(x):\tau_A(q)\leq x \leq \gamma_B(q)-\tau_B(q)\}.\end{equation}
		
		For estimating the lower bound of $\gamma(q)$, we consider the following three cases.
		\vspace{0.2cm}
		
		\textit{Case II-1: } $\sum_{e\in E} f^{(q)}_e(\tau_A(q)) \log (a_e/b_e) \geq  0$.
		\vspace{0.2cm}
		
		For   $k \geq 1$, $w\in E^k$, note that 
		$$
		\frac{\log a_w}{k}= \sum_{e\in E}\frac{\#(w,e)}{k} \log a_e
		$$
		and 
		$$
		\frac{\log b_w}{k}= \sum_{e\in E}\frac{\#(w,e)}{k} \log b_e.
		$$
		Let $\epsilon>0$
		and $\eta = - \epsilon \sum_{e\in E} \log a_eb_e>0.$ For  $w\in E^k \setminus \mathcal{B}_{\tau_A(q),k}(\epsilon)$,
		we have 
		\begin{equation}\nonumber
			\frac{\log a_w}{k} -\frac{\log b_w}{k} \geq \sum_{e\in E} f^{(q)}_e(\tau_A(q))\log a_e +\epsilon \sum_{e\in E} \log a_e -\sum_{e\in E} f^{(q)}_e(\tau_A(q))\log b_e +\epsilon \sum_{e\in E} \log b_e \geq -\eta,
		\end{equation}
		which gives that $$\frac{a_w}{b_w} \geq   e^{-\eta k}.$$
		By Lemma \ref{leclt}, picking $x=\tau_A(q)$ and an  initial distribution $\lambda=\left (\frac{1}{\#E}\right )_{e\in E}$,  there exists a positive constant $C<1$ such that for large enough $k$,
		\begin{equation}\nonumber
			\begin{aligned}
				&\sum_{w=w_1 \cdots w_k\in \mathcal{B}_{\tau_A(q),k}(\epsilon)} \lambda_{w_1} \tilde{F}_{\tau_A(q)}(w_1,w_2 )\cdots \tilde{F}_{\tau_A(q)}(w_{k-1},w_k)\\
				=&\sum_{w=w_1 \cdots w_k\in \mathcal{B}_{\tau_A(q),k}(\epsilon)} 
				\frac{1}{\#E} p_{w_2 \cdots w_k}^q a_{w_2 \cdots w_k}^{\tau_A(q)} b_{w_2 \cdots w_k}^{\gamma_A(q)-\tau_A(q)} \frac{u_{w_k}(\tau_A(q))}{u_{w_1}(\tau_A(q))}<C,
			\end{aligned}
		\end{equation}
		and so
		\begin{equation}
			\nonumber
			\begin{aligned}
				\sum_{w=w_1 \cdots w_k\in E^k \setminus \mathcal{B}_{\tau_A(q),k}(\epsilon)} \frac{1}{\#E} p_{w_2\cdots w_k} ^{q} a_{w_2 \cdots w_k}^{\tau_A(q)} b_{w_2 \cdots w_k} ^{\gamma_A(q)-\tau_A(q)} \frac{u_{w_k}(\tau_A(q))}{u_{w_1}(\tau_A(q))} \geq 1-C>0.
			\end{aligned}
		\end{equation}
		This yields for large enough $k$,
		\begin{equation}
			\nonumber
			\begin{aligned}
				\Vert A_{k}^{\gamma_A(q),q} \Vert &= \sum_{w\in E^{k}} p_w^q \alpha_1(w)^{\tau_w(q)} \alpha_2(w)^{\gamma_A(q)-\tau_w(q)}\\
				&\geq \sum_{w\in E^k \setminus \mathcal{B}_{\tau_A(q),k}(\epsilon)} p_w^q a_w^{\tau_A(q)} b_w^{\gamma_A(q)-\tau_A(q)} \min \left \{1, \left( \frac{a_w}{b_w} \right) ^{\gamma_A(q)-t(q)} \right \}
				\\
				&\geq \sum_{w\in E^k \setminus \mathcal{B}_{\tau_A(q),k}(\epsilon)} p_w^q a_w^{\tau_A(q)} b_w^{\gamma_A(q)-\tau_A(q)}   e  ^{-\eta k (\gamma_A(q)-t(q))} 
				\\
				&\geq \left(\sum_{w\in E^k \setminus \mathcal{B}_{\tau_A(q),k}(\epsilon)} \frac{1}{\#E} p_{w_2\cdots w_k} ^{q} a_{w_2 \cdots w_k}^{\tau_A(q)} b_{w_2 \cdots w_k} ^{\gamma_A(q)-\tau_A(q)} \frac{u_{w_k}(\tau_A(q))}{u_{w_1}(\tau_A(q))} \right)  e  ^{-\eta k (\gamma_A(q)-t(q))}  C' \\
				&\geq (1-C) e  ^{-\eta k (\gamma_A(q)-t(q))} C' ,
			\end{aligned}
		\end{equation}
		where $$C'=(\#E) p_*^q \min \left \{ \alpha_*^{\gamma_A(q)}, \alpha_*^{\tau_A(q)} \alpha^{*(\gamma_A(q)-\tau_A(q))}, \alpha^{*\tau_A(q)} \alpha_*^{\gamma_A(q)-\tau_A(q)}, \alpha^{*\gamma_A(q)}\right \} \cdot \min_{e,e'\in E} \frac{u_e(\tau_A(q))}{u_{e'}(\tau_A(q))} $$ is a positive number.
		Thus $P(\gamma_A(q),q) \geq e^{-\eta (\gamma_A(q)-t(q))}$. Let $\epsilon\to 0$, then $\eta \to 0$, so we have $P(\gamma_A(q),q)\geq 1$ which implies  $\gamma(q)\geq \gamma_A(q)$ by Lemma \ref{propertyP}. Hence, combining this with  \eqref{eq39}, we obtain  $\gamma(q)=\gamma_A(q)=\min\{x+y(x): \tau_A(q)\leq x\leq \gamma_B(q)-\tau_B(q)\}.$
		\vspace{0.2cm}
		
		\textit{Case II-2: } $\sum_{e\in E} f^{(q)}_e(\gamma_B(q)-\tau_B(q)) \log (a_e/b_e) \leq 0$.
		\vspace{0.2cm}
		
		Using a similar argument as above, it can be obtained that $\gamma(q)=\gamma_B(q)=\min\{x+y(x): \tau_A(q)\leq x\leq \gamma_B(q)-\tau_B(q)\}$.
		\vspace{0.2cm}
		
		\textit{Case II-3: Otherwise.}
		\vspace{0.2cm}
		
		In this case, $\sum_{e\in E} f^{(q)}_e(\tau_A(q)) \log (a_e/b_e) < 0$ and $\sum_{e\in E} f^{(q)}_e(\gamma_B(q)-\tau_B(q)) \log (a_e/b_e) >  0$. By Lemma \ref{lemma43} there exists $x \in [\tau_A(q),\gamma_B(q)-\tau_B(q)]$  such that 
		\begin{equation}\label{s24}
			\sum_{e\in E} f^{(q)}_e(x) \log (a_e /b_e) =0.
		\end{equation}
		It follows that $\tau_B(q)\leq y(x) \leq \gamma_A(q)-\tau_A(q).$ Fix this $x$ and again for $\epsilon>0$, let $\eta = - \epsilon \sum_{e\in E} \log a_eb_e$.
		For large enough $k$,  for $w\in E^k\setminus \mathcal{B}_{x,k}(\epsilon)$, using  \eqref{s24}, we have 
		$$
		-\eta \leq  \frac{\log a_w}{k} -\frac{\log b_w}{k} \leq \eta,
		$$
		which gives that 
		\begin{equation}\nonumber
			e^{-\eta k } \leq \frac{a_w}{b_w} \leq e^{ \eta k}.
		\end{equation}
		By Lemma \ref{leclt}, for an initial distribution $\lambda=\left (\frac{1}{\#E}\right )_{e\in E}$, there exists a positive constant $C<1$ such that 
		\begin{equation}
			\nonumber
			\sum_{w=w_1\cdots w_k \in E^k \setminus \mathcal{B}_{x,k}(\epsilon)} \frac{1}{\#E} p_{w_2 \cdots w_k}^q a_{w_2 \cdots w_k}^x b_{w_2 \cdots w_k}^{y(x)} \frac{u_{w_k}(x)}{u_{w_1}(x)} \geq 1-C>0.
		\end{equation}
		Thus 
		\begin{equation}
			\nonumber
			\begin{aligned}
				\Vert A_k^{x+y(x),q} \Vert &\geq \sum_{w\in E^k\setminus \mathcal{B}_{x,k}(\epsilon):a_w\geq b_w} p_w^q a_w^{\tau_A(q)} b_w^{x+y(x)-\tau_A(q)}
				+\sum_{w\in E^k\setminus \mathcal{B}_{x,k}(\epsilon):a_w< b_w} p_w^q b_w^{\tau_B(q)} a_w^{x+y(x)-\tau_B(q)}
				\\
				&= \sum_{w\in E^k\setminus \mathcal{B}_{x,k}(\epsilon):a_w\geq b_w} p_w^q a_w^{x} b_w^{y(x)} \left(\frac{a_w}{b_w}\right)^{\tau_A(q)-x}
				+\sum_{w\in E^k\setminus \mathcal{B}_{x,k}(\epsilon):a_w< b_w} p_w^q b_w^{y(x)} a_w^{x} \left(\frac{b_w}{a_w} \right)^{\tau_B(q)-y(x)}\\
				&\geq \sum_{w\in E^k\setminus \mathcal{B}_{x,k}(\epsilon)} p_w^q a_w^{x} b_w^{y(x)}  \min \left \{e^{\eta k(\tau_A(q)-x)},e^{\eta k
					(\tau_B(q)-y(x))} \right \}\\
				&\geq \sum_{w\in E^k\setminus \mathcal{B}_{x,k}(\epsilon)} \frac{1}{\#E}p_{w_2\cdots w_k}^q a_{w_2 \cdots w_k}^{x} b_{w_2 \cdots w_k}^{y(x)}  \frac{u_{w_k}(x)}{u_{w_1}(x)}\min \left \{e^{\eta k(\tau_A(q)-x)},e^{\eta k
					(\tau_B(q)-y(x))} \right \} C'\\
				&\geq (1-C)	\min \left \{e^{\eta k(\tau_A(q)-x)},e^{\eta k
					(\tau_B(q)-y(x))} \right \} C',
			\end{aligned}
		\end{equation}
		where $$C'=(\#E) p_*^q \min \left \{ \alpha_*^{x+y(x)}, \alpha_*^{x} \alpha^{*y(x)}, \alpha^{*x} \alpha_*^{y(x)}, \alpha^{*(x+y(x)) }\right \} \cdot \min_{e,e'\in E} \frac{u_e(x)}{u_{e'}(x)} $$ is a positive number.
		Therefore 
		$$
		P(x+y(x),q) \geq \min \left \{e^{\eta (\tau_A(q)-x)},e^{\eta 
			(\tau_B(q)-y(x))} \right \}.
		$$
		Letting $\epsilon \to 0$, $P(x+y(x),q)\geq 1$ by $\eta \to0 $. This yields $\gamma(q)\geq x+y(x)$ by Lemma \ref{propertyP}. Combining this with \eqref{eq39}, we obtain
		$\gamma(q)=\min \{x+y(x):\tau_A(q)\leq x \leq  \gamma_B(q)-\tau_B(q)\}.$
		\vspace{0.2cm}
		
		Combining \eqref{eq39} and all these three cases, using Lemma \ref{l1}, we finally obtain \eqref{s21}. 
		
	\end{proof}

	\begin{proof}[Proof of Corollary \ref{cor2}]
		By Theorem \ref{thmain2},
		it suffices to prove that 
		$$
		\max \{ \gamma_A(0),\gamma_B(0)\} \leq t(0).
		$$
		
		Suppose that  $\max \{ \gamma_A(0),\gamma_B(0)\} > t(0).$
		Assume $\gamma_A(0) \geq \gamma_B(0)$ without loss of generality. Noticing that $\rho(F^{(0)}_{\tau_A(0),\gamma_A(0)-\tau_A(0)})=1$, by Lemma \ref{lemma41}, we have 
		$\rho(F^{(0)}_{\tau_A(0),\tau_B(0)})>1$. Also, since $\rho(F^{(0)}_{\gamma_B(0)-\tau_B(0),\tau_B(0)})=1$, again using Lemma \ref{lemma41}, we know that $\gamma_B(0)>t(0)$. So $\min \{\gamma_A(0),\gamma_B(0)\}>t(0).$ Now 
		using Theorem \ref{thmain2}-(b), we have 
		$$
		\gamma(0)=\min \{x+y(x):\tau_A(0)\leq x \leq \gamma_B(0)-\tau_B(0)\}.
		$$
		
		On the other hand,
		by the product formula, for each $v\in V$, we have 
		$$
		{\dim}_B X_v \leq {\dim}_B \big(\pi_x(X_v) \times \pi_y(X_v)\big) \leq {\dim}_B  \pi_x(X_v) + {\dim}_B \pi_y(X_v).
		$$
		Since $\gamma(0)=\dim_B X_v$, $\tau_A(0)=\dim_B \pi_x(X_v)$ and $\tau_B(0) =\dim_B \pi_y( X_v)$, we have 
		$$
		\gamma(0) \leq t(0).
		$$ 
		So there exists a $x \in (\tau_A(0),\gamma_B(0)-\tau_B(0))$ such that $y(x)\leq t(0)-x<\tau_B(0)$. However, $y(\gamma_B(0)-\tau_B(0))=\tau_B(0)$, a contradiction to the fact that  $y(x)$ is decreasing by Lemma \ref{lemma43}.

	\end{proof}
	
	Now we  consider the degenerated case, i.e.  $V$ is a singleton. We will prove Corollary \ref{cor1} as a consequence of  Theorems \ref{thmain1} and \ref{thmain2}. At this time, $(V,E,\Psi)$ degenerates to a box-like self-affine IFS. The directed edge set $E$ can be written as $\{1,\cdots, N\}$. 
	\begin{lemma}\label{l2}
		Let $\{\psi_i\}_{i=1}^N$ be a planar  box-like self-affine IFS, and $\mathcal{P}=(p_i)_{i=1}^N$ be a positive probability vector.
		The function $y(x)$ is  uniquely determined by
		\begin{equation}\nonumber
			\sum_{i=1}^Np_i^q a_i^x b_i^{y(x)}=1.
		\end{equation}
		In addition, 
		$$
		y'(x)=-\frac{\sum_{i=1}^N p_i^q a_i^x b_i^{y(x)}\log a_i}{\sum_{i=1}^N p_i^q a_i^x b_i^{y(x)}\log b_i} \quad \text{ and }\quad y''(x)\geq 0.
		$$
	\end{lemma}
	\begin{proof}
		Note that all rows of  the $N\times N$ matrix $F^{(q)}_{x,y(x)}$ are same, so $$\rho(F^{(q)}_{x,y(x)})=\sum_{i=1}^N p_i^q a_i^x b_i^{y(x)}=1.$$
		Furthermore,
		by theorem of implicit function, we have
		\begin{equation}\label{e12}
			\sum_{i=1}^N p_i^q a_i^x b_i^{y(x)}(\log a_i +y'(x) \log b_i)=0,
		\end{equation}
		which implies that 
		$$
		y'(x)=-\frac{\sum_{i=1}^N p_i^q a_i^x b_i^{y(x)}\log a_i}{\sum_{i=1}^N p_i^q a_i^x b_i^{y(x)}\log b_i}.
		$$
		Differentiating \eqref{e12} implicitly with respect to $x$ gives 
		$$
		\sum_{i=1}^N p_i^q a_i^x b_i^{y(x)} ( (\log a_i +y'(x) \log b_i)^2 + y''(x) \log b_i)=0,
		$$
		which implies that $y''(x)\geq 0$.
	\end{proof}

	\begin{proof}[Proof of Corollary \ref{cor1}]
		Note that all rows of the matrix $F^{(q)}_{x,y(x)}$ are same, so its right Perron vector $u(x) =(\frac{1}{N})_{i=1}^N$ and left Perron vector $v(x)=(N p_i^q a_i^x b_i^{y(x)})_{i=1}^N$. Therefore, $f^{(q)}(x)=(p_i^q a_i^x b_i^y)_{i=1}^N$. Using Theorems \ref{thmain1} and \ref{thmain2}, the results of the corollary hold except in (b3) we need to prove  \begin{equation}\nonumber
			\tau_\mu(q) < \min \{\gamma_A(q),\gamma_B(q)\}.
		\end{equation}
		
		Since in (b3) neither \eqref{c1} or \eqref{c2} holds, i.e. $$\sum_{i=1}^N p_i^q a_i^{\tau_A(q)} b_i^{\gamma_A(q)-\tau_A(q)} \log a_i/b_i <0, \quad \sum_{i=1}^N p_i^q a_i^{\gamma_B(q)-\tau_B(q)} b_i^{\tau_B(q)} \log a_i/b_i>0,
		$$
		by Lemma \ref{l2},
		we have $y'(\tau_A(q))<-1$, $y'(\gamma_B(q)-\tau_B(q)) >-1$ and $\tau_A(q)<\gamma_B(q)-\tau_B(q)$. Hence for function $g(x):=x+y(x)$, $x \in [\tau_A(q),\gamma_B(q)-\tau_B(q)]$, it holds  $g'(\tau_A(q))<0$ and $g'(\gamma_B(q)-\tau_B(q))>0$. This gives that 
		$$
		\tau_\mu(q)=\min \{x+y(x): \tau_A(q) \leq x \leq \gamma_B(q)-\tau_B(q)\}<\min \{\gamma_A(q) , \gamma_B(q) \}.$$ 
		
	\end{proof}
	
	\subsection{Another proof of Corollary \ref{cor1}}\label{subsec32}
	In this subsection, we provide another proof of Corollary \ref{cor1} by using  a result of Feng and Wang \cite[Theorem 1]{FW05}. The main ingredient is to prove the following lemma.
	\begin{lemma}\label{lemma45}
		Let $\{\psi_i\}_{i=1}^N$, $\mathcal{P}$ and $\mu$ be same as in Corollary \ref{cor1}.  
		
		(a). If \eqref{c1} and \eqref{c2} both hold, $\tau_\mu(q)=\max \{\gamma_A(q),\gamma_B(q)\}.$
		
		(b). If \eqref{c1} holds, \eqref{c2} does not hold, $\tau_\mu(q)=\gamma_A(q)$.
		
		(c). If \eqref{c1} does not hold, \eqref{c2} holds, $\tau_\mu(q)=\gamma_B(q)$.
		
		(d). Otherwise, there exists a unique pair $x,y \in \mathbb{R}$ satisfying $\sum_{i=1}^N p_i^q a_i^x b_i^y=1$ and $\sum_{i=1}^N p_i^q a_i^x b_i^y\log a_i/b_i=0$ such that  $\tau_\mu(q)=x+y < \min \{\gamma_A(q),\gamma_B(q)\}$.
	\end{lemma}
	
	\begin{proof}
		
		For any vector $(d_1, \cdots, d_N)$, using $\Gamma(d_1, \cdots, d_N)$ to denote 
		$$
		\Gamma(d_1, \cdots, d_N):= \left \{ (t_1, \cdots, t_N): t_i\geq 0, \sum_{i=1}^N t_i =1, \sum_{i=1}^N t_i d_i \geq 0 \right \}  .
		$$
		Define two functions $f_A,f_B :\Omega=\{(t_i)_{i=1}^N: t_i\geq 0, \sum_{i=1}^{N} t_i=1\} \to \mathbb{R}$ by
		$$ 
		f_A\big((t_i)_{i=1}^N\big)= \frac{\sum_{i=1}^N t_i\big(\log t_i -\tau_A(q)(\log a_i-\log b_i) -q\log p_i\big)}{\sum_{i=1}^N t_i \log b_i},
		$$
		$$
		f_B\big((t_i)_{i=1}^N\big)=\frac{\sum_{i=1}^N t_i\big(\log t_i -\tau_B(q)(\log b_i-\log a_i) -q\log p_i\big)}{\sum_{i=1}^N t_i \log a_i}.
		$$
		Let 
		$$
		\theta_A= \max_{(t_i)_{i=1}^N \in \Gamma\big((\log a_i/b_i)_{i=1}^N\big)}  f_A\big((t_i)_{i=1}^N\big)
		$$
		and 
		$$
		\theta_B= \max_{(t_i)_{i=1}^N\in \Gamma\big((\log b_i /a_i)_{i=1}^N\big)} f_B\big((t_i)_{i=1}^N\big).
		$$
		By  Feng and Wang \cite[Theorem 1]{FW05},
		we have 
		\begin{equation}\label{c3}
			\tau_\mu(q)=\max \{ \theta_A,\theta_B\}.
		\end{equation}

		Now we analyze the extreme points of the $f_A$, $f_B$ on $\Omega$. By Lagrange multipliers method (following a similar calculation in \cite[Proposition 3.4]{LG92}), when $(t_i)_{i=1}^N=(p_i^q a_i^{\tau_A(q)} b_i^{\gamma_A(q)-\tau_A(q)})_{i=1}^N$, $f_A$ reaches a maximal value  $\gamma_A(q)$, and if $(p_i^q a_i^{\tau_A(q)} b_i^{\gamma_A(q)-\tau_A(q)})_{i=1}^N \notin \Gamma\big((\log a_i/b_i)_{i=1}^N\big)$,
		$$
		\theta_A= \max_{(t_i)_{i=1}^N\in \Gamma\big((\log a_i /b_i)_{i=1}^N\big) \cap \Gamma\big((\log b_i /a_i)_{i=1}^N\big)} f_A\big((t_i)_{i=1}^N\big).
		$$ 
		Also when $(t_i)_{i=1}^N= (p_i^q b_i^{\tau_B(q)} a_i^{\gamma_B(q)-\tau_B(q)})_{i=1}^N$, $f_B$ reaches a maximal value $\gamma_B(q)$, and if $(p_i^q b_i^{\tau_B(q)} a_i^{\gamma_B(q)-\tau_B(q)})_{i=1}^N\notin \Gamma\big((\log b_i/a_i)_{i=1}^N\big)$, 
		$$
		\theta_B= \max_{(t_i)_{i=1}^N\in \Gamma\big((\log a_i /b_i)_{i=1}^N\big) \cap \Gamma\big((\log b_i /a_i)_{i=1}^N\big)} f_B\big((t_i)_{i=1}^N\big).
		$$
		Note that \eqref{c1} is equivalent to $(p_i^q a_i^{\tau_A(q)} b_i^{\gamma_A(q)-\tau_A(q)})_{i=1}^N \in \Gamma\big((\log a_i/b_i)_{i=1}^N\big)$ and \eqref{c2} is equivalent to $(p_i^q b_i^{\tau_B(q)} a_i^{\gamma_B(q)-\tau_B(q)})_{i=1}^N\in \Gamma\big((\log b_i/a_i)_{i=1}^N\big)$.
		\vspace{0.2cm}
		
		If \eqref{c1} and \eqref{c2}  both hold, $\theta_A=\gamma_A(q)$ and  $\theta_B=\gamma_B(q)$, so (a) holds by \eqref{c3}.
		
		If \eqref{c1} holds and \eqref{c2} does not hold,
		noticing that when $$(t_i)_{i=1}^N \in \Gamma\big((\log a_i/b_i)_{i=1}^N\big) \cap \Gamma\big((\log b_i /a_i)_{i=1}^N\big),$$ $f_A=f_B$,  we have $\theta_A=\gamma_A(q) \geq \theta_B$, which  still by \eqref{c3} yields (b). Also,
		(c) follows by a same argument.
		
		If \eqref{c1} and \eqref{c2} both do not hold,
		by \eqref{c3}, we have 
		$$
		\tau_\mu(q)=\theta_A=\theta_B=\max_{(t_i)_{i=1}^N\in \Gamma\big((\log a_i /b_i)_{i=1}^N\big) \cap \Gamma\big((\log b_i /a_i)_{i=1}^N\big)} f_A\big((t_i)_{i=1}^N\big).
		$$
		Again using the Lagrange multipliers method, there exists a unique pair  $x,y \in \mathbb{R}$ satisfying $\sum_{i=1}^N p_i^q a_i^x b_i^y=1$ and $\sum_{i=1}^N p_i^q a_i^x b_i^y\log a_i/b_i=0$, such that  
		$$
		\tau_\mu(q)=x+y< \min \{ \gamma_A(q),\gamma_B(q)\}.
		$$ 
		This gives (d).
		
	\end{proof}
	
	Now we prove Corollary \ref{cor1} without using Theorems \ref{thmain1} and \ref{thmain2}.
	
	\begin{proof}[Another proof of Corollary \ref{cor1}]
		By Lemma \ref{l2},
		the function $y(x)$ is determined  by $\sum_{i=1}^{N} p_i^q a_i^x b_i^{y(x)}=1$,  \eqref{c1} is equivalent to $y'(\tau_A(q)) \geq -1$, \eqref{c2} is equivalent to $y'(\gamma_B(q)-\tau_B(q))\leq -1$ and  $y'(x)$ is increasing.
		\vspace{0.2cm}
		
		\textit{Case 1: \eqref{e3} holds, i.e. $\max \{ \gamma_A(q),\gamma_B(q)\} \leq t(q)$}. \vspace{0.2cm}
		
		In this case,
		$y'(\tau_A(q)) \geq y'(\gamma_B(q)-\tau_B(q))$. By Lemma \ref{l1}, it suffices to prove that 
		\begin{equation}\label{c5}
			\tau_\mu(q)=\max \{\gamma_A(q),\gamma_B(q)\}=\max \{x+y(x):\gamma_B(q)-\tau_B(q)\leq x\leq \tau_A(q)\}.
		\end{equation}
		
		If $y'(\gamma_B(q)-\tau_B(q))>-1$, we have \eqref{c1} holds, \eqref{c2} does not hold, and $x+y(x)$ is increasing in $[\gamma_B(q)-\tau_B(q),\tau_A(q)]$. So \eqref{c5} holds by Lemma \ref{lemma45}-(b).
		If $y'(\tau_A(q)) \geq -1$ and $y'(\gamma_B(q)-\tau_B(q))\leq -1$, we have \eqref{c1} and \eqref{c2} both hold, and $\max \{\gamma_A(q),\gamma_B(q)\}=\max \{x+y(x):\gamma_B-\tau_B(q)\leq x\leq \tau_A(q)\}$ which gives \eqref{c5} by Lemma \ref{lemma45}-(a).
		If $y'(\tau_A(q))<-1$, we have \eqref{c1} does not holds, \eqref{c2} holds, and $x+y(x)$ is decreasing, so   $\gamma_B(q)=\max \{x+y(x):\gamma_B-\tau_B(q)\leq x\leq \tau_A(q)\}$ which gives \eqref{c5} 
		by Lemma \ref{lemma45}-(c). 
		\vspace{0.2cm}
		
		\textit{Case 2: \eqref{e4} holds, i.e. $\min \{ \gamma_A(q),\gamma_B(q)\} \geq t(q)$}. \vspace{0.2cm}
		
		In this case, 
		$y'(\tau_A(q)) \leq y'(\gamma_B(q)-\tau_B(q))$. We aim to prove that 
		\begin{equation} \label{e7}
			\tau_\mu(q)=\min \{x+y(x):\tau_A(q)\leq x \leq \gamma_B(q)-\tau_B(q)\}.
		\end{equation}

		If $y'(\tau_A(q)) \geq -1$ and $y'(\gamma_B(q)-\tau_B(q))>-1$, $x+y(x)$ is increasing in $[\tau_A(q),\gamma_B(q)-\tau_B(q)]$, we have $\gamma_A(q)$ equals to the right of \eqref{e7}. Noticing that \eqref{c1} holds, \eqref{c2} does not hold, we have  \eqref{e7} holds by Lemma \ref{lemma45}-(b). If $y'(\tau_A(q)) \geq -1$ and $y'(\gamma_B(q)-\tau_B(q))\leq -1$, for all $x\in [\tau_A(q),\gamma_B(q)-\tau_B(q)]$,  $y'(x)=-1$, which yields that $\gamma_A(q)=\gamma_B(q)=\min \{x+y(x):\tau_A(q)\leq x \leq \gamma_B(q)-\tau_B(q)\}.$ Therefore \eqref{e7} holds by noticing that \eqref{c1} and \eqref{c2} both hold and using Lemma \ref{lemma45}-(a).
		If $y'(\tau_A(q)) <  -1$ and $y'(\gamma_B(q)-\tau_B(q)) > -1$, both \eqref{c1} and \eqref{c2} do not hold. By Lemma \ref{lemma45}-(d), there exists $x \in \mathbb{R}$ such that $y'(x)=-1$, so $x \in  (\tau_A(q),\gamma_B(q)-\tau_B(q))$. At this time,  $x+y(x)$ equals to the right side of \eqref{e7}, which then by Lemma \eqref{lemma45}-(d) implies  \eqref{e7} holds and 
		$$
		\tau_\mu(q)<\min \{\gamma_A(q),\gamma_B(q)\}.
		$$
		If $y'(\tau_A(q)) < -1$ and $y'(\gamma_B(q)-\tau_B(q))\leq  -1$, \eqref{c1} does not hold, \eqref{c2} holds, and $x+y(x)$ is decreasing, so  $\gamma_B(q)$ equals to the right of \eqref{e7}. By Lemma \ref{lemma45}-(c), we know that \eqref{e7} holds.
		
	\end{proof}
	
	\section{Closed forms in general setting}\label{subsec33}
	In this section, we turn to  the   general setting. We will prove Theorem \ref{thmain3}, Corollaries \ref{cor3} and \ref{cor4}. Still as above,  
	we always let  $(V,E,\Psi)$ be a strongly connected planar box-like self-affine  GIFS, but allowing some maps in $\Psi$ to be  anti-diagonal. Let $\mathcal{P}$, $\{\mu_v\}_{v\in V}$ be the associated positive vector and measures as before.
	We will present the closed form expression for $\gamma(q)$ ($L^q$-spectra of $\{\mu_v\}_{v\in V}$ by Theorem \ref{thmain1}). Also as above, when we define new variables, we may omit $q$.

 \subsection{Notations and lemmas}\label{subsec51}
	
	For $x,y\in \mathbb{R}$, $e,e'\in E$, we define a $2\times 2$ matrix $G_{x,y,e,e'}^{(q)}$ by
	\begin{equation}\nonumber
		G^{(q)}_{x,y,e,e'} = \left \{
		\begin{aligned}
			&\left(\begin{array}{cc}
				p_{e'}^q a_{e'}^{x+\tau_{x,e'}(q)} b_{e'}^{y-\tau_{x,e'}(q)} & 0 \\
				0   & p_{e'}^q b_{e'}^{x+\tau_{y,e'}(q)} a_{e'}^{y-\tau_{y,e'}(q)} \\	\end{array}\right) &\quad  &\text{ if }  t(e)=i(e') \text{ and }T_{e'} \text{ is diagonal, } \\
			&\left(\begin{array}{cc}
				0 & p_{e'}^q a_{e'}^{x+\tau_{y,e'}(q)} b_{e'}^{y-\tau_{y,e'}(q)} \\ p_{e'}^q
				b_{e'}^{x+\tau_{x,e'}(q)} a_{e'}^{y-\tau_{x,e'}(q)}   &0 \\	\end{array}\right)  &\quad&\text{ if } t(e)=i(e') \text{ and } T_{e'} \text{ is anti-diagonal,}
			\\
			&\left(\begin{array}{cc}
				0 & 0 \\
				0   &0 \\	\end{array}\right) & \quad & \textit{otherwise. }
		\end{aligned} \right .
	\end{equation}
        Denote $\{e(1),e(2)\}\times \{e'(1),e'(2)\}$ the collection of indices of matrix $G^{(q)}_{x,y,e,e'}$.
	We introduce a $\# E \times \# E$ block matrix $\mathcal{G}_{x,y}^{(q)}$ with entries  defined by
	\begin{equation}\nonumber
		\mathcal{G}_{x,y}^{(q)}[\{e(1),e(2)\},\{e'(1),e'(2)\}] =
		G_{x,y,e,e'}^{(q)}.
	\end{equation}

	Let  $\hat{E}=\{e(1),e(2):e\in E\}$, and 
	let $\iota:\hat{E} \to E$ be a \textit{projection map} so that  for $\hat{e}\in \hat{E}$, $\iota(\hat{e})=e\in E$ satisfying either $\hat{e}=e(1)$ or $\hat{e}=e(2)$.
	Let  $\kappa: \hat{E}\to \hat{E}$ be a one-to-one \textit{permutation map} so that $\kappa(e(1))=e(2)$, $\kappa(e(2))=e(1)$ for $e\in E$. 
	For $w=w_1 \cdots w_k \in E^*$ with $|w|\geq 2$, write $$G^{(q)}_{x,y,w}:=G^{(q)}_{x,y,w_1,w_2} \cdots G^{(q)}_{x,y,w_{k-1},w_k}.$$ 
	
	Define 
	$$
	\hat{E}^*=\{\hat{w}_1 \cdots \hat{w_k}:\hat{w}_i \in \hat{E}, \mathcal{G}^{(q)}_{x,y}(\hat{w}_{i-1},\hat{w}_i)>0, \forall 1<i\leq k, k\in \mathbb{N} \}.
	$$
	Note that the definition of $\hat{E}^*$ is independent of $q,x,y$. For $\hat{w}=\hat{w}_1 \cdots \hat{w}_k\in \hat{E}^*$, denote $|\hat{w}|=k$ the \textit{length} of $\hat{w}$. Denote $\hat{E}^k$ the collection of elements in $\hat{E}^*$ with length $k$.
	For $w=w_1 \cdots w_k \in E^*$, define $$\hat{E}(w)=\{ \hat{w}=\hat{w}_1 \cdots \hat{w}_k\in \hat{E}^k: \iota(\hat{w}_i)=w_i, \forall 1\leq i\leq k\}.$$
	
	\begin{lemma}\label{lemma46}
		For $w \in E^*$, we always have $\# \hat{E}(w)=2$.   Also for $k\in \mathbb{N}$,
    \begin{equation}\label{s26}
			\hat{E}^k=\bigcup_{w\in E^k} \hat{E}(w),
		\end{equation}
		where the  union is disjoint.
	\end{lemma}
	\begin{proof}
		For $w=w_1 \cdots w_k\in E^*$, it suffices to assume $k\geq 2$.
		Noticing that 
		\begin{equation}\label{s25}
			\hat{E}(w)=\{\hat{w}=\hat{w}_1 \cdots \hat{w}_{k} \in \prod_{i=1}^{k} \{ w_i(1),w_i(2)\}: \hat{w}\in \hat{E}^k\}
		\end{equation}
		and $G^{(q)}_{x,y,w}$ is a diagonal or anti-diagonal non-zero matrix, $\# \hat{E}(w)\geq 2$. Suppose $\# \hat{E}(w) \geq 3$, there must exist $\hat{w}\neq \hat{w}'\in \hat{E}(w)$ with  $\hat{w}_1= \hat{w}'_1$, such that  there exists  $1< i\leq k$ with 
		$$
		\hat{w}_{i-1} = \hat{w}_{i-1}' \quad \text{and} \quad
		\hat{w}_i\neq \hat{w}_i',
		$$
		which contradicted to that $G^q_{x,y,\iota(\hat{w}_{i-1}),\iota(\hat{w}_i)}$ is either  diagonal or anti-diagonal. So $\# \hat{E}(w)=2.$
		
		On the other hand, since 
		$$
		\hat{E}^k=\{\hat{w}=\hat{w}_1 \cdots \hat{w}_k \in  \prod_{i=1}^k \{w_i(1),w_i(2)\}:w \in E^k,\hat{w} \in \hat{E}^*\},$$
		using  \eqref{s25}, we obtain 
		\eqref{s26} and  the union in \eqref{s26} is disjoint.
	\end{proof}
	For $\hat{w}=\hat{w}_1\cdots \hat{w}_k \in \hat{E}^*$,
	let $\kappa(\hat{w})=\kappa(\hat{w}_1) \cdots \kappa(\hat{w}_k)$. It follows from  the proof of Lemma \ref{lemma46}, we know that $\kappa$ extends to a one-to-one map from $\hat{E}^*$ to $\hat{E}^*$.
	Still due to Lemma \ref{lemma46}, for $w=w_1 \cdots w_k\in E^*$, we could always  denote 
    $$
    \hat{E}(w)=\{{w}(1),{w}(2)\}$$ 
    with $${w}{(1)}=\hat{w}_1{} \cdots \hat{w}_k{}, \quad  \hat{w}_k{}={w}_k(1)$$ and $${w}{(2)}=\hat{w}_1'{} \cdots \hat{w}_k'{},\quad \hat{w}_k'{}={w}_k(2).$$ It is easy to see that $\kappa({w}{(1)})={w}{(2)}$ and $\kappa({w}{(2)})={w}{(1)}$.
	
	For $e\in E$, we write
	\begin{equation}\nonumber
		\begin{aligned} 
			a_{e(1)}&=a_e,\quad& b_{e(1)}&=b_e,\quad&
			a_{e(2)}&=b_e,\quad&
			b_{e(2)}&=a_e,\quad& & \text{ if } T_{e} \text{ is diagonal},\\
			a_{e(1)}&=b_e,\quad& b_{e(1)}&=a_e,\quad&
			a_{e(2)}&=a_e,\quad&
			b_{e(2)}&=b_e,\quad& & \text{ if } T_{e} \text{ is anti-diagonal},
		\end{aligned}
	\end{equation}
	also write
	$$
	\begin{aligned}
		&p_{e(1)}=p_{e(2)}=p_{e}, \\  &\tau_{e(1)}(q):=\tau_{x,e}(q),\\ &\tau_{e(2)}(q) :=\tau_{y,e}(q).
	\end{aligned}
	$$
	
	\begin{lemma} \label{lemma49}
		Matrices $G^{(q)}_{x,y,e,e'}$ and $\mathcal{G}^{(q)}_{x,y}$ could be written as
		\begin{equation}\nonumber
			G^{(q)}_{x,y,e,e'} = \left \{
			\begin{aligned}
				&\left(\begin{array}{cc}
					p_{e'(1)}^q a_{e'(1)}^{x+\tau_{e'(1)}(q)} b_{e'(1)}^{y-\tau_{e'(1)}(q)} & 0 \\
					0   & p_{e'(2)}^q a_{e'(2)}^{x+\tau_{e'(2)}(q)} b_{e'(2)}^{y-\tau_{e'(2)}(q)} \\	\end{array}\right) &  &\text{ if } t(e)=i(e') \text{, } T_{e'} \text{ is diagonal, } \\
				&\left(\begin{array}{cc}
					0 & p_{e'(2)}^q a_{e'(2)}^{x+\tau_{e'(2)}(q)} b_{e'(2)}^{y-\tau_{e'(2)}(q)} \\ p_{e'(1)}^q
					a_{e'(1)}^{x+\tau_{e'(1)}(q)} b_{e'(1)}^{y-\tau_{e'(1)}(q)}   &0 \\	\end{array}\right)  & &\text{ if }  t(e)=i(e') \text{, } T_{e'} \text{ is anti-diagonal,}
				\\
				&\left(\begin{array}{cc}
					0 & 0 \\
					0   &0 \\	\end{array}\right) &  & \textit{otherwise, }
			\end{aligned} \right .
		\end{equation}
		and 
		$$
		\mathcal{G}^{(q)}_{x,y}(\hat{e},\hat{e}')=
		\left \{
		\begin{aligned}
			&p_{\hat{e}'}^q a_{\hat{e}'}^{x+\tau_{\hat{e}'}(q)} b_{\hat{e}'}^{y-\tau_{\hat{e}'}(q)} &\quad&  \text{ if } \hat{e}\hat{e}'\in \hat{E}^*, 
		  \\
			&0 & \quad & \textit{otherwise. }
		\end{aligned}
		\right.
		$$
	\end{lemma}
	\begin{proof}
		This can be directly seen by the  definitions of $p_{\hat{e}},a_{\hat{e}}, b_{\hat{e}}$, $\tau_{\hat{e}}(q)$ and matrices $G^{(q)}_{x,y,e,e'}$, $\mathcal{G}^{(q)}_{x,y}$.
	\end{proof}
	\begin{lemma}\label{lemma6}
		For $\hat{e},\hat{e}'\in \hat E$ with $\hat{e}\hat{e}'\in \hat{E}^*$, we have 
		\begin{equation}\nonumber
			\begin{aligned}
				&\tau_{\hat{e}}(q)=\tau_{\hat{e}'}(q).
			\end{aligned}
		\end{equation}
	\end{lemma}
	\begin{proof}
   Note that when $T_{\iota(\hat{e}')}$ is diagonal, 
        $$
        \text{ either } \quad \tau_{\hat{e}}(q)=\tau_{x,\iota(\hat{e})}(q), \tau_{\hat{e}'}(q)=\tau_{x,\iota(\hat{e}')}(q) \quad \text{ or } \quad \tau_{\hat{e}}(q)=\tau_{y,\iota(\hat{e})}(q), \tau_{\hat{e}'}(q)=\tau_{y,\iota(\hat{e}')}(q),
        $$
        and when $T_{\iota(\hat{e}')}$ is anti-diagonal,
        $$
        \text{ either } \quad \tau_{\hat{e}}(q)=\tau_{x,\iota(\hat{e})}(q), \tau_{\hat{e}'}(q)=\tau_{y,\iota(\hat{e}')}(q) \quad \text{ or } \quad \tau_{\hat{e}}(q)=\tau_{y,\iota(\hat{e})}(q), \tau_{\hat{e}'}(q)=\tau_{x,\iota(\hat{e}')}(q).
        $$
        The lemma follows from the proof of Lemma \ref{lemma7}.
	\end{proof}
	For $\hat{w}=\hat{w}_1 \cdots \hat{w}_k\in \hat{E}^*$, write $p_{\hat{w}}:= p_{\hat{w}_1} \cdots p_{\hat{w}_k}$, $a_{\hat{w}}:=a_{\hat{w}_1} \cdots a_{\hat{w}_k}$ and $b_{\hat{w}}:=b_{\hat{w}_1} \cdots b_{\hat{w}_k}$. For $w\in E^*$, recall  that $c_w$,$d_w$ are the width and height of the rectangle $\psi_w([0,1]^2)$. 
	
	\begin{lemma}\label{lemma3}
		For $w\in E^*$, we have 
		$$
		\begin{aligned}
			&c_w=a_{{w}{(1)}}=b_{{w}{(2)}}, \quad d_w=a_{{w}{(2)}}=b_{{w}{(1)}} \quad & & \text{ if } T_w \text{ is diagonal,}\\
			&c_w=a_{{w}{(2)}}=b_{{w}{(1)}}, \quad d_w=a_{{w}{(1)}}=b_{{w}{(2)}} \quad & & \text{ if } T_w \text{ is anti-diagonal.}
		\end{aligned}
		$$
	\end{lemma}
	\begin{proof}
		For $e\in E$ and $w=w_1 \cdots w_k \in E^*$ with $e w\in E^*$, note that  $G^{(q)}_{x,y,ew}= G^{(q)}_{x,y,e,w_1} \cdots G^{(q)}_{x,y,w_{k-1},w_k}$. By Lemmas \ref{lemma49}, \ref{lemma6}, and the definitions of ${w}{(1)}$ ,${w}{(2)}$, it is directly to check that  
		\begin{equation}\label{s31}
			G^{(q)}_{x,y,ew}=\left \{
			\begin{aligned}
				&\left(\begin{array}{cc}
					p_{{w}{(1)}}^q a_{{w}{(1)}}^{x+\tau_{w_k(1)}(q)} b_{{w}{(1)}}^{y-\tau_{w_k(1)}(q)} & 0 \\
					0   & p_{{w}{(2)}}^q a_{{w}{(2)}}^{x+\tau_{w_k(2)}(q)} b_{{w}{(2)}}^{y-\tau_{w_k(2)}(q)} \\	\end{array}\right) &\quad  &\text{ if }  T_{w}  \text{ is diagonal,} \\
				&\left(\begin{array}{cc}
					0 & p_{{w}{(2)}}^q a_{{w}{(2)}}^{x+\tau_{w_k(2)}(q)}b_{{w}{(2)}}^{y-\tau_{w_k(2)}(q)} \\ p_{{w}{(1)}}^q a_{{w}{(1)}}^{x+\tau_{w_k(1)}(q)}b_{{w}{(1)}}^{y-\tau_{w_k(1)}(q)}   &0 \\	\end{array}\right)  &\quad&\text{ if }  T_{w} \text{ is anti-diagonal.}
			\end{aligned} \right .
		\end{equation}  
		On the other hand,
		noticing that the absolution values of  nonzero element of the matrix $T_w$ in the first (resp. second) row is $c_w$ (resp. $d_w$). So 
		\begin{equation}\label{s32}
			G^{(q)}_{x,y,ew}=\left \{
			\begin{aligned}
				&\left(\begin{array}{cc}
					p_w^q c_{w}^{x+\tau_{x,w_k}(q)} d_{w}^{y-\tau_{x,w_k}(q)} & 0 \\
					0   & p_w^q d_{w}^{x+\tau_{y,w_k}(q)} c_{w}^{y-\tau_{y,w_k}(q)} \\	\end{array}\right) &\quad  &\text{ if }  T_{w} =\left(\begin{array}{cc}
					\pm c_w & 0 \\ 
					0 & \pm d_w \\	\end{array}\right), \\
				&\left(\begin{array}{cc}
					0 & p_w^q c_w^{x+\tau_{y,w_k}(q)} d_{w}^{y-\tau_{y,w_k}(q)} \\ p_w^q
					d_{w}^{x+\tau_{x,w_k}(q)} c_{w}^{y-\tau_{x,w_k}(q)}   &0 \\	\end{array}\right)  &\quad&\text{ if }  T_{w} =\left(\begin{array}{cc}
					0 & \pm c_w \\ 
					\pm d_w   &0 \\	\end{array}\right).
			\end{aligned} \right .
		\end{equation}
		The lemma follows immediately by 
		comparing \eqref{s31} with \eqref{s32}.
	\end{proof}

	\begin{lemma}\label{lemma47}
		The function $\rho(\mathcal{G}^{(q)}_{x,y})$ is continuous in $x,y\in \mathbb{R}$. For fixed $y\in \mathbb{R}$, $\rho(\mathcal{G}^{(q)}_{x,y})$ is strictly decreasing in $x \in \mathbb{R}$, and there exists a unique $x \in \mathbb{R}$ such that $\rho(\mathcal{G}^{(q)}_{x,y})=1$. This is also true for $\rho(\mathcal{G}^{(q)}_{x,y})$ as a function of $y\in \mathbb{R}$ for fixed $x \in \mathbb{R}$.
	\end{lemma}
	\begin{proof}
		Using a same argument in the proof of  Lemma \ref{lemma41}, this lemma follows.
	\end{proof}
	By Lemma \ref{lemma47}, we can define a function $\hat{y}(x):\mathbb{R}\to \mathbb{R}$ by
	$$
	\rho(\mathcal{G}^{(q)}_{x,\hat{y}(x)})=1.
	$$

	By Proposition \ref{thmain4}, for each $e\in E$, we always have either 
	$$
	\tau_{x,e}(q)=\tau_A(q),\quad \tau_{y,e}(q)=\tau_B(q),
	$$
	or 
	$$
	\tau_{x,e}(q)=\tau_B(q),\quad \tau_{y,e}(q)=\tau_A(q),
	$$
	where $A,B$ are same as in Proposition \ref{thmain4}.
	This means that \begin{equation}\label{s76}
 \begin{aligned}
	\tau_{x,e}(q)+\tau_{y,e}(q)&=t(q), \quad & &\text{ for all } e\in E,\\
 \tau_{\hat{e}}(q)+\tau_{\kappa(\hat{e})}(q)&=t(q), \quad& &  \text{ for all }\hat{e}\in \hat{E},
 \end{aligned}
 \end{equation}
 where $t(q)=\tau_A(q)+\tau_B(q).$
	
	Recall that  a square matrix $Z$ is a \textit{permutation matrix} if 
	every row and column of $Z$ contains precisely one $1$ with other entries $0$.
	\begin{lemma}\label{lemma48}
		There exists a permutation matrix $Z$ such that 
		for $x,y\in \mathbb{R}$, we always have 
		$$
		\mathcal{G}^{(q)}_{x,y}=Z\mathcal{G}^{(q)}_{y-t(q),x+t(q)}Z.
		$$
	\end{lemma}
	\begin{proof}
		Let ${Z}=\text{diag}\left \{\left(\begin{array}{cc}
			0 & 1 \\
			1 & 0 \\	\end{array}\right), \cdots, \left(\begin{array}{cc}
			0 & 1 \\
			1 & 0 \\	\end{array}\right)\right \}$ be a $\#E\times \#E$ block diagonal matrix.
		So ${Z}$ is a $(2\#E)\times (2\#E)$ permutation  matrix.
		Note that for each $e \in E$,  we always have $$\tau_{x,e}(q)+\tau_{y,e}(q)=t(q).$$
		For $e,e'\in E$, note that by definition,
		$$
		G^{(q)}_{x,y,e,e'}
		=\left(\begin{array}{cc}
			0 & 1 \\
			1 & 0 \\	\end{array}\right)
		G^{(q)}_{y-t(q),x+t(q),e,e'}
		\left(\begin{array}{cc}
			0 & 1 \\
			1 & 0 \\	\end{array}\right).
		$$
		By the definition of $\mathcal{G}^{(q)}_{x,y}$ and $Z$, the lemma follows.
	\end{proof}
	\begin{remark}\label{s77}
		By Lemma \ref{lemma48}, we know that $\mathcal{G}^{(q)}_{x,y}$ and $\mathcal{G}^{(q)}_{y-t(q),x+t(q)}$ are only  different by permutations, so  $
		\rho(\mathcal{G}^{(q)}_{x,y}) =\rho(\mathcal{G}^{(q)}_{y-t(q),x+t(q)})$. 
		For $q\geq 0$, define   $$\hat{\gamma}(q)=\hat{y}(0).$$ Then $\hat{\gamma}(q)$  satisfies $$
		\rho(\mathcal{G}^{(q)}_{0,\hat{\gamma}(q)})=\rho(\mathcal{G}^{(q)}_{\hat{\gamma}(q)-t(q),t(q)})=1,
		$$ so that   $\hat{\gamma}(q)$ is well-defined in \eqref{s28}.
	\end{remark}
	\begin{lemma}\label{lemma1}
		Either $\mathcal{G}^{(q)}_{x,y}$ is irreducible, or there exist $\hat{E}',\hat{E}''\subseteq \hat{E}$ with $\#\hat{E}'=\# \hat{E}''=\#E$, $\hat{E}'\cup\hat{E}''=\hat{E}$, $\hat{E}'\cap \hat{E}''=\emptyset$ and $\kappa(\hat{E}')=\hat{E}''$ such that both  $\mathcal{G}^{(q)}_{x,y}[\hat{E}',\hat{E}']$, $\mathcal{G}^{(q)}_{x,y}[\hat{E}'',\hat{E}'']$ are irreducible, and 
		\begin{equation}\label{s40}
			\mathcal{G}^{(q)}_{x,y}[\hat{E}',\hat{E}'](\hat{e},\hat{e}')=\mathcal{G}^{(q)}_{y-t(q),x+t(q)}[\hat{E}'',\hat{E}''](\kappa(\hat{e}),\kappa(\hat{e}')),\quad  \text{ for all } \hat{e},\hat{e}' \in \hat{E}'.
		\end{equation}
		In the later case, there exists a $\# \hat{E} \times \# \hat{E}$  permutation matrix $Z'$ such that  
		\begin{equation}\label{s39}        \mathcal{G}^{(q)}_{x,y}=Z' \left(\begin{array}{cc}
				\mathcal{G}^{(q)}_{x,y}[\hat{E}',\hat{E}'] & 0 \\
				0 & \mathcal{G}^{(q)}_{x,y}[\hat{E}'',\hat{E}''] \\	\end{array}\right) Z'.
		\end{equation}
		
	\end{lemma}
	\begin{proof}
		The proof is basing on a same idea as that of Proposition \ref{thmain4}. 
		
		It suffices to assume that  $\mathcal{G}^{(q)}_{x,y}$ is not irreducible. So we can pick $\hat{e}', \hat{e}''\in \hat{E}$ such that 
		\begin{equation}\label{s35}
			\left(\mathcal{G}^{(q)}_{x,y}\right)^k(\hat{e}',\hat{e}'')=0, \quad  \text{ for all } k\in \mathbb{N}.
		\end{equation}
		Let 
		\begin{equation}\nonumber
			\begin{aligned}
				\hat{E}'=\left \{\hat{e}\in \hat{E}:\left(\mathcal{G}^{(q)}_{x,y}\right)^k(\hat{e}',\hat{e})>0 \text{ for some } k\in \mathbb{N}\right \},
				\\
				\hat{E}''=\left \{\hat{e}\in \hat{E}:\left(\mathcal{G}^{(q)}_{x,y}\right)^k(\hat{e},\hat{e}'')>0 \text{ for some } k\in \mathbb{N}\right \}.
			\end{aligned}
		\end{equation}
		By \eqref{s35}, $\hat{E}' \cap \hat{E}''=\emptyset.$ Since $(V,E)$ is strongly connected, for each $e\in E$, there exist $k_1,k_2 \in \mathbb{N}$ such that 
		$$
		\left(\mathcal{G}^{(q)}_{x,y}\right)^{k_1}\left[\{\hat{e}'\},\{e(1),e(2)\}\right] \text{ and } \left(\mathcal{G}^{(q)}_{x,y}\right)^{k_2}\left[\{e(1),e(2)\}, \{\hat{e}''\}\right] \text{ are non-zero matrices. }
		$$
		Noticing that $$
		0=\left(\mathcal{G}^{(q)}_{x,y}\right)^{k_1+k_2}(\hat{e}',\hat{e}'') \geq \left(\mathcal{G}^{(q)}_{x,y}\right)^{k_1}\left[\{\hat{e}'\},\{e(1),e(2)\}\right] \cdot \left(\mathcal{G}^{(q)}_{x,y}\right)^{k_2}\left[\{e(1),e(2)\}, \{\hat{e}''\}\right],
		$$
		we have
		\begin{equation} \label{s41}
			\text{ either } e(1)\in \hat{E}', e(2)\in \hat{E}'' \text{ or }  e(1)\in \hat{E}'', e(2) \in \hat{E}'.
		\end{equation} 
		Thus  \begin{equation}\nonumber
			\hat{E}'\cup \hat{E}''=\hat{E}, \# \hat{E}'=\#\hat{E}''=\# E \text{ and } \kappa(\hat{E}')=\hat{E}'' .
		\end{equation}
		
		For  $e\in E$, by \eqref{s35}, \eqref{s41} and $\hat{E}'\cap \hat{E}''=\emptyset$, we also have
		\begin{equation}\label{s62}
			\begin{aligned}\left(\mathcal{G}^{(q)}_{x,y}\right)^k(e(1),e(2))&=0, \quad  \text{ for all } k\in \mathbb{N},\quad & &\text{ if } e(1)\in \hat{E}',e(2)\in \hat{E}'',
				\\
				\left(\mathcal{G}^{(q)}_{x,y}\right)^k(e(2),e(1))&=0, \quad \text{ for all } k\in \mathbb{N}, \quad & &\text{ if } e(1)\in \hat{E}'',e(2)\in \hat{E}'.
			\end{aligned}
		\end{equation}
		When $e(1)\in \hat{E}',e(2)\in \hat{E}''$,
		let 
		\begin{equation}
			\nonumber
			\begin{aligned}
				\hat{E}_{e}'&=\left \{\hat{e}\in \hat{E}:\left(\mathcal{G}^{(q)}_{x,y}\right)^k(e(1),\hat{e})>0 \text{ for some } k\in \mathbb{N}\right \},
				\\
				\hat{E}_{e}''&=\left \{\hat{e}\in \hat{E}:\left(\mathcal{G}^{(q)}_{x,y}\right)^k(\hat{e},e(2))>0 \text{ for some } k\in \mathbb{N}\right \}.
			\end{aligned}
		\end{equation}
		When $e(1)\in \hat{E}'',e(2)\in \hat{E}'$,
		let 
		\begin{equation}
			\nonumber
			\begin{aligned}
				\hat{E}_{e}'&=\left \{\hat{e}\in \hat{E}:\left(\mathcal{G}^{(q)}_{x,y}\right)^k(e(2),\hat{e})>0 \text{ for some } k\in \mathbb{N}\right \},
				\\
				\hat{E}_{e}''&=\left \{\hat{e}\in \hat{E}:\left(\mathcal{G}^{(q)}_{x,y}\right)^k(\hat{e},e(1))>0 \text{ for some } k\in \mathbb{N}\right \}.
			\end{aligned}
		\end{equation}
		By \eqref{s41} and \eqref{s62}, we  know that
		$\hat{E}_e'\cap\hat{E}''=\emptyset$, which implies that 
		\begin{equation}\nonumber
			\hat{E}_{e}'=\hat{E}',  \quad \hat{E}_{e}''=\hat{E}''.
		\end{equation}
		Also, it is direct to check that $\hat{e}'\in \hat{E}_{e}'=\hat{E}' $ by a contradiction argument.

		Combining \eqref{s41} and the definition of $\hat{E}',\hat{E}'', \hat{E}_e',\hat{E}_e''$, we have $\mathcal{G}^{(q)}_{x,y}[\hat{E}',\hat{E}']$ is  irreducible and 
		\begin{equation}\label{s79}
		\mathcal{G}^{(q)}_{x,y}(\hat{e}''',\hat{e}'''')=0, \quad  \text{ for all } \hat{e}'''\in \hat{E}',\hat{e}''''\in \hat{E}''.
		\end{equation}
		On the other hand,
		for each $\hat{e}''',\hat{e}''''\in \hat{E}'$, using the proof of  Lemma \ref{lemma48}, 
		\begin{equation}\label{s80}
		\mathcal{G}^{(q)}_{x,y}(\hat{e}''',\hat{e}'''')=\mathcal{G}^{(q)}_{y-t(q),x+t(q)}(\kappa(\hat{e}'''),\kappa(\hat{e}'''')).
		\end{equation}
		So by \eqref{s41}, \eqref{s40} holds and as a consequence, $\mathcal{G}^{(q)}_{x,y}[\hat{E}'',\hat{E}'']$ is irreducible. Finally, it follows from \eqref{s79}, \eqref{s80}, $\kappa(\hat{E}')=\hat{E}''$ and $\kappa(\hat{E}'')=\hat{E}'$, we  obtain \eqref{s39}.
	\end{proof}
	Now we will introduce a vector-valued function $g^{(q)}(x):\mathbb{R} \to \mathbb{R}^{\# \hat{E}}$ analogous to $f^{(q)}(x)$ in \eqref{s71} in the non-rotational setting.
	
	When $\mathcal{G}^{(q)}_{x,\hat{y}(x)}$ is irreducible, let $\hat{u}(x)=(\hat{u}_{\hat{e}}(x))_{\hat{e}\in \hat{E}}$ (resp. $\hat{v}(x)=(\hat{v}_{\hat{e}}(x))_{\hat{e}\in \hat{E}}$) be the right (resp. left) Perron vector of $\mathcal{G}^{(q)}_{x,\hat{y}(x)}$. Define a $\# \hat{E} \times \# \hat{E}$ irreducible stochastic matrix $\tilde{\mathcal{G}}_x$  with entries
	\begin{equation}\label{s72}
	\tilde{\mathcal{G}}_x(\hat{e},\hat{e}')=\mathcal{G}^{(q)}_{x,\hat{y}(x)}(\hat{e},\hat{e}') \frac{\hat{u}_{\hat{e}'}(x)}{\hat{u}_{\hat{e}}(x)},
	\end{equation}
	and a positive probability row vector $$g^{(q)}(x):=(g^{(q)}_{\hat{e}}(x))_{\hat{e}\in \hat{E}}=(\hat{v}_{\hat{e}}(x) \hat{u}_{\hat{e}}(x))_{\hat{e}\in \hat{E}}.$$ So $g^{(q)}(x)$ is an invariant distribution associated with $\tilde{\mathcal{G}}_x$, i.e.
	\begin{equation}\label{s73}
	g^{(q)}(x) \tilde{\mathcal{G}}_x =g^{(q)}(x).
	\end{equation}
	
	When $\mathcal{G}^{(q)}_{x,\hat{y}(x)}$ is not irreducible, using Lemma \ref{lemma1}, we see that $\mathcal{G}^{(q)}_{x,\hat{y}(x)}[\hat{E}',\hat{E}']$ and $\mathcal{G}^{(q)}_{x,\hat{y}(x)}[\hat{E}'',\hat{E}'']$ are two irreducible matrices. Use similar argument as above to $\mathcal{G}^{(q)}_{x,\hat{y}(x)}[\hat{E}',\hat{E}']$ and $\mathcal{G}^{(q)}_{x,\hat{y}(x)}[\hat{E}'',\hat{E}'']$ respectively, let $\hat{u}(x,1)$ (resp. $\hat{u}(x,2)$)  be the right Perron vector of $\mathcal{G}^{(q)}_{x,\hat{y}(x)} [\hat{E}',\hat{E}']$ (resp. $\mathcal{G}^{(q)}_{x,\hat{y}(x)}[\hat{E}'',\hat{E}'']$), and let $\hat{v}(x,1)$ (resp. $\hat{v}(x,2)$) be the left Perron vector of $\mathcal{G}^{(q)}_{x,\hat{y}(x)} [\hat{E}',\hat{E}']$ (resp. $\mathcal{G}^{(q)}_{x,\hat{y}(x)}[\hat{E}'',\hat{E}'']$). Then we may define two stochastic matrices $\tilde{\mathcal{G}}_{x,1}$, $\tilde{\mathcal{G}}_{x,2}$ and two positive probability row vector $g^{(q)}(x,1)$, $g^{(q)}(x,2)$ as above  satisfying 
	$$
	g^{(q)}(x,1)\tilde{\mathcal{G}}_{x,1} =g^{(q)}(x,1) \text{ and }
	g^{(q)}(x,2) \tilde{\mathcal{G}}_{x,2} =g^{(q)}(x,2).
	$$
	Let $\hat{u}(x)=Z'\left(\begin{aligned}
		\hat{u}(x,1)\\
		\hat{u}(x,2)
	\end{aligned}\right)$, $\hat{v}(x)=\left(\hat{v}(x,1),\hat{v}(x,2)\right) Z'$,
	\begin{equation}\label{s74}
	\tilde{\mathcal{G}}_{x}=Z' \left(\begin{array}{cc}
		\tilde{\mathcal{G}}_{x,1} & 0 \\
		0 & \tilde{\mathcal{G}}_{x,2} \\	\end{array}\right) Z',
	\end{equation}
	where $Z'$ is the same in Lemma \ref{lemma1}. Define  $$g^{(q)}(x)=(g^{(q)}(x,1),g^{(q)}(x,2))Z'.$$ Then $\tilde{\mathcal{G}}_x$ is a stochastic matrix satisfying  \eqref{s72} and $g^{(q)}(x)$ is an invariant vector satisfying \eqref{s73}.

	\begin{lemma}\label{lemma2}
		The function $\hat{y}(x)$ is continuous, decreasing, and $g^{(q)}(x)$ is continuous in $x \in \mathbb{R}$.
	\end{lemma}
	\begin{proof}
		The lemma  follows from a same proof of Lemma \ref{lemma43} by using Lemma  \ref{lemma47}.
	\end{proof}
	
	\begin{lemma}\label{lemma4}
		There exists $x\in [\min \{ 0,\hat{\gamma}(q)-t(q)\},\max \{ 0,\hat{\gamma}(q)-t(q)\}]$ such that 
		$$
		\sum_{\hat{e}\in \hat{E}} g^{(q)}_{\hat{e}}(x) \log a_{\hat{e}}/b_{\hat{e}}=0.
		$$
	\end{lemma}
	\begin{proof}
		We only prove the case $\hat{\gamma}(q)\geq t(q)$. By Lemma \ref{lemma2}, it suffices to prove that 
		\begin{equation}\label{s44}
			\left(  \sum_{\hat{e}\in \hat{E}} g^{(q)}_{\hat{e}}(0) \log a_{\hat{e}}/b_{\hat{e}} \right) \left(  \sum_{\hat{e}\in \hat{E}} g^{(q)}_{\hat{e}}(\hat{\gamma}(q)-t(q)) \log a_{\hat{e}}/b_{\hat{e}} \right) \leq 0.
		\end{equation}
  Due to Lemma \ref{lemma1}, we consider two possible cases.
		\vspace{0.2cm}
		
		\textit{Case 1: $\mathcal{G}^{(q)}_{x,\hat{y}(x)}$ is irreducible.}
		\vspace{0.2cm}
		
		Using  $\mathcal{G}^{(q)}_{0,\hat{\gamma}(q)}=Z \mathcal{G}^{(q)}_{\hat{\gamma}(q)-t(q),t(q)}Z$ from Lemma \ref{lemma48}, 
		 and noticing that 
		 $\mathcal{G}^{(q)}_{0,\hat{\gamma}(q)} \hat{u}(0)=\hat{u}(0)$ and  $\mathcal{G}^{(q)}_{\hat{\gamma}(q)-t(q),t(q)} \hat{u}(\hat{\gamma}(q)-t(q)) =\hat{u}(\hat{\gamma}(q)-t(q))$, we have $Z \hat{u}(0)=\hat{u}(\hat{\gamma}(q)-t(q))$,
		 where  ${Z}=\text{diag}\left \{\left(\begin{array}{cc}
		 	0 & 1 \\
		 	1 & 0 \\	\end{array}\right), \cdots, \left(\begin{array}{cc}
		 	0 & 1 \\
		 	1 & 0 \\	\end{array}\right)\right \}$ is a $\#E\times \#E$ block diagonal matrix. Similarly, $\hat{v}(0)Z=\hat{v}(\hat{\gamma}(q)-t(q)).$ So $g^{(q)}(0)Z =g^{(q)}(\hat{\gamma}(q)-t(q))$.
		Note that for each $e\in E$, $\log a_{e(1)}/b_{e(1)} =-\log a_{e(2)}/b_{e(2)}.$ Thus 
		$$
		\begin{aligned}
			\sum_{\hat{e}\in \hat{E}} g_{\hat{e}}^{(q)}(0) \log a_{\hat{e}}/b_{\hat{e}}&=\sum_{e\in E}\left( g^{(q)}_{e(1)}(0) \log a_{e(1)}/b_{e(1)} + g^{(q)}_{e(2)}(0) \log a_{e(2)}/b_{e(2)}\right) 
			\\
			&=-\sum_{e\in E}\left( g^{(q)}_{e(2)}(\hat{\gamma}(q)-t(q)) \log a_{e(2)}/b_{e(2)} + g^{(q)}_{e(1)} (\hat{\gamma}(q)-t(q))\log a_{e(1)}/b_{e(1)} \right)\\
			&=-\sum_{\hat{e}\in \hat{E}} g_{\hat{e}}^{(q)}(\hat{\gamma}(q)-t(q)) \log a_{\hat{e}}/b_{\hat{e}},
		\end{aligned}
		$$
		which implies that \eqref{s44} holds.
		\vspace{0.2cm}
		
		\textit{Case 2: $\mathcal{G}^{(q)}_{x,\hat{y}(x)}$ is not irreducible.}
		\vspace{0.2cm}
		
		By Lemma \ref{lemma1}, we have $\kappa|_{\hat{E}'}$ is a one-to-one map from $\hat{E}'$ to $\hat{E}''$ and  for each $\hat{e},\hat{e}'\in \hat{E}'$, 
		$$\mathcal{G}^{(q)}_{0,\hat{\gamma}(q)}[\hat{E}',\hat{E}'](\hat{e},\hat{e}')=\mathcal{G}^{(q)}_{\hat{\gamma}(q)-t(q),t(q)}[\hat{E}'',\hat{E}''](\kappa(\hat{e}),\kappa(\hat{e}')).
		$$
		So $\hat{u}_{\kappa(\hat{e})}(\hat{\gamma}(q)-t(q),2)= \hat{u}_{\hat{e}}(0,1)$ and $\hat{v}_{\kappa(\hat{e})}(\hat{\gamma}(q)-t(q),2)=\hat{v}_{\hat{e}}(0,1)$, which gives that $g^{(q)}_{\kappa(\hat{e})}(\hat{\gamma}(q)-t(q),2)= g^{(q)}_{\hat{e}}(0,1)$ for $\hat{e}\in \hat{E}'$. Similarly, $g^{(q)}_{\kappa(\hat{e})}(\hat{\gamma}(q)-t(q),1)=g^{(q)}_{\hat{e}}(0,2)$ for $\hat{e}\in \hat{E}''$. Then noticing that $\log a_{\hat{e}}/ b_{\hat{e}} =-\log a_{\kappa({\hat{e}})} /b_{\kappa(\hat{e})}$, we have  
		$$
		\begin{aligned}
			\sum_{\hat{e}\in \hat{E}} g_{\hat{e}}^{(q)}(0) \log a_{\hat{e}}/b_{\hat{e}}
			&= \sum_{\hat{e}\in \hat{E}'} g_{\hat{e}}^{(q)}(0,1) \log  a_{\hat{e}}/b_{\hat{e}}
			+\sum_{\hat{e}\in \hat{E}''} g_{\hat{e}}^{(q)}(0,2) \log  a_{\hat{e}}/b_{\hat{e}}\\
			&=-\sum_{\hat{e}\in \hat{E}'} g^{(q)}_{\kappa(\hat{e})}(\hat{\gamma}(q)-t(q),2) \log a_{\kappa(\hat{e})}/ b_{\kappa(\hat{e})}-\sum_{\hat{e}\in \hat{E}''}g^{(q)}_{\kappa(\hat{e})}(\hat{\gamma}(q)-t(q),1) \log a_{\kappa(\hat{e})}/ b_{\kappa(\hat{e})} \\
			&=-\sum_{\hat{e}\in \hat{E}} g^{(q)}_{\hat{e}}(\hat{\gamma}(q)-t(q)) \log a_{\hat{e}}/ b_{\hat{e}}.
		\end{aligned}
		$$
		So \eqref{s44} holds.
	\end{proof}
        Analogous to that in Section \ref{sec4}, for $\hat{w}\in \hat{E}^*$, $\hat{e}\in \hat{E}$, denote $\#(\hat{w},\hat{e}) :=\# \{i:\hat{w}_i=\hat{e}, 1\leq i\leq |\hat{w}|\}$ the number of times $\hat{e}$ appears in $\hat{w}$.
	\begin{lemma}\label{lemma5}
		For fixed $x\in \mathbb{R}$, for any positive probability vector $\lambda=(\lambda_{\hat{e}})_{\hat{e}\in \hat{E}}$,
		for $\epsilon>0$, there exists $0<C<1$ independent of $\epsilon$  such that $\sum_{\hat{w}=\hat{w}_1\cdots \hat{w}_k \in \hat{\mathcal{B}}_{x,k}(\epsilon) }  \lambda_{\hat{w}_1} \tilde{\mathcal{G}}_x(\hat{w}_1,\hat{w}_2) \cdots \tilde{\mathcal{G}}_x(\hat{w}_{k-1},\hat{w}_k) <C$ for all large enough $k$, where 
		$$	 		\hat{\mathcal{B}}_{x,k}(\epsilon):=\left \{  \hat{w}\in \hat{E}^k: \sum_{\hat{e}\in \hat{E}}\left| \frac{\#(\hat{w},\hat{e})}{k}-g^{(q)}_{\hat{e}}(x)\right|>\epsilon \right \}.
		$$		
	\end{lemma}
	\begin{proof}
		First we suppose  $\mathcal{G}^{(q)}_{x,\hat{y}(x)}$ is irreducible, so $\tilde{\mathcal{G}}_x $ is irreducible.
		Let $X=\{X_i\}_{i\geq 1}$ be a Markov chain on a finite state space $\hat{E}$ associated with a transition probability matrix $\tilde{\mathcal{G}}_x$, and an invariant distribution $g^{(q)}(x)$. By Proposition \ref{prop2}, using a  same argument in the proof of Lemma \ref{leclt}, the lemma follows.
		
		It remains to prove the case that 
		$\mathcal{G}^{(q)}_{x,\hat{y}(x)}$ is not irreducible. By Lemma \ref{lemma1}, $\mathcal{G}^{(q)}_{x,\hat{y}(x)}[\hat{E}',\hat{E}']$ and $\mathcal{G}^{(q)}_{x,\hat{y}(x)}[\hat{E}'',\hat{E}'']$ are two irreducible matrices. For $k\in \mathbb{N}$, denote $(\hat{E}')^k=\{\hat{w}=\hat{w}_1\cdots \hat{w}_k\in \hat{E}^k: \hat{w_i} \in \hat{E}'\}$ and $(\hat{E}'')^k=\{\hat{w}=\hat{w}_1\cdots \hat{w}_k\in \hat{E}^k: \hat{w_i} \in \hat{E}''\}$ . So $\hat{E}^k =(\hat{E}')^k \cup (\hat{E}'')^k$ by Lemma \ref{lemma1}. Define 
		$$
		\hat{\mathcal{B}}_{x,k}^{(1)}(\epsilon):=\left \{  \hat{w}\in (\hat{E}')^k: \sum_{\hat{e}\in \hat{E}'}\left| \frac{\#(\hat{w},\hat{e})}{k}-g^{(q)}_{\hat{e}}(x,1)\right|>\epsilon \right \},
		$$
		and 
		$$
		\hat{\mathcal{B}}_{x,k}^{(2)}(\epsilon):=\left \{  \hat{w}\in (\hat{E}'')^k: \sum_{\hat{e}\in \hat{E}''}\left| \frac{\#(\hat{w},\hat{e})}{k}-g^{(q)}_{\hat{e}}(x,2)\right|>\epsilon \right \}.
		$$
		By the definition of $g^{(q)}(x)$, $\hat{\mathcal{B}}_{x,k}(\epsilon)=\hat{\mathcal{B}}_{x,k}^{(1)}(\epsilon) \cup \hat{\mathcal{B}}_{x,k}^{(2)}(\epsilon)$.
		Let $X'=\{X'_i\}_{i\geq 1}$ (resp. $X''=\{X''_i\}_{i\geq 1}$) be a Markov chain on a finite state space $\hat{E}'$ (resp. $\hat{E}''$) associated with a transition probability matrix $\tilde{\mathcal{G}}_{x,1}$ (resp. $\tilde{\mathcal{G}}_{x,2}$), and an invariant distribution $g^{(q)}(x,1)$ (resp. $g^{(q)}(x,2)$). Taking an initial distribution $(\lambda_{\hat{e},1})_{\hat{e}\in \hat{E}'}=(\lambda_{\hat{e}}/\sum_{\hat{e}\in \hat{E}'}\lambda_{\hat{e}})_{\hat{e}\in \hat{E}'} $ (resp. $(\lambda_{\hat{e},2})_{\hat{e}\in \hat{E}''}=(\lambda_{\hat{e}}/\sum_{\hat{e}\in \hat{E}''}\lambda_{\hat{e}})_{\hat{e}\in \hat{E}''})$, by Propostion \ref{prop2}, we see that there exists $0<C<1$ such that 
		\begin{equation}\label{s46}
			\sum_{\hat{w}=\hat{w}_1 \cdots \hat{w}_k \in \hat{\mathcal{B}}_{x,k}^{(1)}(\epsilon)} \lambda_{\hat{w}_1,1} \tilde{\mathcal{G}}_{x,1}(\hat{w}_1, \hat{w}_2) \cdots \tilde{\mathcal{G}}_{x,1}(\hat{w}_{k-1}, \hat{w}_k) <C,
		\end{equation}
		and 
		\begin{equation}\label{s47}
			\sum_{\hat{w}=\hat{w}_1 \cdots \hat{w}_k \in \hat{\mathcal{B}}_{x,k}^{(2)}(\epsilon)} \lambda_{\hat{w}_1,2} \tilde{\mathcal{G}}_{x,2}(\hat{w}_1, \hat{w}_2) \cdots \tilde{\mathcal{G}}_{x,2}(\hat{w}_{k-1}, \hat{w}_k) <C. 
		\end{equation}
		Combining \eqref{s74}, \eqref{s46}, \eqref{s47} and the definition of $(\lambda_{\hat{e},1})_{\hat{e}\in \hat{E}'}$,$(\lambda_{\hat{e},2})_{\hat{e}\in \hat{E}''}$, 
		we have 
		$$
		\sum_{\hat{w}=\hat{w}_1 \cdots \hat{w}_k \in \hat{\mathcal{B}}_{x,k}(\epsilon)} \lambda_{\hat{w}_1} \tilde{\mathcal{G}}_{x}(\hat{w}_1, \hat{w}_2) \cdots \tilde{\mathcal{G}}_{x}(\hat{w}_{k-1}, \hat{w}_k) < C(\sum_{\hat{e}\in \hat{E}'} \lambda_{\hat{e}} +\sum_{\hat{e}\in \hat{E}''} \lambda_{\hat{e}}) =C.
		$$
	\end{proof}

 \subsection{Proofs of Theorem \ref{thmain3} and Corollaries \ref{cor3}, \ref{cor4}}

 With all these lemmas in hand, now we come to the proofs of the main results in this section.
	
	\begin{proof}[Proof of Theorem \ref{thmain3}]
		By Lemmas \ref{lemma47}, \ref{lemma2} and Remark \ref{s77}, we see 
		\begin{equation}
			\begin{aligned}\label{s49}
				&\left \{(x,y):\rho(\mathcal{G}^{(q)}_{x,y})=1, \min \{ 0,\hat{\gamma}(q)-t(q)\} \leq x \leq \max \{ 0,\hat{\gamma}(q)-t(q)\} \right\}\\
				=&\left \{(x,y):\rho(\mathcal{G}^{(q)}_{x,y})=1, \min \{ t(q),\hat{\gamma}(q)\} \leq y \leq \max \{ t(q),\hat{\gamma}(q)\} \right \}.
			\end{aligned}
		\end{equation}
		For $w\in E^*,$ Note that 
		\begin{equation}\label{s48}
			\alpha_1(w)= \left \{ \begin{aligned}
				&c_w \quad &\text{ if } c_w\geq d_w,\\
				&d_w &\text{ if } c_w <d_w,
			\end{aligned}\right.
		\end{equation}
		and 
		\begin{equation} \label{s75}
		\tau_w(q)=\left \{ \begin{aligned}
			&\tau_{x,w_k}(q) \quad & &\text{ if } c_w\geq d_w \text{ and }  T_w  \text{ is diagonal, }\\
			&\tau_{y,w_k}(q) & &\text{ if } c_w <d_w \text{ and } T_w \text{ is diagonal, }
            \\
            &\tau_{y,w_k}(q) \quad & &\text{ if } c_w\geq d_w \text{ and }  T_w  \text{ is anti-diagonal, }
            \\	&\tau_{x,w_k}(q) & &\text{ if } c_w <d_w \text{ and } T_w \text{ is anti-diagonal.}\\
		\end{aligned}\right .
		\end{equation}
		We divide the proof into two parts.
		\vspace{0.2cm}
		
		\textit{Part I: When $\hat{\gamma}(q)\leq t(q)$.}  \vspace{0.2cm}
  
		Combining \eqref{s48}, \eqref{s75} and Lemma \ref{lemma3}, for $x,y\in \mathbb{R}$, we have 
  \begin{equation} \label{s78}
		p_w^q \alpha_1(w)^{x+\tau_w(q)} \alpha_2(w)^{y-\tau_w(q)}=\left \{ \begin{aligned}
			&p_{{w}(1)}^q a_{{w}(1)}^{x+\tau_{{w}_k(1)}(q)} b_{{w}(1)}^{y-\tau_{{w}_k(1)}(q)}\quad & &\text{ if } c_w\geq d_w \text{ and }  T_w  \text{ is diagonal, }\\
			&p_{{w}(2)}^q a_{{w}(2)}^{x+\tau_{{w}_k(2)}(q)} b_{{w}(2)}^{y-\tau_{{w}_k(2)}(q)} & &\text{ if } c_w <d_w \text{ and } T_w \text{ is diagonal, }
            \\
            &p_{{w}(2)}^q a_{{w}(2)}^{x+\tau_{{w}_k(2)}(q)} b_{{w}(2)}^{y-\tau_{{w}_k(2)}(q)} \quad & &\text{ if } c_w\geq d_w \text{ and }  T_w  \text{ is anti-diagonal, }
            \\	
            &p_{{w}(1)}^q a_{{w}(1)}^{x+\tau_{{w}_k(1)}(q)} b_{{w}(1)}^{y-\tau_{{w}_k(1)}(q)} & &\text{ if } c_w <d_w \text{ and } T_w \text{ is anti-diagonal.}\\
		\end{aligned}\right .
		\end{equation}
		Recalling the definition of matrix $A_k^{s,q}$ in Section \ref{sec2}, using Lemmas \ref{lemma46}, \ref{lemma49} and  \eqref{s78}, taking $x=0,y=\gamma(q)$,
		we have 
		\begin{equation}
			\begin{aligned}\nonumber
				\Vert A_k^{\gamma(q),q} \Vert 
				&=\sum_{w\in E^k} p_w^q \alpha_1(w)^{\tau_w(q)} \alpha_2(w)^{\gamma(q)-\tau_w(q)}
                \\
				&\leq  \sum_{w\in E^k} \left(p_{{w}{(1)}}^q a_{{w}{(1)}}^{\tau_{{w}_k{(1)}}(q)} b_{{w}{(1)}}^{\gamma(q)-\tau_{{w}_k{(1)}}(q)} + p_{{w}{(2)}}^q a_{{w}{(2)}}^{\tau_{{w}_k{(2)}}(q)} b_{{w}{(2)}}^{\gamma(q)-\tau_{{w}_k{(2)}}(q)} \right)
				\\
	      &=\sum_{\hat{w}\in \hat{E}^k} p_{\hat{w}}^q a_{\hat{w}}^{\tau_{\hat{w}_k}(q)} b_{\hat{w}}^{\gamma(q)-\tau_{\hat{w}_k}(q)}
				\\
				&\asymp \Vert \left (\mathcal{G}^{(q)}_{0, \gamma(q)} \right )^k\Vert.
			\end{aligned}
		\end{equation}
		By the definition of $\gamma(q)$, we have 
		$$
		\begin{aligned}
			1 &= \lim_{k\to \infty} \Vert A_k^{\gamma(q),q} \Vert ^{1/k}
			\leq  \rho\left (\mathcal{G}^{(q)}_{0,\gamma(q)}\right ),
		\end{aligned}
		$$
		which gives that $\gamma(q) \leq \hat{\gamma}(q)$ by Lemma \ref{lemma47}. The upper bound  estimate of $\gamma(q)$ follows. 
  \vspace{0.2cm}

		For any $x \in [\hat{\gamma}(q)-t(q), 0]$, we have $\hat{y}(x)\in [ \hat{\gamma}(q),t(q) ]$. Using \eqref{s76}, we have 
		\begin{equation}\nonumber
			\begin{aligned}
				\Vert (\mathcal{G}^{(q)}_{x,\hat{y}(x)})^k \Vert &\asymp & & \sum_{w\in E^k} \left( p_{{w}{(1)}}^q a_{{w}{(1)}}^{x+\tau_{{w}_k{(1)}}(q)} b_{{w}{(1)}}^{\hat{y}(x)-\tau_{{w}_k{(1)}}(q)} + p_{{w}{(2)}}^q a_{{w}{(2)}}^{x+\tau_{{w}_k{(2)}}(q)} b_{{w}{(2)}}^{\hat{y}(x)-\tau_{{w}_k{(2)}}(q)} \right)
				\\
				&= & & \sum_{w\in E^k:T_w \text{ is diagonal}} \left(p_w^q c_w^{x+\tau_{x,w_k}(q)}d_w^{\hat{y}(x)-\tau_{x,w_k}(q)} +p_w^q d_w^{x+\tau_{y,w_k}(q)}c_w^{\hat{y}(x)-\tau_{y,w_k}(q)}\right)
    \\
    & & &+\sum_{w\in E^k:T_w \text{ is anti-diagonal}} \left(p_w^q d_w^{x+\tau_{x,w_k}(q)}c_w^{\hat{y}(x)-\tau_{x,w_k}(q)} +p_w^q c_w^{x+\tau_{y,w_k}(q)}d_w^{\hat{y}(x)-\tau_{y,w_k}(q)}\right)
    \\
				&= & &\sum_{w\in E^k} p_w^q \alpha_1(w)^{\tau_w(q)} \alpha_2(w)^{x+\hat{y}(x)-\tau_w(q)} \left( \left( \frac{\alpha_1(w)}{\alpha_2(w)} \right) ^x+ \left( \frac{\alpha_1(w)}{\alpha_2(w)} \right)^{\hat{y}(x)-t(q)} \right)
    \\
    &\lesssim & &\sum_{w\in E^k} p_w^q \alpha_1(w)^{\tau_w(q)} \alpha_2(w)^{x+\hat{y}(x)-\tau_w(q)}
    \\
				&= & &\Vert A_k^{x+\hat{y}(x),q} \Vert,
			\end{aligned}
		\end{equation}
  where the second equality follows from a check through $c_w\geq d_w$ or $c_w <d_w$ separately.
    This yields $P(x+\hat{y}(x),q)\geq 1$. It follows that $\gamma(q)\geq x+\hat{y}(x)$ by Lemma \ref{propertyP}, which gives that \begin{equation}\nonumber
			\gamma(q)\geq \max \{x+\hat{y}(x): \hat{\gamma}(q)-t(q)\leq x \leq 0\},
		\end{equation}
		a lower bound estimate of $\gamma(q)$.
		Combining  this with \eqref{s49} and the upper bound estimate $\gamma(q)\leq \hat{\gamma}(q)$, we obtain  \eqref{s50}.
		
		\vspace{0.2cm}

		\textit{Part II: When $\hat{\gamma}(q)>t(q)$.}  \vspace{0.2cm}
		
		We firstly prove that $\gamma(q)\leq \min \{x+\hat{y}(x): 0\leq x\leq \hat{\gamma}(q)-t(q)\}.$
		For $x\in [0,\hat{\gamma}(q)-t(q)]$, we have $\hat{y}(x)\in [t(q), \hat{\gamma}(q)]$. Note that 
		\begin{equation}\nonumber
			\begin{aligned}
				\Vert A_k^{\gamma(q),q} \Vert &= & & \sum_{w\in E^k:c_w \geq d_w} p_w^q c_w^{\tau_w(q)} d_w^{\gamma(q)-\tau_w(q)} +\sum_{w\in E^k:c_w<d_w} p_w^q d_w^{\tau_w(q)} c_w^{\gamma(q)-\tau_w(q)}\\
				&\leq& & \sum_{w\in E^k:c_w\geq d_w}p_w^q c_w^{x+\tau_w(q)} d_w^{\gamma(q)-x-\tau_w(q)} +\sum_{w\in E^k:c_w<d_w} p_w^q c_w^{\gamma(q)-\hat{y}(x)+t(q)-\tau_w(q)} d_w^{\hat{y}(x)+\tau_w(q)-t(q)}\\
				&\leq & &\sum_{w\in E^k:c_w\geq d_w} \left( p_{{w}{(1)}}^q a_{{w}{(1)}}^{x+\tau_{{w}_k{(1)}}(q)} b_{{w}{(1)}}^{\gamma(q)-x-\tau_{{w}_k{(1)}}(q)}+p_{{w}{(2)}}^q a_{{w}{(2)}}^{x+\tau_{{w}_k{(2)}}(q)}b_{{w}{(2)}}^{\gamma(q)-x-\tau_{{w}_k{(2)}}(q)} \right)
    \\
        & & &+\sum_{w\in E^k:c_w<d_w} \left(p_{{w}{(1)}}^q a_{{w}{(1)}}^{\gamma(q)-\hat{y}(x)+\tau_{{w}_k{(1)}}(q)} b_{{w}{(1)}}^{\hat{y}(x)-\tau_{{w}_k{(1)}}(q)}+p_{{w}{(2)}}^q a_{{w}{(2)}}^{\gamma(q)-\hat{y}(x)+\tau_{{w}_k{(2)}}(q)} b_{{w}{(2)}}^{\hat{y}(x)-\tau_{{w}_k{(2)}}(q)}\right)
				\\
				&\lesssim& & \max \left \{\Vert \left(\mathcal{G}^{(q)}_{x,\gamma(q)-x}\right)^k \Vert ,\Vert \left (\mathcal{G}^{(q)}_{\gamma(q)-\hat{y}(x),\hat{y}(x)}\right)^k \Vert \right \},
			\end{aligned}
		\end{equation}
  where the second inequality follows from a  check through $T_w$ is diagonal or not separately. Then
		it is easy to see that $\max\left \{\rho\left (\mathcal{G}^{(q)}_{x,\gamma(q)-x}\right), \rho\left (\mathcal{G}^{(q)}_{\gamma(q)-\hat{y}(x),\hat{y}(x)}\right)\right\} \geq 1$, which gives that $\gamma(q)\leq x+\hat{y}(x)$. Thus \begin{equation}\label{s51}
			\gamma(q)\leq \min \{x+\hat{y}(x):0\leq x \leq \hat{\gamma}(q)-t(q)\},
		\end{equation}
  an upper bound estimate of $\gamma(q)$.
  \vspace{0.2cm}
		
		By \eqref{s49}, it remains to prove the reverse inequality of $\eqref{s51}$.
        Recall that by Lemma \ref{lemma4}, there exists $x\in [0,\hat{\gamma}(q)-t(q)]$ such that 
		\begin{equation}\label{s52}
			\sum_{\hat{e}\in \hat{E}} g^{(q)}_{\hat{e}}(x) \log a_{\hat{e}}/b_{\hat{e}}=0.
		\end{equation}
	Fix this  $x$.
		Noticing that $x+\hat{y}(x)\geq t(q)$,  by \eqref{s78}, we always have 
		\begin{equation}\label{s54}
			p_w^q \alpha_1(w)^{\tau_w(q)} \alpha_2(w)^{x+\hat{y}(x)-\tau_w(q)} = \min \left \{ p_{{w}{(1)}}^q a_{{w}{(1)}}^{\tau_{{w}_k{(1)}}(q)} b_{{w}{(1)}}^{x+\hat{y}(x)-\tau_{{w}_k{(1)}}(q)}, p_{{w}{(2)}}^q a_{{w}{(2)}}^{\tau_{{w}_k{(2)}}(q)} b_{{w}{(2)}}^{x+\hat{y}(x)-\tau_{{w}_k{(2)}}(q)}\right \}.
		\end{equation}
		Therefore, by Lemma  \ref{lemma3}, \eqref{s54} and noticing that  for $\hat{w}\in \hat{E}^*$, $a_{\kappa(\hat{w})}=b_{\hat{w}}$, $b_{\kappa(\hat{w})}=a_{\hat{w}}$, we have 
		\begin{equation}
			\label{s57}
			\begin{aligned}
				\Vert A_k^{x+\hat{y}(x),q} \Vert
				&=\sum_{w\in E^k} p_w^q \alpha_1(w)^{\tau_w(q)} \alpha_2(w)^{x+\hat{y}(x)-\tau_w(q)}
				\\
				&= 
				\sum_{w\in E^k} \min \left \{ p_{{w}{(1)}}^q a_{{w}{(1)}}^{\tau_{{w}_k{(1)}}(q)} b_{{w}{(1)}}^{x+\hat{y}(x)-\tau_{{w}_k{(1)}}(q)}, p_{{w}{(2)}}^q a_{{w}{(2)}}^{\tau_{{w}_k{(2)}}(q)} b_{{w}{(2)}}^{x+\hat{y}(x)-\tau_{{w}_k{(2)}}(q)}\right \}\\
				&\asymp \sum_{\hat{w}\in \hat{E}^k} \min \left \{ p_{\hat{w}}^q a_{\hat{w}}^{\tau_{\hat{w}_k}(q)} b_{\hat{w}}^{x+\hat{y}(x)-\tau_{\hat{w}_k}(q)}, p_{\kappa(\hat{w})}^q a_{\kappa(\hat{w})}^{\tau_{\kappa(\hat{w}_k)}(q)} b_{\kappa(\hat{w})}^{x+\hat{y}(x)-\tau_{\kappa(\hat{w}_k)}(q)} \right \}\\
				&= \sum_{\hat{w}\in \hat{E}^k :a_{\hat{w}}\geq b_{\hat{w}}} p_{\hat{w}}^q a_{\hat{w}}^{\tau_{\hat{w}_k}(q)} b_{\hat{w}}^{x+\hat{y}(x)-\tau_{\hat{w}_k}(q)} +\sum_{\hat{w}\in \hat{E}^k :a_{\hat{w}}< b_{\hat{w}}} p_{\hat{w}}^q b_{\hat{w}}^{\tau_{\kappa(\hat{w}_k)}(q)} a_{\hat{w}}^{x+\hat{y}(x)-\tau_{\kappa(\hat{w}_k)}(q)}.
			\end{aligned}
		\end{equation}
		
		For $\hat{w}\in \hat{E}^k$, note that 
		$$
		\frac{\log a_{\hat{w}}}{k}= \sum_{\hat{e}\in \hat{E}}\frac{\#(\hat{w},\hat{e})}{k} \log a_{\hat{e}}
		$$
		and 
		$$
		\frac{\log b_{\hat{w}}}{k}= \sum_{\hat{e}\in \hat{E}}\frac{\#(\hat{w},\hat{e})}{k} \log b_{\hat{e}}.
		$$
		For $\epsilon>0$, let $\eta = - \epsilon \sum_{\hat{e}\in \hat{E}} \log a_{\hat{e}}b_{\hat{e}}$.
		For large enough $k$,  for $\hat{w}\in \hat{E}^k\setminus \hat{\mathcal{B}}_{x,k}(\epsilon)$, using  \eqref{s52}, we have 
		$$
		-\eta \leq  \frac{\log a_{\hat{w}}}{k} -\frac{\log b_{\hat{w}}}{k} \leq \eta,
		$$
		which gives that 
		\begin{equation}\label{s55}
			e^{-\eta k } \leq \frac{a_{\hat{w}}}{b_{\hat{w}}} \leq e^{ \eta k}.
		\end{equation}
		By Lemma \ref{lemma5}, picking $\lambda=(\frac{1}{\#\hat{E}})_{\hat{e}\in \hat{E}}$, using Lemma \ref{lemma6} and \eqref{s72}, we can see 
		\begin{equation}
			\label{s56}
			\begin{aligned}
				&\sum_{\hat{w}=\hat{w}_1\cdots \hat{w}_k \in \hat{E}^k \setminus \hat{\mathcal{B}}_{x,k}(\epsilon)} \frac{1}{\#\hat{E}} p_{\hat{w}_2 \cdots \hat{w}_k}^q a_{\hat{w}_2 \cdots \hat{w}_k}^{x+\tau_{\hat{w}_k}(q)} b_{\hat{w}_2 \cdots \hat{w}_k}^{\hat{y}(x)-\tau_{\hat{w}_k}(q)} \frac{\hat{u}_{\hat{w}_k}(x)}{\hat{u}_{\hat{w}_1}(x)}
				\\		=&\sum_{\hat{w}=\hat{w}_1\cdots \hat{w}_k \in \hat{E}^k \setminus \hat{\mathcal{B}}_{x,k}(\epsilon)} \lambda_{\hat{w}_1} \tilde{\mathcal{G}}_x(\hat{w}_1,\hat{w}_2)\cdots \tilde{\mathcal{G}}_x(\hat{w}_{k-1},\hat{w}_k)
				\geq 1-C>0,
			\end{aligned}
		\end{equation}
		for some $C<1$.
		Therefore, using \eqref{s57},  \eqref{s55}  and \eqref{s56}, we have 
		\begin{equation}
			\nonumber
			\begin{aligned}
				\Vert A_k^{x+\hat{y}(x),q} \Vert
				&\gtrsim & &\sum_{\hat{w}\in \hat{E}^k \setminus \hat{\mathcal{B}}_{x,k}(\epsilon):a_{\hat{w}}\geq b_{\hat{w}}} p_{\hat{w}}^q a_{\hat{w}}^{\tau_{\hat{w}_k}(q)} b_{\hat{w}}^{x+\hat{y}(x)-\tau_{\hat{w}_k}(q)} +\sum_{\hat{w}\in \hat{E}^k \setminus \hat{\mathcal{B}}_{x,k}(\epsilon):a_{\hat{w}}< b_{\hat{w}}} p_{\hat{w}}^q b_{\hat{w}}^{\tau_{\kappa(\hat{w}_k)}(q)} a_{\hat{w}}^{x+\hat{y}(x)-\tau_{\kappa(\hat{w}_k)}(q)}
				\\
				&= & &\sum_{\hat{w}\in \hat{E}^k\setminus \hat{\mathcal{B}}_{x,k}(\epsilon):a_{\hat{w}}\geq b_{\hat{w}}} p_{\hat{w}}^q a_{\hat{w}}^{x+\tau_{\hat{w}_k}(q)} b_{\hat{w}}^{\hat{y}(x)-\tau_{\hat{w}_k}(q)} \left(\frac{a_{\hat{w}}}{b_{\hat{w}}}\right)^{-x}\\
    & & &
				+\sum_{\hat{w}\in \hat{E}^k\setminus \hat{\mathcal{B}}_{x,k}(\epsilon):a_{\hat{w}}< b_{\hat{w}}} p_{\hat{w}}^q b_{\hat{w}}^{\hat{y}(x)-\tau_{\hat{w}_k}(q)} a_{\hat{w}}^{x+\tau_{\hat{w}_k}(q)} \left(\frac{b_{\hat{w}}}{a_{\hat{w}}} \right)^{t(q)-\hat{y}(x)}\\
				&\geq & &\sum_{\hat{w}\in\hat{E}^k\setminus \hat{\mathcal{B}}_{x,k}(\epsilon)} p_{\hat{w}}^q a_{\hat{w}}^{x+\tau_{\hat{w}_k}(q)} b_{\hat{w}}^{\hat{y}(x)-\tau_{\hat{w}_k}(q)}  \min \left \{e^{\eta k(-x)},e^{\eta k
					(t(q)-\hat{y}(x))} \right \}\\
				&\geq & &\sum_{\hat{w}\in \hat{E}^k\setminus \hat{\mathcal{B}}_{x,k}(\epsilon)} \frac{1}{\#\hat{E}}p_{\hat{w}_2\cdots \hat{w}_k}^q a_{\hat{w}_2 \cdots \hat{w}_k}^{x+\tau_{\hat{w}_k}(q)} b_{\hat{w}_2 \cdots \hat{w}_k}^{\hat{y}(x)-\tau_{\hat{w}_k}(q)}  \frac{\hat{u}_{\hat{w}_k}(x)}{\hat{u}_{\hat{w}_1}(x)}\min \left \{e^{\eta k(-x)},e^{\eta k
					(t(q)-\hat{y}(x))} \right \} C'\\
				&\geq & & (1-C)	\min \left \{e^{\eta k(-x)},e^{\eta k
					(t(q)-\hat{y}(x))} \right \} C',
			\end{aligned}
		\end{equation}
		where $$
  \begin{aligned}
  C'=(\#\hat{E}) p_*^q \min \big \{ \alpha_*^{x+\hat{y}(x)}, &\alpha_*^{x+\tau_A(q)} \alpha^{*(\hat{y}(x)-\tau_A(q))},\alpha_*^{x+\tau_B(q)} \alpha^{*(\hat{y}(x)-\tau_B(q))},\\ &\alpha^{*(x+\tau_A(q))} \alpha_*^{\hat{y}(x)-\tau_A(q)},
  \alpha^{*(x+\tau_B(q))} \alpha_*^{\hat{y}(x)-\tau_B(q)}, \alpha^{*(x+\hat{y}(x)) } \big \} \cdot \min_{\hat{e},\hat{e}'\in \hat{E}} \frac{\hat{u}_{\hat{e}}(x)}{\hat{u}_{\hat{e}'}(x)} 
  \end{aligned}
  $$ is a positive number.
		Thus, 
		$$
		P(x+\hat{y}(x),q) \geq \min \left \{e^{- \eta x},e^{\eta 
			(t(q)-\hat{y}(x))} \right \}.
		$$
		Letting $\epsilon \to 0$,   $\eta \to0 $ gives $P(x+\hat{y}(x),q)\geq 1$. This yields $\gamma(q)\geq x+\hat{y}(x)$ by Lemma \ref{propertyP}, a lower bound estimate of $\gamma(q)$. Combining this with \eqref{s49} and \eqref{s51}, we have  $\gamma(q)= x+\hat{y}(x)$ and \eqref{s58} holds. 
		
	\end{proof}
	
	\begin{proof}[Proof of Corollary \ref{cor3}]
		The proof is analogous to that of Corollary \ref{cor2}.
		
		By Theorem \ref{thmain3},
		it suffices to prove that 
		$$
		\hat{\gamma}(0) \leq    t(0).
		$$
		
		Suppose that  $\hat{\gamma}(0)>t(0)$. Using Theorem \ref{thmain3}-(b), we have 
		$$
		\gamma(0)=\min \{x+\hat{y}(x):0\leq x \leq \hat{\gamma}(0)-t(0)\}.
		$$
		By the product formula, we have 
		$
		\gamma(0) \leq t(0).
		$
		So there exists a $x \in (0,\hat{\gamma}(0)-t(0))$ such that $\hat{y}(x)\leq t(0)-x<t(0)$. However, $\hat{y}(\hat{\gamma}(0)-t(0))=t(0)$, a contradiction to the fact that  $\hat{y}(x)$ is decreasing by Lemma \ref{lemma2}.
	\end{proof}

	\begin{proof}[Proof of Corollary \ref{cor4}]
		It suffice to prove that 
		\begin{equation}\label{s60}
			\rho(\mathcal{G}^{(q)}_{x,y})=\rho(H^{(q)}_{x,y}).
		\end{equation}
		At this time, for $1\leq i,j\leq N$,
		$$
		{G}^{(q)}_{x,y,i,j} =\left \{
		\begin{aligned}
			&\left(\begin{array}{cc}
				p_j^q a_{j}^{x+\tau_{\mu^x}(q)} b_{j}^{y-\tau_{\mu^x}(q)} & 0 \\
				0   & p_j^q b_{j}^{x+\tau_{\mu^y}(q)} a_{j}^{y-\tau_{\mu^y}(q)} \\	\end{array}\right) &\quad  &\text{ if }  1\leq j<k, \\
			&\left(\begin{array}{cc}
				0 & p_j^q a_j^{x+\tau_{\mu^y}(q)} b_{j}^{y-\tau_{\mu^y}(q)} \\ p_j^q
				b_{j}^{x+\tau_{\mu^x}(q)} a_{j}^{y-\tau_{\mu^x}(q)}   &0 \\	\end{array}\right)  &\quad&\text{ if }  k\leq j\leq N.
		\end{aligned} \right .
		$$
		So by the definition of $H^{(q)}_{x,y}$ in \eqref{s61}, we have $$H^{(q)}_{x,y}=\sum_{j=1}^N G^{(q)}_{x,y,i,j}, \quad \text{ for all } 1\leq i\leq N.
        $$
        
		Note that for two non-negative matrices $A$ and $A'$, $\Vert A+A' \Vert =\Vert A \Vert +\Vert A' \Vert $.
		Using Gelfand formula, for any $1\leq i\leq N$, we have 
		$$
		\begin{aligned}
			\rho(\mathcal{G}^{(q)}_{x,y}) &=\lim_{k\to \infty} \left(\sum_{j=1}^N \Vert  \sum_{l_1, \cdots, l_{k-1} =1}^N  G^{(q)}_{x,y,i,l_1} \cdots G^{(q)}_{x,y,l_{k-1},j} \Vert  \right) ^  {1/k}\\
			&=\lim_{k\to \infty} \big \Vert \sum_{l_1, \cdots, l_k =1}^N  G^{(q)}_{x,y,i,l_1} \cdots G^{(q)}_{x,y,i,l_k} \big \Vert ^  {1/k} \\
			&=\lim_{k\to \infty} \big \Vert 
			\big ( \sum_{l=1}^N G^{(q)}_{x,y,i,l} \big )^k \big \Vert ^{1/k}
			\\
			&=\rho(H^{(q)}_{x,y}).
		\end{aligned}
		$$
		Therefore, \eqref{s60} holds.
	\end{proof}

	\section{Proof of Theorem \ref{thmain1}}\label{sec3}
	In this section, we are going to prove Theorem \ref{thmain1}. The main idea is basing on Fraser's work \cite{F16} for the IFS setting.  Let $(V,E,\Psi)$, $\mathcal{P}$ and $\{\mu_v\}_{v\in V}$ be same as before. Let  $\gamma(q)$ be the function determined  by $P(\gamma(q),q)=1$ (see Remark \ref{re1}), where $P$ is the pressure function introduced in Lemma \ref{lemd}. 
	
	Denote the collection of all \textit{infinite admissible words} by
	$$
	E^\infty=\{\omega=\omega_1\omega_2\cdots:t(\omega_{i-1})=i(\omega_i), \forall i>1\}.
	$$
	For $w=w_1\cdots w_k\in E^*$, denote $[w]=\{\omega=\omega_1 \omega_2 \cdots \in E^\infty:\omega_i=w_i \text{ for } 1\leq i\leq k\}$ the \textit{cylinder set} of $w$.
	For $\delta>0$ , write
	$$
	E^*_\delta=\{w=w_1 \cdots w_k \in E^*: \alpha_2(w)<\delta \leq \alpha_2(w_1 \cdots w_{k-1})\}.
	$$
	Roughly speaking, $E^*_\delta$ consists of all $w\in E^*$ for which the shortest side of the rectangle $\psi_w([0,1]^2)$ is comparable to $\delta$. So  for each $w=w_1 \cdots w_k\in E^*_\delta$, we have 
	\begin{equation}\label{t5}
		\delta>\alpha_2(w) \geq \alpha_2(w_1\cdots w_{k-1}) \alpha_2(w_k) \geq  \alpha_* \delta,
	\end{equation}
	where $\alpha_*$ is defined in \eqref{s16}.
	It is easy to see that $E^*_\delta$ is a finite partition of $E^\infty$, i.e. $\#E_\delta^*<\infty$ and
	$$
	E^\infty=\bigcup_{w\in E^*_\delta} [w],
	$$
	where the union is disjoint.
	
	Before proving Theorem \ref{thmain1}, we need the following lemma, which is adapted from \cite[Lemma 7.1]{F16}.
	
	\begin{lemma}\label{lem31}
		Let $q\geq 0$ and $\delta>0$.
		
		(a). For $s> \gamma(q)$,
		$$
		\sum_{w\in E^*_\delta} \varphi^{s,q}(w) \lesssim_{s,q} 1.
		$$
		
		(b). For $s<\gamma(q)$,
		$$
		\sum_{w\in E^*_\delta} \varphi^{s,q}(w) \gtrsim_{s,q} 1.
		$$
	\end{lemma}
	\begin{proof}
		(a). For $s>\gamma(q)$, it follows that 
		$$
		\sum_{w\in E^*_\delta}\varphi^{s,q}(w)\leq \sum_{w\in E^*}\varphi^{s,q}(w) =\sum_{k\geq 1}\sum_{w\in E^k} \varphi^{s,q}(w)
		=\sum_{k\geq 1} \Vert A^{s,q}_k \Vert <\infty,
		$$
		since $\lim_{k\to \infty} \Vert A_k^{s,q} \Vert ^{1/k} =P(s,q)<1$ by Lemma \ref{propertyP}.
		\vspace{0.2cm}
		
		(b). We divide the proof into two cases, $s\leq t(q)$ or $s>t(q)$. 
		\vspace{0.2cm}
		
		\textit{Case 1: $s\leq t(q).$} \vspace{0.2cm}
		
		We will prove $\sum_{w\in E^*_\delta} \varphi^{s,q}(w)\geq 1$ through a 
		contradiction argument. Suppose  
		\begin{equation}\label{s1}
			\sum_{w\in E^*_\delta} \varphi^{s,q}(w)< 1.
		\end{equation}
		
		Fix  a $\delta>0$. Note that for all $w\in E^*$, by Lemma \ref{lem23} and \eqref{s1},  we have 
		\begin{equation}\label{s2}
			\sum_{w'\in E^*_\delta:t(w)=i(w')} \varphi^{s,q}(ww') \leq \varphi^{s,q}(w) \sum_{w'\in E^*_\delta} \varphi^{s,q}(w')<\varphi^{s,q}(w).
		\end{equation}
		Now for large enough $k\in \mathbb{N}$, define 
		$$
		E^k_\delta =\{w=w^{(1)} \cdots w^{(m)}\in E^*:w^{(j)}\in E^*_\delta, \text{ with } |w|>k \text{ and } |w^{(1)}\cdots w^{(m-1)}|\leq k\}.
		$$
		Clearly, $E^k_\delta$ is a finite partition of $E^\infty$. Repeatedly using \eqref{s2}, for large enough $k$, we have 
		\begin{equation}\label{t4}
			\sum_{w\in E^k_\delta} \varphi^{s,q}(w)<1.
		\end{equation}
		
		Let $n_0(\delta)=\max\{|w|:w\in E^*_\delta\}$. For $w\in E^{k+n_0(\delta)}$, we have $w=w^{(1)}w^{(2)}$ for some $w^{(1)}\in E^k_\delta$, $w^{(2)} \in E^*\cup \{\emptyset\}$ with $t(w^{(1)})=i(w^{(2)})$ and $|w^{(2)}|< n_0(\delta)$. Let $c(\delta)=\max\{\varphi^{s,q}(w):|w|<n_0(\delta)\}>0$ which is independent of $k$. Since there are at most $(\#E)^{n_0(\delta)}$ choices of $w^{(2)}\in E^*\cup \{\emptyset\}$, by Lemma \ref{lem23}  and using \eqref{t4}, we have
		$$
		\Vert A^{s,q}_{k+n_0(\delta)} \Vert =\sum_{w\in E^{k+n_0(\delta)}}\varphi^{s,q}(w)\leq  (\#E)^{n_0(\delta)} c_\delta \sum_{w\in E^k_\delta} \varphi^{s,q}(w) <(\#E)^{n_0(\delta)} c_\delta.
		$$
		It follows that 
		$$
		P(s,q)=\lim_{k\to \infty}\Vert A_k^{s,q} \Vert^{1/k} \leq 1.
		$$
		which gives that $s\geq \gamma(q)$ by Lemma \ref{propertyP}, a contradiction.
		\vspace{0.2cm}
		
		\textit{Case 2: $s>t(q)$.} \vspace{0.2cm}
		
		Noticing that  the directed graph $(V,E)$ is strongly connected, for each pair $v,v'\in V$,  define $L(v,v'):=\min \{|w|:w\in E^*,v \stackrel{w}{\rightarrow} v'\}$ and let $L:=\max \{L(v,v'):v,v'\in V\}$ be the maximal length of shortest paths between any two vertices.
		
		Let $C=\min\{\varphi^{s,q}(w):|w|\leq L\}>0.$ Replace the  matrix norm $\Vert \cdot \Vert $ with the \textit{maximum row sum  norm} $\Vert \cdot \Vert_1$, i.e. $\Vert A \Vert_1=\max_{1\leq i\leq N}\sum_{1\leq j\leq N} |a_{ij}|$ for a $N\times N$ matrix $A=(a_{ij})_{1\leq i,j\leq N}$. Due to the equivalence of matrix norms and  $s<\gamma(q)$, it follows that $\Vert A_k^{s,q} \Vert_1 \to \infty$ as $k\to \infty$. So we could find  $k\in \mathbb{N}$ and  $v_0\in V$  such that 
		\begin{equation}\label{s3}
			\sum_{w\in E^k:i(w)=v_0} \varphi^{s,q}(w)=\Vert A_k^{s,q}\Vert_1 >1/C.
		\end{equation}
		Fix such $k$ and $v_0$,
		for small enough $\delta>0$, let 
		\begin{align}
			E_{k,\delta} = \big\{ w^{(1)} \tilde{w}^{(1)} w^{(2)} & \tilde{w}^{(2)}\cdots \tilde{w}^{(m-1)} w^{(m)}\in E^*: |w^{(j)}| =k, i(w^{(j)})=v_0, t(w^{(j)}) \stackrel{ \tilde{w}^{(j)}}{\longrightarrow} v_0,  \nonumber\\		
			& |\tilde{w}^{(j)} |=L(t(w^{(j)}),v_0) ,\alpha_2(w^{(1)} \tilde{w}^{(1)}\cdots w^{(m)}) <\delta \leq \alpha_2(w^{(1)} \tilde{w}^{(1)}\cdots w^{(m-1)})\big\}.\nonumber
		\end{align}
		We can directly check that the cylinder sets of elements in $E_{k,\delta}$ are all disjoint (but $E_{k,\delta}$ is not a finite partition of $E^\infty$).
		For any $w\in E^*$, by Lemma \ref{lem23} and using \eqref{s3}, we have 
		\begin{equation}\label{s4}
			\sum_{w'\in E^{L(t(w),v_0)}} \sum_{w''\in E^{k}:i(w'')=v_0} \varphi^{s,q}(ww'w'')\geq C \varphi^{s,q}(w) \sum_{w''\in E^k:i(w'')=v_0} \varphi^{s,q}(w'')\\
			>\varphi^{s,q}(w).
		\end{equation}
		Repeatedly using \eqref{s4}, we have 
		\begin{equation}\label{t6}
			\sum_{w\in E_{k,\delta}} \varphi^{s,q}(w) >1/C.
		\end{equation}
		
		Also, note that for $w=w^{(1)} \tilde{w}^{(1)} \cdots w^{(m)}\in E_{k,\delta}$, 
		\begin{equation}\label{s5}
			\delta >\alpha_2(w) \geq \alpha_2(w^{(1)} \tilde{w}^{(1)}\cdots w^{(m-1)}) \alpha_2(\tilde{w}^{(m-1)} w^{(m)}) \geq \alpha_*^{L+k} \delta.
		\end{equation}
		Take $\delta'=\alpha_*^{L+k} \delta$, then for $w\in E_{\delta'}^*$, by \eqref{t5} we have
		\begin{equation}\label{s6}
			\alpha_*^{L+k} \delta =\delta' >\alpha_2(w) \geq \alpha_* \delta' =\alpha_*^{L+k+1} \delta.
		\end{equation}
		For $w \in E^*_{\delta'}$,
		if $w$ has a prefix in $E_{k,\delta}$, we can write $w=w'w''$ with $t(w')=i(w'')$ for some $w'\in E_{k,\delta}$ and $w'' \in E^*$. Combining \eqref{s5}, \eqref{s6} and $\alpha_2(w'w'') \leq \alpha_2(w')\alpha_1(w'') \leq \alpha_2(w') \alpha^{*|w''|}$, we have 
		$$
		\alpha_2(w') <\delta \alpha_*^{L+k+1}  \alpha_*^{-L-k-1} \leq \alpha_2(w) \alpha_*^{-L-k-1} \leq \alpha_2(w') \alpha^{*|w''|} \alpha_*^{-L-k-1},
		$$
		which gives that $|w''|\leq (L+k+1)\frac{\log \alpha_*}{\log \alpha^*}$. 
		
		Let $C'=\min \{\varphi^{s,q}(w):|w|\leq(L+k+1)\frac{\log \alpha_*}{\log \alpha^*} \}>0$, then by \eqref{t6} and Lemma \ref{lem23}, we have 
		\begin{align}
			\sum_{w\in E^*_{\delta'}} \varphi^{s,q}(w) &\geq \sum_{w\in E^*_{\delta'}:w \text{ has a prefix in } E_{k,\delta}} \varphi^{s,q}(w) \nonumber\\
			&\geq C' \sum_{w\in E_{k,\delta}} \varphi^{s,q}(w) >C'/C.\nonumber
		\end{align}
		The constant $C'/C$ only depends on $k$ and the choice of $k$ depends on $s,q$.
	\end{proof}
	
	Before proceeding, we mention a fact \cite[Lemma 4.1]{F16} that will be used in the proof of Theorem \ref{thmain1}.
	\begin{equation}\label{s7}
		\text{For }k\in \mathbb{N}, a_1,\cdots, a_k\geq 0 \text{ and } q\geq 0,
		\left( \sum_{i=1}^k a_i\right)^q \asymp_{k,q} \sum_{i=1}^{k} a_i^q.
	\end{equation}
	\begin{proof}[Proof of Theorem \ref{thmain1}]
		
		For $q\geq 0$ and $\delta >0$, recall that we use $\mathcal{M}_\delta$ to denote the collection of $\delta$-mesh on $\mathbb{R}^2$, and for a measure $\mu$ we write $\mathcal{D}_\delta^q(\mu)=\sum_{Q\in \mathcal{M}_\delta} \mu(Q)^q$.
		\vspace{0.2cm}
		
		First of all, due to ROSC, there exists $M>0$ independent of $\delta$ such that for each $Q\in \mathcal{M}_\delta$,
		we have 
		\begin{equation}\label{s30}
			\#\{w\in E^*_\delta:\mu_{t(w)}\circ \psi^{-1}_w \left( Q\cap \psi_w([0,1]^2)\right)>0\} \leq M.
		\end{equation}
		Using this and \eqref{s7}, for $v\in V$, we have 
		\begin{equation}
			\begin{aligned}\label{s8}
				\mathcal{D}_\delta^q(\mu_v)&=\sum_{Q\in \mathcal{M}_\delta} \mu_v(Q)^q =\sum_{Q\in \mathcal{M}_\delta} \left( \sum_{w\in E^*_\delta:i(w)=v} p_w \mu_{t(w)}\circ \psi_w^{-1}\big(Q\cap \psi_w([0,1]^2)\big)\right)^q\\
				&\asymp_{q} \sum_{Q\in \mathcal{M}_\delta}  \sum_{w\in E^*_\delta:i(w)=v} p_w^q \mu_{t(w)}\circ \psi_w^{-1}\big(Q\cap \psi_w([0,1]^2)\big)^q\\
				&=\sum_{w\in E^*_\delta:i(w)=v} p_w^q \mathcal{D}_\delta^q(\mu_{t(w)}\circ \psi_w^{-1}).
			\end{aligned}
		\end{equation}
		Using \eqref{s7} again and the definition of $\pi_w$ in \eqref{s10}, we have 
		\begin{equation}\nonumber
			\mathcal{D}_\delta^q(\mu_{t(w)}\circ\psi_w^{-1})\asymp_q \mathcal{D}^q_{\delta/\alpha_1(w)} \left(\pi_w(\mu_{t(w)})\right).
		\end{equation}
		So equation \eqref{s8} becomes
		\begin{equation}\label{s12}
			\mathcal{D}_\delta^q(\mu_v) \asymp_q \sum_{w\in E^*_\delta:i(w)=v} p_w^q  \mathcal{D}^q_{\delta/\alpha_1(w)} \left(\pi_w(\mu_{t(w)})\right).
		\end{equation}
		\vspace{0.2cm}
		
		On the other hand, recall that for $q\geq 0$, $w\in E^*$, $\tau_w(q)=\tau_{\pi_w(\mu_{t(w)})}(q)$ and equals to either $\tau_A(q)$ or $\tau_B(q)$. By the definition of $L^q$-spectra, for all $w\in E^*$, $\epsilon>0$, $q\geq 0$, small enough $\delta>0$, we have 
		\begin{equation}\label{s13}
			\delta^{-\tau_w(q)+\epsilon/2} \lesssim_{\epsilon,q} \mathcal{D}_\delta^q(\pi_w\left(\mu_{t(w)})\right) \lesssim_{\epsilon,q} \delta^{-\tau_w(q)-\epsilon/2}.
		\end{equation}
		
		For each fixed $v\in V$, 
		combining \eqref{s12}, \eqref{s13} and Lemma \ref{lem31}, noticing that $\delta\asymp \alpha_2(w)$ for $w\in E_\delta^*$, we have 
		\begin{align}
			\delta^{\gamma(q)+\epsilon}\mathcal{D}_\delta^q(\mu_v) &\asymp_q \delta^{\gamma(q)+\epsilon} \sum_{w\in E^*_\delta:i(w)=v} p_w^q  \mathcal{D}^q_{\delta/\alpha_1(w)} \left(\pi_w(\mu_{t(w)})\right)\nonumber\\
			&\lesssim_{\epsilon,q} \delta^{\gamma(q)+\epsilon} \sum_{w\in E^*_\delta:i(w)=v} p_w^q \left( \frac{\delta}{\alpha_1(w)} \right) ^{-\tau_w(q)-\epsilon/2} \nonumber \\
			&\lesssim_{\epsilon,q}  \sum_{w\in E^*_\delta:i(w)=v} p_w^q \alpha_1(w)^{\tau_w(q)} \alpha_2(w)^{\gamma(q)+{\epsilon}/{2}-\tau_w(q)} \cdot \alpha_1(w)^{{\epsilon}/{2}} \nonumber\\
			&\leq \sum_{w\in E^*_\delta:i(w)=v} \varphi^{\gamma(q)+{\epsilon}/{2},q}(w)\nonumber\\
			&\lesssim_{\epsilon,q} 1. \nonumber
		\end{align}
		Letting $\delta,\epsilon\to 0$, we have $\overline{\tau}_{\mu_v}(q)\leq \gamma(q)$.\vspace{0.2cm}
		
		A similar argument will imply $\gamma(q) \leq \underline{\tau}_{\mu_v}(q)$. Indeed, noting that for any $v'\in V$, it always holds  $v \stackrel{w}{\rightarrow} v'$ for some $w$ and $\mu_v \geq p_w \mu_{v'}\circ\psi_w^{-1}$, we have $ \mathcal{D}_\delta^q(\mu_v)\gtrsim_q \mathcal{D}_\delta^q(\mu_{v'})$. Using this, still by \eqref{s12}, \eqref{s13} and Lemma \ref{lem31}, we have 
		
		\begin{align}
			\delta^{\gamma(q)-\epsilon}\mathcal{D}_\delta^q(\mu_v)&\gtrsim_q \sum_{v'\in V}\delta^{\gamma(q)-\epsilon}\mathcal{D}_\delta^q(\mu_{v'})\nonumber\\
			&\asymp_q \delta^{\gamma(q)-\epsilon} \sum_{w\in E^*_\delta} p_w^q  \mathcal{D}^q_{\delta/\alpha_1(w)} \left(\pi_w(\mu_{t(w)})\right)\nonumber\\
			&\gtrsim_{\epsilon,q} \delta^{\gamma(q)-\epsilon} \sum_{w\in E^*_\delta} p_w^q \left( \frac{\delta}{\alpha_1(w)} \right) ^{-\tau_w(q)+{\epsilon}/{2}} \nonumber \\
			&\gtrsim_{\epsilon,q}  \sum_{w\in E^*_\delta} p_w^q \alpha_1(w)^{\tau_w(q)} \alpha_2(w)^{\gamma(q)-{\epsilon}/{2}-\tau_w(q)} \cdot \alpha_1(w)^{-{\epsilon}/{2}} \nonumber\\
			&\geq \sum_{w\in E^*_\delta} \varphi^{\gamma(q)-{\epsilon}/{2},q}(w)\nonumber\\
			&\gtrsim_{\epsilon,q} 1, \nonumber
		\end{align}
		which yields that $\gamma(q)\leq \underline{\tau}_{\mu_v}(q)$. 
		
		Therefore for each $v\in V$, $\tau_{\mu_v}(q)$ exists and equals to $ \gamma(q)$. 
	\end{proof}
	

	\section{Examples}\label{sec5}
\begin{figure}[htbp]
		\centering
		\subfigure{
			\includegraphics[width=0.3\textwidth]{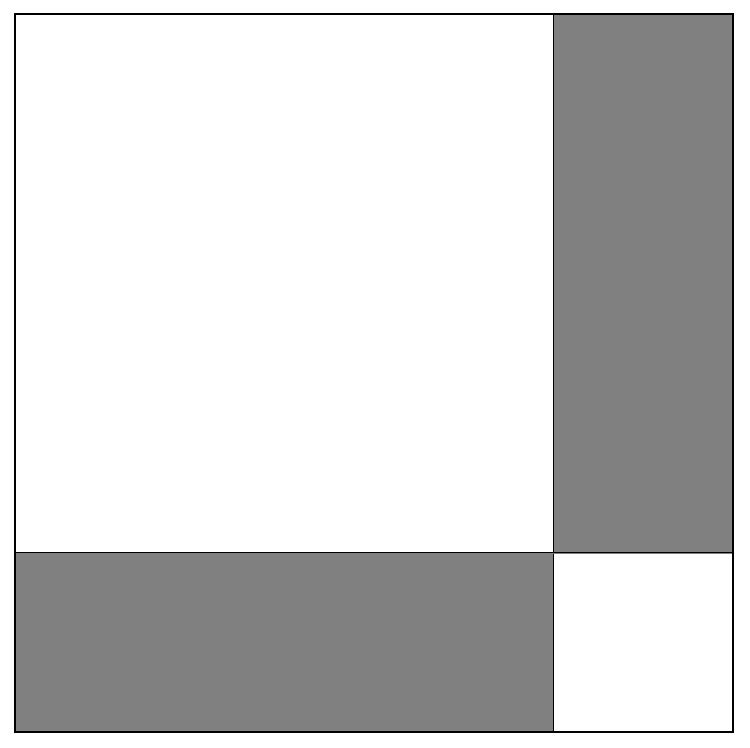}}
		\subfigure{\includegraphics[width=0.3\textwidth]{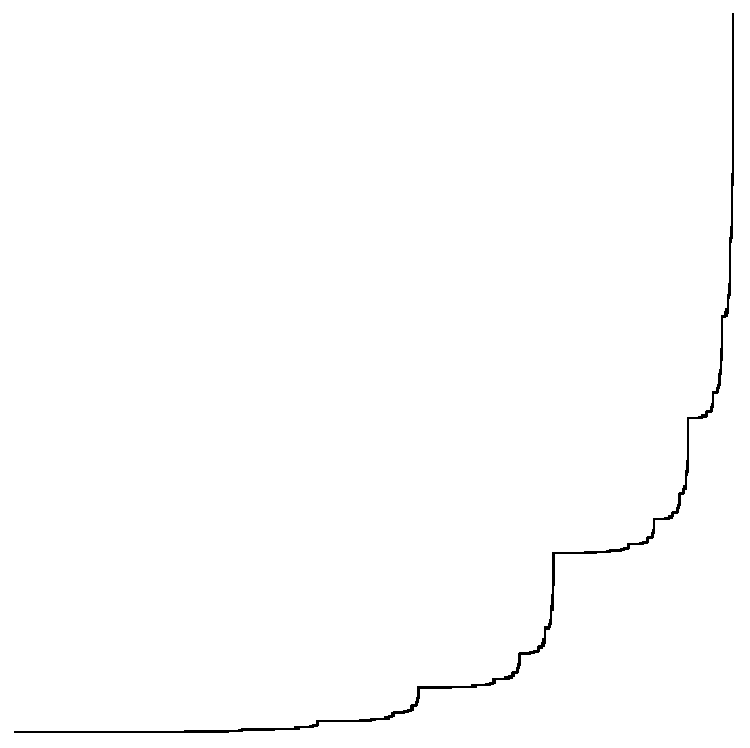}}
        \subfigure{\includegraphics[width=0.3\textwidth]{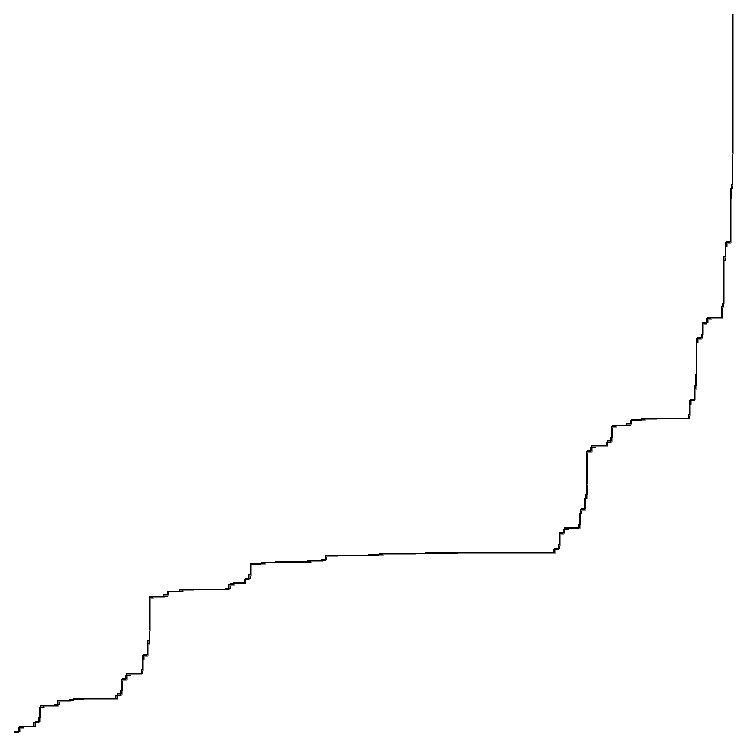}}
		\caption{Left: the shaded rectangles are images of the iterated function $\psi_1$ (resp. $\psi_1'$) and $\psi_2$ (resp. $\psi_2'$). Middle: the attractor $X$. Right: the attractor $X'$. }
		\label{f1}
	\end{figure}

 In this section, we provide two examples to illustrate our results. We only look at the  IFS case for simplicity. For some $a,b\in (0,1)$ with $a+b\leq 1$, let 
	$$
	\psi_1(\xi_1,\xi_2) = \left (
	\begin{aligned}
		&a &0\\
		&0 &b
	\end{aligned}
	\right ) 
	\left( \begin{aligned}
		&\xi_1\\
		&\xi_2
	\end{aligned}\right) \quad \text{ and } \quad
	\psi_2(\xi_1,\xi_2) = \left (
	\begin{aligned}
		&b &0\\
		&0 &a
	\end{aligned}
	\right ) 
	\left( \begin{aligned}
		&\xi_1\\
		&\xi_2
	\end{aligned}\right) +\left(
	\begin{aligned}
		&1-b\\
		&1-a
	\end{aligned} \right).
	$$ 
    Then $\{\psi_1, \psi_2\}$ becomes a planar box-like self-affine IFS. Let $X$ be its attractor. Let
     
    $$
	\psi'_1(\xi_1,\xi_2) = \left (
	\begin{aligned}
		&0 &a\\
		&b &0
	\end{aligned}
	\right ) 
	\left( \begin{aligned}
		&\xi_1\\
		&\xi_2
	\end{aligned}\right) \quad \text{ and } \quad
	\psi'_2=\psi_2,	$$
    and $X'$ be the attractor of the planar box-like self-affine IFS $\{\psi_1',\psi_2'\}$. Note that the images of $\psi_1$ and  $\psi_1'$ (resp. $\psi_2$ and  $\psi_2'$) under $[0,1]^2$ are same.
	See Figure \ref{f1} for $X,X'$ when  $a=3/4,b=1/4$.

	Let $\mu$ (resp. $\mu'$) be the self-affine measure associated with  $\{\psi_1,\psi_2\}$ (resp. $\{\psi_1',\psi_2'\}$) and  a probability vector $\mathcal{P}=(1/2,1/2)$. We compute the closed form expression for $L^q$-spectra of $\mu,\mu'$ respectively. 
 \vspace{0.2cm}

   \noindent \textit{For IFS $\{\psi_1,\psi_2\}$:}
    \vspace{0.2cm}
    
    For $q\geq 0$,  $\tau_{\mu^x}(q)=\tau_{\mu^y}(q)$ and equals to the unique solution $s(q)$ of 
	\begin{equation}\label{t10}
		\left(\frac{1}{2}\right)^q a^{s(q)} +	\left(\frac{1}{2}\right)^q b^{s(q)}=1,
	\end{equation}
	$\gamma_A(q)=\gamma_B(q)$ and equals to the unique solution $r(q)$ of 
	\begin{equation}\label{t11}
		\left( \frac{1}{2} \right)^q a^{s(q)} b^{r(q)-s(q)} + \left( \frac{1}{2} \right)^q b^{s(q)} a^{r(q)-s(q)}  =1.
	\end{equation}
	Combining \eqref{t10} and \eqref{t11}, we know that 
	$r(q)=s(q)$.
We can use either Corollary \ref{cor1} or Corollary \ref{cor4} to calculate the closed form expression for $\tau_{\mu}(q)$, $q\geq 0.$
\vspace{0.2cm}
 
\textit{Using Corollary \ref{cor1}. }
Note that $\max \{\gamma_A(q),\gamma_B(q)\} \leq \tau_{\mu^x}(q)+\tau_{\mu^y}(q)$ is equivalent to $s(q)\geq 0$. Combining this with \eqref{t10}, we know that  $\max \{\gamma_A(q),\gamma_B(q)\} \leq \tau_{\mu^x}(q)+\tau_{\mu^y}(q)$ is equivalent to $0\leq q\leq 1$. 
	By Corollary \ref{cor1}, we know that 
\begin{equation}\label{tauq}
	\tau_{\mu}(q)=\left\{  \begin{aligned}
		&s(q) 
		& &\text{ if }  0\leq q \leq  1,
		\\
		&\min \{x+y:\left( \frac{1}{2} \right)^q a^{x} b^{y} + \left( \frac{1}{2} \right)^q b^{x} a^{y}  =1, s(q)\leq x \leq 0\}  & &\text{ if }  q \geq 1.
	\end{aligned}\right.
\end{equation}
 Consider the implicit function $y(x)$ determined by    $\left( \frac{1}{2} \right)^q a^{x} b^{y(x)} + \left( \frac{1}{2} \right)^q b^{x} a^{y(x)}  =1$.
  Take $x=\frac{(q-1)\log 2}{\log ab}$, so that $y'(x)=-1$ and 
 \begin{equation}\label{t20}
 2x=x+y(x)=\min \{ x+y:\left( \frac{1}{2} \right)^q a^{x} b^{y} + \left( \frac{1}{2} \right)^q b^{x} a^{y}  =1\}.
 \end{equation}
 When $q\geq1$, noticing that $s(q)\leq 0$, we get $$2\left(\frac{1}{2}\right)^q (ab)^{x}=1=\left(\frac{1}{2}\right)^q a^{s(q)} +	\left(\frac{1}{2}\right)^q b^{s(q)}\leq 2\left(\frac{1}{2}\right)^q (ab)^{s(q)},
 $$
 which gives that $x\geq s(q)$. Clearly, also we have $x\leq 0$.
 Combining this with \eqref{tauq} and \eqref{t20}, we get
 \begin{equation}\label{t12}
	\tau_{\mu}(q)=\left\{  \begin{aligned}
		&s(q)
		& &\text{ if }  0\leq q \leq  1,
		\\
		&\frac{2(q-1)\log 2}{\log ab}  & &\text{ if }  q \geq 1.
	\end{aligned}\right.
	\end{equation}
\vspace{0.2cm}

\textit{Using Corollary \ref{cor4}.}
 Let
 $$
    H^{(q)}_{x,y}=\left( \begin{array}{cc}
				\left(\frac{1}{2}\right)^q a^{x+s(q)} b^{y-s(q)}+\left(\frac{1}{2}\right)^q b^{x+s(q)}a^{y-s(q)} \quad &  0\\ 0 \quad    & \left(\frac{1}{2}\right)^q b^{x+s(q)} a^{y-s(q)}+\left(\frac{1}{2}\right)^q a^{x+s(q)} b^{y-s(q)}\\	\end{array}\right).
    $$
 Then
 $$
 \rho(H^{(q)}_{x,y})=\left(\frac{1}{2}\right)^q a^{x+s(q)} b^{y-s(q)}+\left(\frac{1}{2}\right)^q b^{x+s(q)}a^{y-s(q)}.
 $$
 Thus by taking $x=0$ in the above equation,  $\hat{\gamma}(q)=r(q)$ by using \eqref{t11}. Note that when $q\geq 1$,
 $$
 \begin{aligned}&\min \{x+y: \left(\frac{1}{2}\right)^q a^{x+s(q)} b^{y-s(q)}+\left(\frac{1}{2}\right)^q b^{x+s(q)}a^{y-s(q)}=1,0\leq x \leq -s(q)\}
 \\
 = &\min \{ x+y:\left( \frac{1}{2} \right)^q a^{x} b^{y} + \left( \frac{1}{2} \right)^q b^{x} a^{y}  =1,s(q)\leq x \leq 0\}.
 \end{aligned}
 $$
Therefore we still get \eqref{tauq}, and so \eqref{t12} also follows by Corollary \ref{cor4}.

\begin{remark}
    The box-like self-affine IFS $(\psi_1,\psi_2)$ was considered in \cite{FFL21} which illustrated that $ \tau_{\mu}(q) <\min \{\gamma_A(q),\gamma_B(q)\}
	$ may happen. Recently, Kolossv\'{a}ry \cite[Proposition 4.4]{K23} calculated the same expression \eqref{t12} for $\tau_{\mu}(q)$ in the setting that  IFS's under consideration have grid structure.
\end{remark}

\noindent\textit{For IFS $\{ \psi'_1,\psi'_2\}$:}
\vspace{0.2cm}

    Noticing that the linear part of $\psi_1'$ is anti-diagonal, $\{\mu'^x,\mu'^y\}$ is a strongly connected self-similar graph-directed measure family, i.e.
    $$
    \mu'^x(I)=\frac{1}{2} \mu'^y(aI)+\frac{1}{2} \mu'^x(bI+1-b)
    $$
    and 
    $$
    \mu'^y(I)=\frac{1}{2} \mu'^x(bI)+\frac{1}{2} \mu'^y(aI+1-a),
    $$
    for all Borel sets $I \subseteq \mathbb{R}.$

    Let  $\beta(q)$ be the unique solution of    $$
    \rho \left( \begin{aligned}
				\left(\frac{1}{2}\right)^q b^{\beta(q)} \quad & & \left(\frac{1}{2}\right)^q a^{\beta(q)} \\ \left(\frac{1}{2}\right)^q b^{\beta(q)} \quad    & &\left(\frac{1}{2}\right)^q a^{\beta(q)}\\	\end{aligned}\right)=\left(\frac{1}{2}\right)^q b^{\beta(q)} + \left(\frac{1}{2}\right)^q a^{\beta(q)} =1.
    $$
    It follows from the result in \cite{S93} that $\tau_{\mu'^x}(q)=\tau_{\mu'^y}(q)=\beta(q).$ So $\tau_{\mu'^x}(q)=\tau_{\mu'^y}(q)=\tau_{\mu^x}(q)=\tau_{\mu^y}(q)$ by \eqref{t10} and $\beta(q)=s(q)$. Take 
    $$
    H'^{(q)}_{x,y}=\left( \begin{aligned}
				\left(\frac{1}{2}\right)^q b^{x+\beta(q)}a^{y-\beta(q)} \quad & & \left(\frac{1}{2}\right)^q a^{x+\beta(q)} b^{y-\beta(q)}\\ \left(\frac{1}{2}\right)^q b^{x+\beta(q)} a^{y-\beta(q)} \quad    & &\left(\frac{1}{2}\right)^q a^{x+\beta(q)} b^{y-\beta(q)}\\	\end{aligned}\right),
    $$
which implies that 
$$
\rho(H'^{(q)}_{x,y})=\left(\frac{1}{2}\right)^q b^{x+\beta(q)}a^{y-\beta(q)} + \left(\frac{1}{2}\right)^q a^{x+\beta(q)} b^{y-\beta(q)}=\rho(H^{(q)}_{x,y}).
$$
So $\tau_{\mu'}(q)=\tau_{\mu}(q)$ for $q\geq 0$ by Corollary \ref{cor4}.
    

 \vspace{0.2cm}
	


	\bibliographystyle{amsplain}
	
\end{document}